\documentclass[11pt]{article}

\reversemarginpar

\newcommand{\mylabel}
{\label}

\usepackage{amssymb}
\usepackage{amsthm}
\usepackage{amsmath}
\usepackage{amscd}
\usepackage{graphicx}
 
\allowdisplaybreaks
    
\newcommand{\R}{\mathbb{R}} 
\newcommand{\C}{\mathbb{C}}
\newcommand{\Sph}{\mathbb{S}}

\newcommand{\Rtn}{\R^{2n}}
\newcommand{\Cno}{\C^{n+1}}

\newcommand{\cp}{\mathbb{CP}}
\newcommand{\eps}{\varepsilon} 
\newcommand{\me}{\mathrm{e}} 
\newcommand{\mi}{\mathrm{i}}
\newcommand{\elbow}{\, \mbox{\rule[-0.1ex]{1.125ex}{0.15ex}\rule[-0.1ex]{0.15ex}{1.25ex}}\;}
\newcommand{\dif}{\mathrm{d}} 
\newcommand{\Dif}{\mathrm{D}}
\renewcommand{\Re}{\mathbf{Re}} 
\renewcommand{\Im}{\mathbf{Im}}

\newcommand{\linspanr}{\mathrm{span_{\, \R}}}

\newcommand{\un}{U(n+1)}
\newcommand{\sun}{SU(n+1)}

\newcommand{\approxsol}{\tilde{\Lambda}_{U,\zeta}}
\newcommand{\lcm}{\mathit{l.c.m.}}
\renewcommand{\gcd}{\mathit{g.c.d.}}
\newcommand{\ext}{\mathit{ext}}
\newcommand{\neck}{\mathit{neck}}

\newcommand{\annr}{\mathit{Ann}(r/2, r)}
\newcommand{\lcs}{\mathit{LCS}}
\newcommand{\clbg}{C^{l,\beta}_\gamma}
\newcommand{\cobg}{C^{0,\beta}_{\gamma-4}}
\newcommand{\cfbg}{C^{4, \beta}_\gamma}
\newcommand{\ctbg}{C^{2, \beta}_{\gamma-2}}

\newtheorem*{mainthm}{Main Theorem}
\newtheorem{thm}{Theorem}
\newtheorem{lemma}[thm]{Lemma}
\newtheorem{cor}[thm]{Corollary}
\newtheorem{prop}[thm]{Proposition}
\newtheorem*{nonumthm}{Theorem}

\theoremstyle{definition}
\newtheorem{defn}[thm]{Definition}

\newtheoremstyle{rmk}{5pt}{5pt}{}{}{\scshape}{:}{.5em}{}
\theoremstyle{rmk}
\newtheorem*{rmk}{Remark}

\begin{document}

\title{Equivariant Gluing Constructions of Contact Stationary \\ Legendrian Submanifolds in $\Sph^{2n+1}$}

\author{Adrian Butscher \\ University of Toronto at Scarborough \\ email: \ttfamily butscher@utsc.utoronto.ca}

\maketitle

\begin{abstract}
A contact-stationary Legendrian submanifold of $\Sph^{2n+1}$ is a Legendrian submanifold whose volume is stationary under contact deformations.  The simplest contact-stationary Legendrian submanifold (actually minimal Legendrian) is the real, equatorial $n$-sphere $S_0$.  This paper develops a method for constructing contact-stationary (but not minimal) Legendrian submanifolds of $\Sph^{2n+1}$ by gluing together configurations of sufficiently many $\un$-rotated copies of $S_0$.  Two examples of the construction, corresponding to finite cyclic subgroups of $\un$ are given.  The resulting submanifolds are very symmetric; are geometrically akin to a `necklace' of copies of $S_0$ attached to each other by narrow necks and winding a large number of times around $\Sph^{2n+1}$ before closing up on themselves; and are topologically equivalent to $\Sph^1 \times \Sph^{n-1}$. 
\end{abstract}

\tableofcontents

\renewcommand{\baselinestretch}{1.25}
\normalsize


\section{Introduction and Statement of Results}

\subsection{Background}

\paragraph{Lagrangian variational problems.}  A minimal submanifold $L$ in a Riemannian manifold $M$ satisfies a classical variational problem, namely that the volume of $L$ is stationary amongst all deformations of $L$.  When $M$ is K\"ahler-Einstein (possessing a Riemannian Einstein metric and a compatible symplectic form), then it is has been known for some time that minimal \emph{and} Lagrangian submanifolds of $M$ possess a rich mathematical structure.  Moreover, it is possible to pose two very natural \emph{restricted} variational problems in the class of Lagrangian submanifolds of $M$ whose critical points are also mathematically quite interesting, and are related in a number of ways to minimal and Lagrangian submanifolds.  First, one can demand that the volume of a Lagrangian submanifold $L$ is stationary with respect to all variations of $L$ which preserve the Lagrangian condition.  A natural sub-class of variations preserving the Lagrangian condition is the set of \emph{Hamiltonian transformations} generated by all real-valued functions of $M$.  Therefore one can also demand that the volume of $L$ is stationary with respect to this sub-class of variations.  In the former case, $L$ is said to be \emph{Lagrangian stationary}; whereas in the latter case, $L$ is said to be \emph{Hamiltonian stationary}.  

The Euler-Lagrange equations for Hamiltonian stationary submanifolds in a K\"ahler-Einstein submanifold are quite simple to write down.  Let $H_L$ be the mean curvature vector of a Lagrangian submanifold $L$.  The K\"ahler-Einstein condition implies that the one-form $ H_L \elbow \omega \big|_L$ satisfies $\dif \left(  H_L \elbow \omega \big|_L \right) = 0$ where $\omega$ is the symplectic form.  Then $L$ is Hamiltonian stationary if in addition $\nabla \cdot  \left(  H_L \elbow \omega \big|_L \right) = 0$, where $\nabla \cdot$ is the divergence operator of the metric of $L$.  Therefore $L$ is Hamiltonian stationary if and only if the one-form $ H_L \elbow \omega \big|_L $ is harmonic.

The stationarity of Lagrangian submanifolds of a K\"ahler-Einstein manifold $M$ under Lagrangian and Hamiltonian deformations has been studied by several authors, most notably Oh \cite{oh2,oh1},  Helein and Romon \cite{heleinromon2,heleinromon3,heleinromon1}, Schoen and Wolfson \cite{sw1,sw2}.  Oh initially posed the Lagrangian stationary and Hamiltonian stationary variational problems and derived the Euler-Lagrange equations above as well as a second variation formula.  H\'elein and Romon show that when the $M$ is a 2-complex-dimensional Hermitian symmetric space, this mean curvature condition can be reformulated as an infinite dimensional integrable system whose solutions possess a Weierstra\ss-type representation similar to the Dorfmeister-Pedit-Wu representation for harmonic surfaces in a real symmetric space \cite{dpw}.  Moreover, they find all Hamiltonian stationary --- and non-minimal --- doubly periodic immersions of $\R^2$ into $\cp^2$ using this representation.  They did not address the question of when these doubly periodic immersions close up to form immersions or embeddings of tori into $\cp^2$, though the existence of closed configurations is certainly predicted by integrable systems methods.  Finally, Schoen and Wolfson initiated the study of the Lagrangian variational problem from the geometric analysis point of view, treating it as a method for constructing minimal and Lagrangian submanifolds as limits of volume-minimizing sequences of Lagrangian submanifolds.  The rationale behind this approach is their observation that a Lagrangian stationary submanifold of $M$ that is smooth must necessarily be minimal (because in this case the mean curvature vector field of $L$ is itself the infinitesimal generator of a Lagrangian variation, allowing the volume to be further decreased unless the mean curvature vanishes).

\paragraph{Singularities and Legendrian variational problems.}

From the work of Schoen and Wolfson it has emerged that questions concerning the regularity of limits of volume-minimizing sequences of Lagrangian submanifolds are very delicate and that the structure of the singularities of Lagrangian stationary submanifolds can be very complicated.  Indeed, general structure theorems for these singularities (such as classifications) are currently well out of reach and the most fruitful approach at the moment is the constructions of classes of examples of Lagrangian stationary, Hamiltonian stationary, and minimal Lagrangian submanifolds possessing singularities of various types.   

A class of singularities that is particularly amenable to study is when the stationary Lagrangian submanifolds possess an \emph{isolated conical singularity}.  Such a singularity occurs when the  stationary Lagrangian submanifold $L$ is singular at an isolated point where the tangent cone (in the sense of geometric measure theory) is modeled on an actual cone, i.e.~the tangent cone $C$ is a complete, homothetically invariant, submanifold of Euclidean space having multiplicity one.  In this case, $C$ is determined by the submanifold $\Lambda$ formed by intersecting $C$ with the unit sphere in Euclidean space, called the \emph{link} of the cone.  The Lagrangian stationary, Hamiltonian stationary, and minimal Lagrangian conditions then translate into conditions satisfied by $\Lambda$.

To express these conditions, recall the geometric structures of Euclidean space $\C^{n+1}$, which is the simplest example of a Calabi-Yau manifold where the metric and symplectic form are the standard ones and the holomorphic volume form is $\Omega = \dif z^1 \wedge \cdots \wedge \dif z^{n+1}$.  If $L$ is a Lagrangian submanifold then it can be shown that $\Omega(E_1, \ldots E_{n+1})$ is a complex number of unit length for every orthonormal basis $E_1, \ldots, E_{n+1}$ for $T_pL$ and $p \in L$.  Furthermore, the value of this complex number is independent of the choice of basis so that the prescription $p \mapsto \exp(\mi \Theta_L(p) )$ defines a possibly multi-valued function $\Theta_L$ on $L$ which is called the \emph{Lagrangian angle function} \cite{hl1}.  It can also be shown that this function satisfies $ H_L \elbow \omega \big|_L  = \dif \Theta_L$.  Hence if $L$ is minimal and Lagrangian then $\Theta_L = \mathit{constant}$ and if $L$ is a Hamiltonian stationary submanifold, then $\Delta \Theta_L = 0$ locally and $\dif \Theta_L$ defines a harmonic one-form globally.   The unit sphere $\Sph^{2n+1}$ in $\C^{n+1}$ inherits the induced metric as well as a \emph{contact structure} --- that is, $\Sph^{2n+1}$ possesses a fully non-integrable  hyperplane field $\Xi$ called the contact distribution.  The hyperplane $\Xi_p$ is equal to the orthogonal complement of the tangent vector of the Hopf fibration of $\Sph^{2n+1}$ obtained by considering all curves of the form $\me^{\mi t} z_0$ where $t \in [0,2 \pi)$ and $z_0 \in \Sph^{2n+1}$.  Finally if $C$ is a Lagrangian cone in $\C^{n+1}$ then its link $\Lambda := C \cap \Sph^{2n+1}$ is a \emph{Legendrian} submanifold, i.e.~an integrable submanifold of $\Xi$ of the largest possible dimension. 

The fundamental observations are that if a Lagrangian cone $C$ is minimal then the Legendrian link  $\Lambda$ is minimal; and if $C$ is Lagrangian or Hamiltonian stationary, then $\Lambda$ is stationary under all variations of $\Lambda$ through Legendrian submanifolds.  It turns out that this is equivalent to $\Lambda$ being stationary under all contact structure-preserving variations (a.k.a~\emph{contactomorphisms}) and $\Lambda$ is said to be \emph{contact-stationary}.  Finally, the Lagrangian angle function of $C$ is also determined by its values on $\Lambda$, so that one can define a \emph{Legendrian angle function} $\Theta_\Lambda$ which is constant when $\Lambda$ is minimal Legendrian and defines a harmonic one-form when $\Lambda$ is contact-stationary.  The relation between $\Theta_\Lambda$ and the mean curvature $H_\Lambda$ becomes $\dif \Theta_\Lambda = H_\Lambda \elbow \omega \big|_{\Sph^{2n+1}}$.  Thus contact-stationary Legendrian submanifolds satisfy the equation $\nabla \cdot \left( H_\Lambda \elbow \omega \big|_{\Sph^{2n+1}} \right) = 0$.

\paragraph{Minimal and contact-stationary Legendrian submanifolds of {\boldmath $\Sph^{2n+1}$}.}

The study of minimal Legendrian and contact-stationary Legendrian submanifolds of the sphere $\Sph^{2n+1}$ is relatively recent endeavour, but a certain number of results exist.  The simplest case is in dimension $n=1$.  It is clear that the minimal Legendrian submanifolds are contact curves that are also great circles.  The contact-stationary submanifolds are non-trivial, however.  They are the so-called $(p,q)$-curves discovered by Schoen and Wolfson in \cite{sw2}, where $p$ and $q$ are relatively prime integers parameterizing these curves.  

In higher dimensions, the situation changes entirely.  In dimension $n=2$, Yau showed in \cite{yauleg1,yauleg2} that genus zero minimal Legendrian submanifolds of $\Sph^5$ are trivial; i.e.~they are equatorial 2-spheres.  Moreover, since 2-spheres possess no non-trivial harmonic one-forms, there are no contact-stationary Legendrian submanifolds of genus zero that are not minimal.  In higher genus, however, the situation is quite different.  In genus one, a huge abundance of minimal Legendrian tori in $\Sph^5$ has recently been discovered that can be studied using techniques of integrable systems theory \cite{carberrymcintosh,haskins2,mcintosh} and work has been done to understand the so-called `geometric complexity' of these torus cones, in the sense of Haskins \cite{haskins2}.  This latter task is important since Joyce has conjectured that the singular minimal Lagrangian submanifolds of $\C^3$ with isolated conical singularities whose cone links have the least geometric complexity are those which occur `most often' in the boundary of the moduli space of all minimal Lagrangian submanifolds \cite{j5}.  There are plenty of contact-stationary, doubly periodic Legendrian immersions of $\R^2$ into $\Sph^5$: all are lifts of H\'elein and Romon's examples.  But the question of which of these close up to form tori is not known explicitly, except for the \emph{homogeneous} tori classified by H\'elein and Romon.  These are the product spheres $\Sph^1(r_1) \times \Sph^1(r_2) \times \Sph^1 (r_3) \subset \C^3$, with $r_1^2 + r_2^2 + r_3^2 = 1$ so that they are submanifolds of $\Sph^5$, and then taken modulo the Hopf map.  In higher genus, Wang has constructed minimal Legendrian submanifolds of various genera in $\Sph^5$ using a reflection technique \cite{sunghowang}, but these are not everywhere smooth.   The question of whether \emph{smooth} and \emph{embedded} higher-genus minimal Legendrian submanifolds exist in $\Sph^5$ has only been settled very recently by Haskins and Kapouleas \cite{haskinskapouleas}.  These authors have constructed odd-genus minimal Legendrian submanifolds by fusing together multiple copies of Haskins' $U(1)$-invariant minimal Legendrian tori \cite{haskins1} in a manner analogous to Kapouleas' fusion of Wente tori \cite{kapouleas5}.  The question of the existence of higher-genus contact-stationary Legendrian submanifolds of $\Sph^5$ is open, though in light of Haskins' and Kapouleas' result, the answer is most certainly that these do exist.

In dimensions $n \geq 3$ very little is known beyond the standard examples of the equatorial $n$-spheres and the Legendrian Clifford tori.   It seems that none of the techniques available in dimension $n=2$ carry over to the $n\geq 3$ case: there is no integrable systems framework for studying the minimal and contact-stationary equations in these dimensions; and the Haskins-Kapouleas construction can not be generalized because Haskins' tori are unique to dimension $n=2$.  A general construction which does exist is by Castro, Li and Urbano \cite{castroliurbano}.  These authors have developed a `Legendrian warped product' which can be used to combine two contact-stationary Legendrian submanifolds contained in lower-dimensional spheres of dimension $2 p+1$ and $2 q+1$ to form a contact-stationary Legendrian submanifold in the sphere of dimension $2(p+q)+3$.

The purpose of this paper is to develop a method of constructing contact-stationary Legendrian submanifolds of $\Sph^{2n+1}$ using a gluing technique starting from the simplest building blocks.  That is, an approximately contact-stationary Legendrian submanifold will be constructed by forming connected sums of simple, minimal Legendrian building blocks; and then this construction will be perturbed (by solving the non-linear PDE satisfied by contact-stationary Legendrian submanifolds) to yield an exactly contact-stationary Legendrian submanifold.  The method is such that the Legendrian angle function $\Theta_\Lambda$ is forced to acquire \emph{periods} and that the harmonic one-form $\dif \Theta_\Lambda$ represents a non-trivial element of the first cohomology of $\Lambda$.  The exactly special Legendrian submanifold that results from the perturbation process can thus never be minimal.

\subsection{The Main Theorem}

\paragraph{Preliminaries.}  

In order to state the Main Theorem to be proved in this paper, it is necessary to first introduce some important terminology.  The first concept that must be put into words is a special way in which two Legendrian submanifolds of the sphere can intersect each other.

\begin{defn}
	Let $\Lambda, \Lambda'$ be two Legendrian submanifolds of $\Sph^{2n+1}$. Then $\Lambda$ and $\Lambda'$ are \emph{contact-transverse} over a Hopf fiber $\mathcal F$ if the following two conditions are satisfied.
	\begin{enumerate}
		\item There exists a point $p \in \Lambda \cap \mathcal F$ and a complex number $\me^{\mi \alpha}$ so that $\me^{\mi \alpha} p \in \Lambda' \cap \mathcal F$.
		
		\item  $T_p \Lambda \oplus T_p (\me^{-\mi \alpha} \Lambda') = \Xi_p$.  
	\end{enumerate}
	The submanifolds $\Lambda$ and $\Lambda'$ are said to have \emph{contact-transverse intersection} at $p$ if $\alpha = 0$ in the definition above, so that $p \in \Lambda \cap \Lambda'$.
\end{defn}
	
\begin{defn}
	If $\Lambda$ and $\Lambda'$ are contact-transverse over a Hopf fiber $\mathcal F$, then the number $\alpha$ is called the \emph{Hopf separation} between $\Lambda$ and $\Lambda'$ over $\mathcal F$.
\end{defn}	

A construction needed in this paper is the \emph{Legendrian connected sum} of two Legendrian submanifolds $\Lambda$ and $\Lambda'$ with contact-transverse intersection at a point $p$.  The details of this construction will be given in Section \ref{subsec:legconnect}.  In the mean time, it suffices to state the condition on the tangent spaces $T_p\Lambda$ and $T_p\Lambda'$ under which the construction is possible.  The condition can be phrased in terms of an \emph{angle criterion} satisfied by the \emph{characteristic angles} of the direct sum of these two tangent spaces. When this angle criterion is satisfied, then there will exist a Legendrian neck that can be used to connect $\Lambda$ to $\me^{\mi \alpha} \Lambda'$ in a neighbourhood of $p$.

The characteristic angles and the angle criterion can best be described as follows. First, realize that $\Xi_p$ is isomorphic to the symplectic vector space $\C^n$, and $T_p \Lambda$ and $T_p \Lambda'$ are isomorphic to Lagrangian $n$-planes in $\C^n$.  It can be shown that for any pair of transversely intersecting Lagrangian $n$-planes $\Pi_{1}$ and $\Pi_{2}$ in $\C^n$, there exists a Hermitian orthonormal  basis $E_{1}, \ldots, E_{n}$ of $\C^n$ and a unique set of angles $\theta_{k} \in (0, \pi/2]$ for $k=1, \ldots , n-1$ and $\theta_n \in [\pi/2, \pi)$ such that
\begin{gather*}
    	\Pi_{1} = \linspanr \big\{ E_{1}, \ldots, E_{n} \big\} \\
    	\Pi_{2} = \linspanr \big\{ \me^{\mi \theta_1} E_{1} , \ldots, \me^{\mi \theta_n} E_{n} \big\} \, .
\end{gather*}
The existence of this basis follows from standard complex linear algebra and can be found in \cite{harvey}. 

\begin{defn}
	The \emph{angle criterion} for the pair of $n$-planes $\Pi_{1}$ and $\Pi_{2}$ is that their characteristic angles $\theta_1, \ldots, \theta_n$ satisfy $\theta_{1} + \cdots + \theta_{n} = \pi$.
\end{defn}

\noindent Note that the angle criterion need not hold for a general pair of intersecting special Lagrangian planes, though it is \emph{always} satisfied in dimensions $n=2,3$ for numerical reasons.

\paragraph{The concept behind  the construction.}  The real equatorial $n$-sphere $S_0 = \Sph^{2n+1} \cap ( \R^{n+1} \times \{0\} )$ is undoubtedly the simplest minimal and Legendrian submanifold of $\Sph^{2n+1}$ that one can conceive.  The idea behind the construction of this paper is to use $S_0$ and isometric copies of itself as building blocks that can be connected together at various points of contact-transverse intersection by means of small Legendrian necks to form a long, closed chain of $n$-spheres having the topology of the closed cylinder $\Sph^1 \times \Sph^{n-1}$ and winding around $\Sph^{2n+1}$.  The construction of the Legendrian necks will also be carried out herein and is based on the Lawlor neck \cite{lawlor2,lawlor1} which has already been used in connected sum constructions involving intersecting minimal Lagrangian submanifolds \cite{me1,leed,leey}.  

Since each $n$-sphere is minimal Legendrian, the first hope is that the long chain will be almost \emph{minimal} Legendrian and can be perturbed using a small contactomorphism into an exactly minimal Legendrian submanifold.  Indeed, a certain limiting subset of the Haskins $U(1)$-invariant minimal Legendrian tori does have a very similar geometric characterization as the configuration proposed above.  However, the deformation to a minimal Legendrian submanifold should not be possible for the following reason.   It turns out that the Legendrian necks have a non-negligible \emph{height} in the Hopf fiber direction.  Thus to use such a neck for connecting two spheres with the least error, it is actually first necessary to create just the right amount of Hopf separation between them.  The problem is that the Legendrian angle function of a translate of $S_0$ by $\me^{\mi \alpha}$ along the Hopf fibers is $\alpha$.  Thus if a large number of $n$-spheres is attached together and the angle function increases by $\alpha$ from one sphere to the next, then the angle function can not remain small.  Moreover, if the long chain of $n$-spheres winds completely around $\Sph^{2n+1}$ and closes up, then the Legendrian angle function must acquire a \emph{period} and cease to be single-valued.  One would therefore not expect to be able to deform this configuration into a Legendrian submanifold with angle zero.  

The solution is to give up attempting to construct a minimal Legendrian chain of $n$-spheres and instead attempt to construct a contact-stationary Legendrian chain of $n$-spheres.  The angle function of such an object must satisfy $\nabla \cdot \dif \Theta_\Lambda = 0$  so a slowly varying angle function acquiring a period around the chain can be incorporated into the construction.  And since there are cohomologically non-trivial harmonic one-forms on a closed cylinder, one might expect to find solutions.

\paragraph{The construction of the approximate solution.}  The Main Theorem of this paper is that the concept posed in the previous paragraphs can be realized provided the building block $n$-spheres and Legendrian necks can be assembled in a sufficiently symmetric manner.  The first step is to construct an \emph{approximate solution} of the problem via the following three steps.

\begin{enumerate}

	\item Choose $U \in \sun$ generating a cyclic subgroup of $\sun$ of integer order $N$.  Suppose further that $S_0 \cap U(S_0)$ has contact-transverse intersection at the points $\pm p \in S_0$ in such a way that the tangent spaces at the intersection points satisfy the angle criterion. 
		
	\item Let $\zeta := \me^{2 \pi \mi a / k N}$ for some integer $a$ and large integer $k$ satisfying $\gcd(a, kN) = 1$.   
		
	\item Glue each $(\zeta U)^s (S_0)$ to its nearest neighbours by inserting Legendrian necks.
	
\end{enumerate}

The resulting submanifold will be denoted by $\approxsol$ and is a smooth Legendrian, and if $U$ and $\zeta$ are chosen properly, also embedded.  Each $(\zeta U)^s (S_0)$ is connected to exactly two of its neighbours in such a way that what is produced is a chain of $kN$ rotated copies of $S_0$ winding $a$ times around $\Sph^{2n+1}$ and eventually closes up.  The Legendrian angle function on the part of $(\zeta U)^s(S_0)$ away from the necks is almost constant equal to $2 \pi s a / kN$.  Hence that angle function increases only in the neck regions from one $n$-sphere to the next; and increases by $2 \pi a$ for every loop around the chain. 

\paragraph{Deforming the approximate solution and the Main Theorem.}  The means for deforming $\approxsol$ into a contact-stationary Legendrian submanifold can be sketched as follows.  It will be shown later how to generate small contactomorphisms starting with functions defined on $\approxsol$.  In the mean time, let $\phi_f: \Sph^{2n+1} \rightarrow \Sph^{2n+1}$ denote the contactomorphism generated by the function $f : \approxsol \rightarrow \R$.  The function $f$ will be selected to satisfy the equation $\nabla \cdot \left( H_{\phi_f(\approxsol) } \elbow \dif \alpha \right)= 0$, in which case $\phi_f(\approxsol)$ is a contact-stationary Legendrian submanifold.   This equation turns out to be a fourth order, nonlinear, elliptic partial differential equation on $\approxsol$.

The equation $\nabla \cdot  \left( H_{ \phi_f(\approxsol) } \elbow \dif \alpha \right) = 0$ will be solved perturbatively near zero and thus it is necessary to study the linearization of this equation at zero, hereinafter denoted $\mathcal L_{U, \zeta}$.  The ability to solve the equations hinges on being able to find a right inverse, bounded above by a constant independent of $k$, for $\mathcal L_{U, \zeta}$.  However, it is a manifestation of a general phenomenon in geometric singular perturbation problems that $\mathcal L_{U, \zeta}$ can have small eigenvalues tending to zero as $k \rightarrow \infty$.  The associated approximate co-kernel constitutes an obstruction to solvability.  The eigenfunctions with small eigenvalue, called \emph{Jacobi fields},  have a geometric origin and can be understood in an approximate sense: as generators of continuous transformation of $\approxsol$ that act by $\un$-rotation of exactly one of the constituent $n$-spheres of $\approxsol$  while leaving all the others fixed.  (See \cite[Appendix B]{kapouleas7} for a proof of this fact.)   

One way to avoid the obstruction described above is to exploit symmetry.  That is, if $\approxsol$ possesses a group of isometries, then one can deform it equivariantly; i.e.~by contactomorphisms that are invariant under the isometries.  If the group of isometries is so large that there are no Jacobi fields invariant under all isometries of $\approxsol$ at once, then equivariant deformation has the effect of eliminating the Jacobi fields.    This is precisely the setting of this paper, where the Main Theorem that will be proved is the following.

\begin{mainthm}
	Let $\approxsol$ be the approximately contact-stationary Legendrian submanifold of $\Sph^{2n+1}$ constructed above.  If $U$ is chosen appropriately, then there exists $k_0$ so that if $k>k_0$, then there is a small contact deformation of $\approxsol$ into a contact-stationary Legendrian submanifold $\Lambda_{U, \zeta}$ possessing the same symmetries as $\approxsol$.  Finally, $\Lambda_{U, \zeta}$ is embedded whenever $\approxsol$ is.
\end{mainthm}

\noindent Two  examples of $\sun$-rotations to which this theorem applies are given in Section \ref{sec:examples}.

\paragraph{Acknowledgements.}   I would like to thank Rick Schoen for steering me towards problems involving minimal and contact-stationary Legendrian submanifolds at the IPAM conference on Lagrangian submanifolds in the spring of 2003.  During the course of this research, I have benefitted from fruitful discussions with many people, including: Justin Corvino, Hansj\"org Geiges, Rafe Mazzeo, Robert McCann, Tommaso Pacini, Dan Pollack, Jesse Ratzkin, Josh Sabloff, Vin da Silva and Jon Wolfson.  I would also like to acknowledge the support of the Fields Institute during the early stages of this project, and to acknowledge my colleagues at the University of Toronto for making this possible.  Finally, I owe an enormous debt of thanks to Frank Pacard for extremely valuable, insightful and stimulating discussions on this topic.

\section{Geometric Preliminaries}
\label{sec:geomprelim}

\subsection{Riemannian Contact Manifolds}

\paragraph{Contact structures.} 
A contact structure on an odd-dimensional manifold $M^{2n+1}$ is a hyperplane field $\Xi$ on $M$, which is fully non-integrable in the sense of the Frobenius theorem.  A 1-form $\alpha$ on $M$ whose kernel at $p \in M$ is the contact hyperplane $\Xi_p \subseteq T_pM$ is called a \emph{contact form} and full non-integrability can be expressed via the requirement that $\alpha \wedge ( \dif \alpha)^{\wedge n} \neq 0$.   In turn, this relationship ensures that $\dif \alpha$ is a non-degenerate, skew-symmetric 2-form on each contact hyperplane   --- in short, a symplectic form.  The local structure theorem in contact geometry is \emph{Gray's Theorem}: it asserts that every contact manifold is locally diffeomorphic to $\R^{2n}  \times \R$ with coordinates $(x,y,t)$ and a `standard' contact form such as $$\alpha := \dif t + \frac{1}{2} \sum_{k=1}^{n} \left( x^k \dif y^k - y^k \dif x^k \right) \qquad \mbox{or} \qquad \alpha' := \dif t - \sum_{k=1}^n y^k \dif x^k \, .$$   Note that $\alpha \big(\frac{\partial}{\partial t} \big) = \alpha' \big(\frac{\partial}{\partial t} \big) = 1$ and $\dif \alpha = \dif \alpha' = \sum_{k=1}^n \dif x^k \wedge \dif y^k$ is the standard symplectic form of $\R^{2n}$.  These are the defining features of a `standard' contact structure.  Gray's Theorem is the analogue of the Darboux Theorem of symplectic geometry.  

A contact manifold $M$ possesses a canonical non-vanishing vector field called the \emph{Reeb vector field} defined by the requirement that $R \elbow \dif \alpha = 0$ and $\alpha(R) = 1$.  For example, the standard contact structures of $\R^{2n+1}$ have Reeb vector field  $R :=\frac{\partial}{\partial t}$.  Denote by $[p]$ the integral curve of the Reeb vector field passing through $p$.  The set of all integral curves of the Reeb field is called the \emph{characteristic foliation}.  When $[p]$ is a smooth, 1-dimensional embedded submanifold for every $p \in M$, then the space of fibers is a smooth manifold $\hat{M} = \{ [p] \, : \, p \in M\}$.  The canonical projection $\pi : M \rightarrow \hat{M}$ is a submersion such that $\mathrm{Ker} \big( (\pi_\ast)_p \big) = \mathrm{span} (R_p)$ and $\pi_\ast ( \Xi_p ) = T_{[p]} \hat M$ at any $p \in M$.  Furthermore, one has a canonical lifting process $\lambda_p : T_{[p]} \hat{M} \rightarrow T_pM$ for any $p \in \pi^{-1}([p])$ by defining $\lambda_p(\hat{X}_{[p]})$ to be the unique vector in $\Xi_p$ satisfying $\pi_\ast ( \lambda_p \big(\hat{X}_{[p]}) \big) = \hat{X}_ {[p]}$.  As a result, whenever a structure on $M$ is given that is equivariant with respect to the Reeb fibration, one obtains a similar structure on $\hat{M}$.  For instance, $\hat{M}$ is a symplectic manifold because one can show that $\dif \alpha$ is equivariant and the prescription $\hat{\omega}(\hat{X}, \hat{Y}) := \dif \alpha \big( \lambda_p(\hat{X}), \lambda_p(\hat{Y}) \big)$ for any pair of vector $\hat X, \hat Y \in T_{[p]} \hat M$ yields a well-defined symplectic form satisfying $\pi^\ast \hat{\omega} = \dif \alpha$. 

The diffeomorphisms of $M$ that preserve the contact structure are called \emph{contactomorphisms}.  A contactomorphism $\phi$ must preserve the kernel of $\alpha$, so that  $\phi^\ast \alpha = \me^F \alpha$ for some function $F: M \rightarrow \R$.  Thus if a one-parameter group of contactomorphisms of $M$ is generated by a vector field $X$, then Lie differentiation shows that $X$ must satisfy $\dif ( X  \elbow \alpha ) + X \elbow \dif \alpha = \dot{F} \alpha$.   One deduces that there is a function $u : M \rightarrow \R$ so that $X := X_u$ satisfies 
\begin{equation}
	\label{eqn:contactdef}
	\alpha(X_u) = u \qquad \mbox{and} \qquad X_u \elbow \dif \alpha \Big|_{\Xi} = -\dif u \, .
\end{equation}
The vector field $X_u$ is the \emph{contact Hamiltonian vector field} generated by $u$.  Note that if $u$ is Reeb-equivariant, then $u$ descends to a function $\hat u : M \rightarrow \R$ and $\pi_\ast ( X_u ) $ is the Hamiltonian vector field of $\hat M$ generated by the function $\hat u$.  Conversely, Hamiltonian vector fields of $\hat M$ can be lifted to equivariant contact-Hamiltonian vector fields of $M$.  Moreover, if $\phi^t$ is the one-parameter group of contactomorphisms generated by $u$ and  $\hat \phi^t$ is the one-parameter group of symplectomorphisms generated by $\hat u$, then it is straightforward to check that $\phi^t$ is Reeb-equivariant and $\pi \circ \phi^t (p) = \hat \phi^t([p])$ for all $p \in [p]$.  When the circumstance described here holds, one says that $\phi^t$ \emph{covers} $\hat \phi^t$.

Suppose that $\Lambda \subseteq M$ is an integral submanifold of $\Xi$ (i.e.~$T_p\Lambda \leq \Xi_p$ for every $p \in \Lambda$), then the full non-integrability condition implies that $\dif \alpha \big|_\Lambda = 0$.  But since $\dif \alpha$ is a symplectic form on each contact hyperplane, the tangent spaces of $\Lambda$ are isotropic subspaces of the contact hyperplanes and hence their dimension can be no larger that $n$.  An integral submanifold of $\Xi$ with this maximal dimension is called a \emph{Legendrian} submanifold.  It is easy to see that any Legendrian submanifold of $M$ projects under $\pi$ to a Lagrangian submanifold of $\hat{M}$.  Conversely, if $L$ is any Lagrangian submanifold of $\hat{M}$, then each of its tangent spaces can be lifted to a family of subspaces of $\Xi$.  If $\hat{X}$ and $\hat{Y}$ are two tangent vectors of $L$, then the calculation $\alpha ( [ \lambda(\hat{X}), \lambda(\hat{Y}) ] \big)  = - \dif \alpha \big( \lambda(\hat{X}), \lambda(\hat{Y}) \big) = \hat{\omega}(\hat{X}, \hat{Y}) = 0$ shows that the sub-distribution of $\Xi$ formed in this way is integrable in the sense of the Frobenius theorem.  The integral submanifolds of this sub-distribution are a family of Legendrian submanifolds of $M$ projecting onto $L$ under $\pi$.  

There is a great abundance of Legendrian submanifolds in $\R^{2n+1}$ with the contact structure $\alpha$ which are graphs of functions. The following proposition describes these fully.

\begin{prop}
	\label{prop:leg}
	 Let $\Pi = \R^n \times \{0 \} \times \{ 0\} \subseteq \Rtn \times \R$ with Lagrangian projection $[\Pi] = \R^n \times \{ 0 \}  \subseteq \Rtn$.	 
	\begin{enumerate}
	
		\item Suppose $L$ is a Lagrangian submanifold of $\Rtn$ that is graphical over $[\Pi]$.  Then there exists  a function $f : \R^n \rightarrow \R$ so that $L = L_f$ where $L_f := \{ (x, \nabla f(x) ) : x \in \R^n \}$.  
		
		\item Furthermore, $L$ can be lifted to the family of graphical Legendrian submanifolds $\Lambda_{f,c}$ in $\Rtn \times \R$ where  $\Lambda_{f,c} := \{ (  x, \nabla f(x) , 2 f(x) + \sum_{k=1}^n x^k \frac{\partial f}{\partial x^k} + c ) : x \in \R^n \}$. 
		
		\item If $\Lambda$ is a Legendrian submanifold in $\Rtn \times \R$ that is graphical over $\Pi$ then there exists $c \in \R$ and $f : \R^n \rightarrow \R$ so that $\Lambda = \Lambda_{f,c}$. 
	\end{enumerate}
\end{prop}

\noindent A consequence of this result is that Gray's Theorem now implies that Legendrian submanifolds of an arbitrary contact manifold $M$ come in great abundance, at least locally, since any sufficiently small neighbourhood of $M$ can be endowed with coordinates in which the contact structure has the standard forms above.  

A final result about Legendrian submanifolds is called  the \emph{Legendrian Neighbourhood Theorem} and it asserts that there is a contactomorphism between a tubular neighbourhood of any Legendrian submanifold $\Lambda \subseteq M$ and a tubular neighbourhood of the zero section in $T^\ast \Lambda \times \R$ endowed with its canonical contact form $\dif t - \sigma$, where $\sigma$ is the canonical one-form of $T^\ast \Lambda$.  Thus \emph{any} nearby Legendrian submanifold $\Lambda'$ that is graphical over $\Lambda$, when pulled back under this contactomorphism, is the one-jet of a function on  $\Lambda$.  One can conclude from this the important fact that \emph{Legendrian deformations} of $\Lambda$ (i.e.~deformations of $\Lambda$ through Legendrian submanifolds via a one-parameter family of contactomorphisms) are generated by functions on $\Lambda$.

\paragraph{Sasakian geometry.} 
An odd-dimensional manifold $M$ equipped with both a contact structure and a Riemannian metric, compatible with each other  in the nicest way, is known as a \emph{Sasakian structure} on $M$.  A manifold with a Sasakian structure is the most general possible arena in which to envisage the equations studied in this paper.   There are several ways of defining a Sasakian structure, the most germane of which is essentially taken from \cite{boyersurvey}, with a few modifications and additions to suit the needs of this paper. 

\begin{defn}
	Suppose $M$ is an odd-dimensional manifold that carries a Riemannian metric $g$ and a Killing field $R$ of constant length 2.  Define the endomorphism $J \in \mathit{End}(TM)$ by $J(X) = \frac{1}{2} \nabla_X R$ (so that $J(R) = 0$).  Then $(g, R, J)$ is a \emph{Sasakian structure} for $M$ if the following condition is met:
	$$(\nabla_X J) (Y)= \frac{1}{2} \big[  g(R, Y) X - g(X, Y) R \big]$$
	for all vector fields $X,Y$ on $M$.   If the metric of $M$ is Einstein and satisfies $\mathrm{Ric}_g = 2 n g$ then $M$ is said to be \emph{Sasaki-Einstein}.

\end{defn}

\begin{rmk}
The factor of 2 is unnecessary, but makes the definition compatible with the standard contact form of $\Sph^{2n+1}$.  
\end{rmk}

Immediate consequences of the definition are the following.

\begin{thm}
	Suppose $M$ carries the Sasakian structure $(g, R, J)$ and let $\alpha = g(R, \cdot)$ be the metric dual one-form of $R$.
	\begin{enumerate}
		
		\item The one-form $\alpha$ is a contact form whose contact structure $\Xi := \mathrm{Ker}(\alpha)$ is orthogonal to $R$ and has Reeb vector field equal to $R$.
	
		\item The characteristic foliation of $M$ consists of geodesics and the metric $g$ is Reeb-equivariant.
		
		\item The endomorphism $J$ is skew-symmetric, Reeb-equivariant, preserves $\Xi$, and satisfies $$\dif \alpha (X, Y) = g( JX, Y)$$ for all vector fields $X, Y$ on $M$.  Furthermore, $J^2 = - I + R \otimes \alpha$ where $I$ is the identity.
	
		\item The Riemann curvature of $g$ satisfies $\mathrm{Riem}_g(X, R)Y =  g(X,Y)R - g(R, Y)X$ for all vector fields $X, Y$ on $M$.
	\end{enumerate}
	
\end{thm}

From now on, suppose that the characteristic foliation of $M$ consists of smooth, 1-dimensional, embedded submanifolds; in this case the Sasakian structure of $M$ is said to be \emph{regular} and the fiber space $\hat M$ is once again a smooth manifold of dimension $2n$.   The canonical lifting process described earlier induces two additional structures on $\hat M$ that come from equivariant structures of $M$, namely  a metric $\hat{g}$ defined by $\hat{g}(\hat{X}, \hat{Y}) := g( \lambda(\hat{X}), \lambda(\hat{Y}) )$ as well as an endomorphism $\hat J$ of the tangent spaces of $\hat M$ defined by $\hat J ( \hat X) = \pi_\ast \circ\Phi \circ \lambda(\hat X)$.

\begin{thm}
	Suppose $(M, g, R, J)$ is a Sasakian manifold and let $\hat M$ be the fiber space of its characteristic foliation.  Let $(\hat g, \hat \omega, \hat J)$ be the structures on $\hat M$ defined above.
	\begin{enumerate}
		\item $(\hat M, \hat g, \hat \omega, \hat J)$ is a K\"ahler manifold of dimension $2n$.
		
		\item If $M$ is Sasaki-Einstein, then $\hat M$ is K\"ahler-Einstein with Ricci curvature $\mathrm{Ric}_{\hat g} = (2n+2) \hat g$.
	\end{enumerate}
\end{thm}

If a contact manifold $M$ possesses a Sasakian structure, then one can ask how the intrinsic and extrinsic geometry of its Legendrian submanifolds compares with that of their Lagrangian images under projection to the fiber space $\hat M$.  The next theorem uses classical facts about Riemannian submersions \cite{oneill} to answer this question.  Denote by $B_{V,W}$ the second fundamental form of a submanifold $V$ in an ambient manifold $W$. 
  
\begin{thm}
\mylabel{thm:meancurvlift}
Suppose $M$ is a Sasakian manifold and $\hat M$ is its fiber space.  If $\Lambda$ is a Legendrian submanifold of $M$ that projects onto a Lagrangian submanifold $[\Lambda] \subseteq \hat{M}$ under $\pi$, then
the following are true.
\begin{enumerate}
	\itemsep = 0ex
	\item The projection $\pi \big|_\Lambda : \Lambda \rightarrow [\Lambda]$ is a local isometry when $\Lambda$ and $[\Lambda]$ are given their induced metrics in $M$ and $\hat{M}$ respectively.
	
	\item The second fundamental forms of $\Lambda$ and $[\Lambda]$ satisfy $$B_{\Lambda, M} (X,Y) = \lambda \left( B_{[\Lambda], \hat{M}}( \pi_\ast (X), \pi_\ast (Y)) \right) \! \quad \mbox{and} \quad \, \pi_\ast \Big( B_{\Lambda, M} (\lambda(\hat X), \lambda (\hat Y) ) \Big) = B_{[\Lambda], \hat M} (\hat X, \hat Y)$$
for all $X,Y$ tangent to $\Lambda$ and $\hat X, \hat Y$ tangent to $[\Lambda]$.
\end{enumerate}
\end{thm}

\noindent An important conclusion to be drawn from this theorem is that any minimal and Legendrian submanifold of $M$ projects isometrically to a Lagrangian submanifold of $\hat M$ that is minimal; and any minimal and Lagrangian submanifold of $\hat M$ lifts isometrically to a family of Legendrian submanifolds of $M$ that are all minimal.

\paragraph{Contact stationary Legendrian submanifolds.}

In a Sasakian manifold, one can consider the variation of the volume of a Legendrian submanifold $\Lambda$ with respect to a restricted class of variations, namely the contact deformations of $\Lambda$.  The Euler-Lagrange equation of this variational problem is given in the following theorem, in essence proved by Schoen and Wolfson \cite{sw1}.

\begin{thm}
	\label{thm:eulerlagrange}
	Let $\Lambda$ be a Legendrian submanifold of a Sasakian manifold $M$.   If the volume of $\Lambda$ is stationary with respect to all contact deformations of $\Lambda$, then $\nabla \cdot \big( H_{\Lambda} \elbow \dif \alpha \big|_{\Lambda} \big) = 0$, where $H_\Lambda$ is the mean curvature of $\Lambda$.
\end{thm}

\begin{proof}
	By Theorem \ref{thm:meancurvlift} and the fact that contact deformations of a Legendrian submanifold $\Lambda$ of $M$ cover symplectic deformations of $\hat \Lambda$ in $\hat M$, it is sufficient to perform the calculation in the fiber space $\hat M$.  Without loss of generality, let $\hat \phi_u$ be a symplectomorphism generated by a function $u : [\Lambda] \rightarrow \R$ with Hamiltonian vector field $X_u = \hat J \nabla u$, and let $[\Lambda_t] = \phi_{tu}([\Lambda])$.  Then
	\begin{align}
		\label{eqn:firstvar}
		\left. \frac{\dif}{\dif t} \mathit{Vol} \big([\Lambda_t]) \big) \right|_{t=0} &= - \int_\Lambda \hat g (H_{[\Lambda]}, \hat J \hat \nabla u) \, \dif \mathrm{Vol}_{\hat g} \notag \\
		&= -\int_\Lambda u \, \hat \nabla \! \cdot \! \big( H_{[\Lambda]} \elbow \hat \omega \big|_{[\Lambda]} \big) \, \dif \mathrm{Vol}_{\hat g} \notag \\
		&= -\int_\Lambda u \, \nabla \! \cdot \! \big( H_{\Lambda} \elbow \dif \alpha \big|_{\Lambda} \big) \, \dif \mathrm{Vol}_{g}
	\end{align} 
	by Stokes' Theorem.  Note: formula \eqref{eqn:firstvar} is called the \emph{first variation of volume}.
\end{proof}

The contact-stationary submanifolds of $M$ are defined as the stationary points of the Euler-Lagrange equation determined in Theorem \ref{thm:eulerlagrange}.

\begin{defn}
	Let $\Lambda$ be a Legendrian submanifold of a Sasakian manifold $M$.  Then $\Lambda$ is called \emph{contact-stationary} if $\nabla \cdot \big( H_{\Lambda} \elbow \dif \alpha \big|_{\Lambda} \big) = 0$.
\end{defn}

In this paper, a solution of the equation $\nabla \cdot \big( H_{\Lambda} \elbow \dif \alpha \big|_{\Lambda} \big) = 0$ will be found in which $\Lambda$ is a contact deformation of an approximate solution $\Lambda_0$.   In a Sasakian ambient manifold, one has the notion of a \emph{normal deformation}; in other words, the infinitesimal vector field of the deformation is everywhere normal to $\Lambda$.  A procedure for generating normal contact deformations from functions on $\Lambda$ will be elaborated in a later section of this paper; it amounts to extending $u$ off $\Lambda$ in a canonical way and integrating the contact Hamiltonian vector field defined by \eqref{eqn:contactdef}.  If $u : \Lambda \rightarrow \R$ is a sufficiently small function (in a $C^1$ norm) then the contact Hamiltonian vector field can be integrated up to time one.  Denote the time-one contactomorphism simply by $\phi_u$.   Since the deformation $\phi_u$ depends on the first derivatives of $u$ and finding the divergence of the mean curvature of $\phi_u(\Lambda)$ takes another three derivatives, the operator 
\begin{equation}
	\label{eqn:defop}
	u \mapsto \nabla \cdot \left(  H_{\phi_u(\Lambda)} \elbow \dif \alpha \Big|_{\phi_u(\Lambda)} \right)
\end{equation}
is a non-linear, fourth-order, partial differential operator.  The following terminology for this operator will be used from now on.

\begin{defn}
	\label{defn:cminleg}
	The \emph{contact-stationary Legendrian deformation operator} of the Legendrian submanifold $\Lambda$ is defined to be the operator $\Phi_\Lambda : C^{4, \beta}(\Lambda) \rightarrow C^{0,\beta}(\Lambda)$ given in equation \eqref{eqn:defop}.
\end{defn}

\noindent Another important definition is of the linearization of $\Phi_\Lambda$ at $f$. 

\begin{defn}
	The \emph{linearization} of $\Phi_\Lambda$ at $f$ is the linear operator defined by the prescription $u \mapsto \frac{\dif}{\dif t} \Phi_\Lambda(f+ tu) \big|_{t=0}$.  It will be denoted $\Dif \Phi_\Lambda(f) (u)$ or $\mathcal L_\Lambda$ if $f$ is clear from the context.
\end{defn} 

\noindent This linear operator will now be calculated at $f=0$ in the general context of Sasakian geometry.   Einstein summation convention is used where necessary below: i.e.~repeated indices are summed, a comma denotes partial differentiation, a semi-colon denotes covariant differentiation with respect to the induced metric of $\Lambda$, indices are lowered and raised with this metric and its inverse, respectively, and so on.  The calculation proceeds in two steps.

\begin{prop}
	\label{prop:linearizationH}
	Let $\Lambda$ be a Legendrian submanifold of a Sasaki-Einstein manifold $M$ and suppose $\phi_u : M \rightarrow \R$ is a normal contact deformation generated by a function $u : \Lambda \rightarrow \R$.   The infinitesimal variation of the mean curvature tensor of $\Lambda$ is given by
	\begin{equation}
		\label{eqn:linearizationH}
		\left. \frac{\dif}{\dif t} \right|_{t = 0} \left(  H_{\Lambda_t} \elbow \dif \alpha \Big|_{\Lambda^t} \right) = \dif \big( \Delta_\Lambda u + (2n+2) u \big)
	\end{equation}
	where $\Lambda_t = \phi_{tu}(\Lambda)$ and $\Delta_\Lambda$ is the Laplace operator of $\Lambda$.
\end{prop}

\begin{proof}
According to Theorem \ref{thm:meancurvlift}, it is sufficient to calculate the infinitesimal variation of the mean curvature of the projected Lagrangian submanifold $[\Lambda]$ in the fiber space $\hat M$.  Moreover, it is easy to see that there is a one-parameter family of symplectomorphisms $\hat \phi^t : \hat M \rightarrow \hat M$ such that $\phi_{tu}$ covers $\hat \phi^t$ in a tubular neighbourhood of $[\Lambda]$, and that the Hamiltonian vector field $\hat X$ is normal to $[\Lambda]$ with $X = \hat J \nabla u$ along $[\Lambda]$, where $\nabla$ is the induced connection of $[\Lambda]$.
 
One begins by identifying a good local frame with which to perform the calculations.  Let $E_1, \ldots, E_n$ be a Riemannian normal coordinate frame for $[\Lambda]$ centered at a point $[p] \in [\Lambda]$.  Extend these vectors by parallelism to a neighbourhood so that the vectors $E_1, \ldots, E_n, \hat J E_1, \ldots, \hat J E_n$ span the tangent spaces of $\hat M$ there.  In this frame, the components of the second fundamental form $B_{[\Lambda]}$ and mean curvature $H_{[\Lambda]}$ are
$$B_{ijk} = \hat g (  \hat \nabla_{E_i} E_j , J E_k)  \qquad \mbox{and} \qquad H_{k} = h^{ij}  \hat g (  \hat \nabla_{E_i} E_j , J E_k)$$
where $\hat \nabla$ is the ambient connection of $\hat M$ and $h^{kl}$ are the components of the inverse of the induced metric $h$ of $[\Lambda]$, whose components are $h_{kl} = \hat  g (E_k, E_l) = \hat g(\hat JE_k, \hat JE_l)$ because $\hat J$ is an isometry.  Note that $B_{ijk}$ is symmetric in all its indices.  

The time derivative of the components $H_k$ will be computed at the point $[p]$ itself, where the Christoffel symbols of $h$ vanish, partial differentiation and covariant differentiation coincide, and therefore $\hat \nabla_{E_i} E_j = B_{ij}^{\hspace{1.5ex} k} \hat J E_k$.  Recall that it is possible to assume $\hat \nabla_X E_i = \hat \nabla_{E_i} X$, so that 
\begin{align}
	\label{eqn:hdot}
	\left. \frac{\dif}{\dif t}  H_k  \right|_{t = 0} &= X (h^{ij}) \hat g(\hat \nabla_{E_i} E_j, \hat J E_k) + h^{ij} X \hat g(\nabla_{E_i} E_j, \hat J E_k)  \notag \\
	&=2 B_{ijk} B^{ijl} u_{; l} + h^{ij} \big( \hat g (  \hat \nabla_X  \hat \nabla_{E_i} E_j, J E_k ) + \hat g(  \hat \nabla_{E_i} E_j, J  \hat \nabla_{X} E_k) \big)   \notag \\
	&=  h^{ij}  \big( \hat g (  \hat \nabla_{E_i}  \hat \nabla_{E_j} X , J E_k ) + u^{;l} \, \mathrm{Riem}_{\hat g} (E_i, J E_l, E_j, JE_k)\big) \notag \\
	&\qquad - H^l B_{lk}^{\hspace{1.5ex} m} u_{,m}  + 2 B_{ijk} B^{ijl} u_{; l}
\end{align}
since $\hat \nabla_{E_k} X = \hat J \hat \nabla_{E_k} \nabla u = u_{;k}^{\hspace{1.5ex} l} J E_l - u^{;l} B_{kl}^{\hspace{1.5ex} m} E_m$.  Here $\mathrm{Riem}_{\hat g}$ is the Riemannian curvature tensor of $\hat M$.  Calculate the remaining term containing $X$ as follows:
\begin{align}
	\label{eqn:hdotone}
	\hat g (  \hat \nabla_{E_i}  \hat \nabla_{E_j} X , J E_k ) &= \hat g (\hat \nabla_{E_i}  \hat \nabla_{E_j} \nabla u ,  E_k) \notag \\
	&= E_i \hat g ( \nabla_{E_j} \nabla u + B_{[\Lambda]}(E_j, \nabla u), E_k) - \hat g (  \nabla u + B_{[\Lambda]}(E_j, \nabla u), \nabla_{E_j} E_k)\notag \\
	&= h( \nabla_{E_i} \nabla_{E_j} \nabla u , E_k) - \hat g( B_{[\Lambda]}(E_j, \nabla u), B_{[\Lambda]}(E_i, E_k) ) \notag \\
	&= u_{;kij} - B_j^{\hspace{1ex} lm} B_{ikl} u_{;m}
\end{align}
Finally, perform the following manipulation on the $\mathrm{Riem}_{\hat g}$ term.  Start by adding and subtracting exactly the quantity required to yield the full Ricci curvature of $\hat M$.  That is,
\begin{align}
	\label{eqn:hdotthree}
	h^{ij}  \mathrm{Riem}_{\hat g} (E_i, \hat J E_l, E_j, \hat JE_k) &= h^{ij}  \mathrm{Riem}_{\hat g} (E_i, \hat J E_l, E_j, \hat JE_k) + h^{ij} \mathrm{Riem}_{\hat g} (\hat J E_i,  \hat J E_l, \hat J E_j, \hat J E_k) \notag \\
	&\qquad - h^{ij} \mathrm{Riem}_{\hat g} (\hat J E_i, \hat J E_l, \hat J E_j, \hat J E_k) \notag \\
	&= (2n +2) h_{kl} - h^{ij} \mathrm{Riem}_{\hat g}(E_i, E_l, E_j, E_k) \, .
\end{align}	
Complete the calculation by applying the Gau\ss\ equation, which reads
\begin{equation}
	\label{eqn:hdotfour}
	h^{ij} \mathrm{Riem}_{\hat g}(E_i, E_l, E_j, E_k) = R_{kl} - H^m B_{mkl} + B_k^{\hspace{1ex} mn} B_{lmn} \, ,
\end{equation}
where $R_{kl}$ are the components of the Ricci tensor of $\Lambda$.  Substituting \eqref{eqn:hdotone}, \eqref{eqn:hdotthree} and \eqref{eqn:hdotfour} into \eqref{eqn:hdot} gives
\begin{equation*}
	\left. \frac{\dif}{\dif t} H_k \right|_{t=0} = \Delta_\Lambda (u_{;k} ) - u_{;l} R^l_k - (2n+2) u_{;k} = \big( \Delta_\Lambda u \big)_{;k} - ( 2 n+2) u_{;k}\, .
\end{equation*}
This is the desired formula.
\end{proof}

\begin{cor}
	\label{cor:linop}
	Let $\Lambda$ be a Legendrian submanifold of a Sasaki-Einstein manifold $M$ and suppose $\phi_u : M \rightarrow \R$ is a normal contact deformation generated by a function $u : \Lambda \rightarrow \R$.   The infinitesimal variation of the divergence of the mean curvature tensor of $\Lambda$ is given by
	\begin{equation}
		\label{eqn:linearization}
		\left. \frac{\dif}{\dif t} \right|_{t = 0} \nabla \cdot \left(  H_{\Lambda_t} \elbow \dif \alpha \Big|_{\Lambda^t} \right) = \Delta_\Lambda \big( \Delta_\Lambda u + (2n+2) u \big) + Q_\Lambda(u)
	\end{equation}
	where $\Lambda_t = \phi_{tu}(\Lambda)$ and $Q_\Lambda$ is the operator given by $Q_\Lambda(u) = 2\,  (\nabla \cdot  B_\Lambda) (H_\Lambda, \nabla u)  - ( H_\Lambda \cdot \nabla )^2 u$.
\end{cor}

\begin{proof}
Compute the time derivative of $h^{ij} H_{i;j} = h^{ij} ( H_{i,j} - H_s \Gamma^s_{ij} )$ in terms of the derivatives of $h^{ij}$ and $\Gamma^s_{ij}$ in the same way as above.
\end{proof}

\begin{rmk}
Observe that the full symmetry of $B_\Lambda$ in all its slots implies that the operator $Q_\Lambda$ is self-adjoint, as one would expect.
\end{rmk}

A standard fact about the linearization of the contact-stationary Legendrian operator is that one-parameter families of isometries that are also contactomorphisms (which shall be called \emph{contact isometries}) produce elements  in the kernel of its linearization.

\begin{cor}
	\label{cor:jacorigin}
	Suppose $\phi_f^t : M \rightarrow M$ is a one-parameter family of contact isometries of $M$ generated by a function $f : M \rightarrow \R$.  If $\Lambda$ is a Legendrian submanifold of $M$ then $\Dif \Phi_\Lambda(0) ( f \big|_\Lambda) = 0$.
\end{cor}

\begin{proof}  
Since $\phi_f^t$ preserves the mean curvature and induced metric of Legendrian submanifolds, it is necessarily the case that the divergence of the mean curvature of  $\phi_f^t (\Lambda)$ is constant.   Hence
$$\left. \frac{\dif}{\dif t} \right|_{t=0} \nabla \cdot \left( H_{\phi_{f}^t (\Lambda)} \elbow \dif \alpha \Big|_{\phi_{f}^t(\Lambda)} \right) = 0  \, .$$
Without loss of generality $\phi_{f}^t$ can be replaced by $\phi_{\hat f}^t$ when $t$ is small, where $\hat f$ is some extension of the function $f \big|_{\Lambda}$ that generates a normal contact deformation whose infinitesimal deformation satisfies $X_{\hat f} = J \nabla f \big|_\Lambda$ on $\Lambda$.  Consequently $\Dif \Phi_\Lambda(0) ( f \big|_\Lambda) = 0$.
\end{proof}

\begin{rmk}
A one-parameter family $\phi^t_f$ of contact isometries preserves \emph{both} the induced metric of $\Lambda$ and the mean curvature of $\Lambda$.  Hence it is actually true that both $( \Delta_\Lambda  + 2(n+1) )( f \big|_{\Lambda} )= \mathit{constant}$ and $Q_\Lambda( f \big|_{\Lambda}) = 0$ which together imply $\Dif \Phi_\Lambda(0) ( f \big|_\Lambda) = 0$.
\end{rmk}

\subsection{The Geometry of the Unit Sphere} 
\label{subsec:cminleg}

\paragraph{Geometric structures.} The Calabi-Yau structure of $\C^{n+1}$ consists of the standard Euclidean metric $\delta$, the standard symplectic form $\omega_0$ and the standard complex structure $J_0$ satisfying  $\omega_0(X,Y) = \delta(J_0 X, Y)$, as well as the canonical holomorphic volume form $\Omega := \dif z^1 \wedge \cdots \wedge \dif z^{n+1}$.  The Calabi-Yau structure of $\C^{n+1}$ induces a Sasaki-Einstein structure on $\Sph^{2n+1}$.  First, denote the position vector field by  
$$P  := \sum_{k=1}^{n+1} \left( x^k \frac{\partial}{\partial x^k} + y^k \frac{\partial}{\partial y^k} \right) = \sum_{k=1}^{n+1} \left( z^k \frac{\partial}{\partial z^k} + \bar{z}^k \frac{\partial}{\partial \bar z^k} \right) \, ,$$
given both in real and complex coordinates.   Then the relevant objects are the following.
\begin{itemize}
	\item The metric is the standard metric of the sphere, induced from the ambient Euclidean metric.

	\item The contact form is $$\alpha := \frac{1}{2} P \elbow \omega_0 \big|_{\Sph^{2n+1}} = \displaystyle \frac{1}{2} \sum_{k=1}^{n+1} \Big( x^k \dif y^k - y^k \dif x^k \Big) \Big|_{\Sph^{2n+1}} = \frac{1}{4\mi} \sum_{k=1}^{n+1} \left( \bar z^k \dif z^k - z^k \dif \bar z^k \right) \Big|_{\Sph^{2n+1}} $$ so that $\dif  \alpha = \omega_0 \big|_{\Sph^{2n+1}}$.
		
	\item The Reeb vector field is $R := 2 J_0 P$ so that the contact structure is $\Xi_p := \big( \mathrm{span} ( J_0 P_p ) \big)^\perp$ for every $p \in \Sph^{2n+1}$.
	
	\item The endomorphism $J$ is defined to equal $J_0$ on $\Xi$ and to vanish in the Reeb direction. 

\end{itemize}

The characteristic foliation of $\Sph^{2n+1}$ as a Sasakian manifold coincides with the \emph{Hopf fibration} and the fiber space of $\Sph^{2n+1}$ coincides with the K\"ahler manifold $\cp^n$.  Furthermore, the Reeb projection coincides with the Hopf projection $\pi_H : \Sph^{2n+1} \rightarrow \cp^n$ which is the Riemannian submersion of the sphere onto $\cp^n$ with the Fubini-Study metric.  To see all this, recall that $\cp^n$ is the space of complex lines in $\C^{n+1}$.  Thus $\cp^n$ is the orbit space of the action of multiplication by non-zero complex numbers restricted to the sphere, i.e.~the $\Sph^1$ action $\theta \cdot p := \me^{\mi \theta} p$ for $p \in \Sph^{2n+1}$.  But the differential of this action is $\frac{\dif}{\dif \theta} (\theta \cdot p) \big|_{\theta = 0} = \mi p$ which is exactly the value of the vector field $JP$ at $p$ (in complex coordinates).  Hence the orbits of the action coincide with the characteristic foliation.

\paragraph{Contact isometries of the sphere.} The isometries of $\Sph^{2n+1}$ that preserve the Sasakian structure derive from the complex structure-preserving isometries of Euclidean space, namely the $\un$-rotations.  The one-parameter subgroups of $\un$ will play a crucial role in the sequel.  

Let $U^t$ be a one-parameter subgroup of $\un$-rotations with $U^0 = \mathit{Id}$ and recall the following facts.  There is a Hermitian matrix $H$ so that $U^t = \exp (\mi H t)$.  Also, $U^t$ is a one-parameter family of symplectomorphisms of $\C^{n+1}$ whose Hamiltonian vector field is given by $X_H(z): = \mi H z$ for any $z \in \C^{n+1}$.   The Hamiltonian function associated to $X_H$ is the Hermitian, harmonic, homogeneous polynomial of degree 2 given by $q_H(z) := z^\ast H z$.   Moreover,  $U^t$ restricts to a contactomorphism of $\Sph^{2n+1}$ with contact vector field $X_H$ and contact Hamiltonian equal to the restriction $q_H \big|_{\Sph^{2n+1}}$.  

Since $U(n+1)$-rotations of $\Sph^{2n+1}$ are both isometries and contactomorphisms, every minimal or contact-stationary Legendrian submanifold remains Legendrian under continuous $U(n+1)$-rotation and the mean curvature is unchanged.  By Corollary \ref{cor:jacorigin} one has $\mathcal L_{\Lambda} \big( q_H \big|_{\Lambda} \big)$ where $\mathcal L_{\Lambda}$ is the linearization of the contact-stationary Legendrian deformation operator.  Hence the kernel of this linearized operator is always non-trivial.

The following explicit specification of the generators of $\un$ and their associated Hermitian, harmonic, homogeneous polynomials of degree 2 on $\C^{n+1}$ will be needed in the sequel. A basis for the set of complex $(n+1) \times (n+1)$ Hermitian matrices is given by:
\begin{gather*}
	I := \left(\rule{0ex}{4ex}
	\begin{smallmatrix}
		1 & & \\[-1ex]
		& \ddots & \\[-0.5ex]
		& & 1
	\end{smallmatrix} 
	\right) 
	\quad
	H_1 :=  \left(\rule{0ex}{5ex}
	\begin{smallmatrix}
		n & && \\
		& -1 & & \\[-1ex]
		& & \ddots & \\
		& & & -1
	\end{smallmatrix} 
	\right)
	\quad H_2 :=  \left(\rule{0ex}{6ex}
	\begin{smallmatrix}
		0 & & & & \\
		& 1 & & & \\
		& & 0 & & \\[-1ex]
		& & & \ddots & \\
		& & & & -1
	\end{smallmatrix} 
	\right)
	\: \: \cdots \: \: \: H_n :=  \left(\rule{0ex}{6ex}
	\begin{smallmatrix}
		0 & & & & \\
		& 0 & & & \\[-1ex]
		& & \ddots & & \\
		& & & 1 & \\
		& & & & -1
	\end{smallmatrix} 
	\right) \, , \\
	\intertext{which are the real, diagonal matrices; along with}
	H_{jk} := \left(\rule{0ex}{5ex}
	\begin{smallmatrix}
		& && &   \\
		&&& 1 & \\
	 	& & & &\\
		& 1 & & &\\
		& & & & 
	\end{smallmatrix} 
	\right) 	
	\qquad \qquad
	H'_{jk} :=\left(\rule{0ex}{5ex}
	\begin{smallmatrix}
		& && &   \\
		&&& \mi & \\
	 	& & & &\\
		& - \mi & & &\\
		& & & & 
	\end{smallmatrix} 
	\right) 	 \, ,
\end{gather*}
which are the symmetric matrices having 1 in the $j^{\mathrm{th}}$ row and $k^{\mathrm{th}}$ column, and the anti-symmetric matrices having $\mi= \sqrt{-1}$ in the $j^{\mathrm{th}}$ row and $k^{\mathrm{th}}$ column.  Their associated polynomials are:
\begin{equation}
	\label{eqn:ambientpolys}
	\begin{aligned}
		q_0(z) &:= z^\ast I z = \sum_{s=0}^{n} |z^s|^2 \\[-1ex]
		q_1(z) &:= z^\ast H_1 z = n |z^1|^2 - \sum_{s=2}^{n+1} |z^s|^2 \\
		q_j (z) &:= z^\ast H_j z = |z^j|^2 - |z^{n+1}|^2 \qquad \mbox{for $j = 2, \ldots, n$} \\[1ex]
		q_{\textsc{re}, jk} (z) &:= z^\ast H_{jk} z = 2\, \Re ( z^j \bar z^k) \qquad \mbox{for $1 \leq j < k \leq n+1$} \\[1ex]	
		q_{\textsc{im}, jk} (z) &:= z^\ast H_{jk}' z = 2\, \Im ( z^j \bar z^k) \qquad \mbox{for $1 \leq j < k \leq n+1$} \, . 
	\end{aligned}
\end{equation}

\section{The Legendrian Connected Sum Procedure}
\label{sec:legconnect}

The present section of the paper describes the \emph{Legendrian connected sum procedure} that will be used to connect a chain of minimal Legendrian $\un$-rotated $n$-spheres together.  What will be presented here are the technical details of the procedure; concrete examples of initial configurations of $n$-spheres to which the procedure can be applied will be presented in Section \ref{sec:approxsol}.  Thus for the moment consider a single $n$-sphere $S_0$ and two neighbouring $n$-spheres $S_1 := U_1(S_0)$ and $S_2 := U_2(S_0)$ with the following properties.  There is a small $c \in (0,2 \pi)$ so that $S_0$ and $S_j$ have Hopf separation $c$ at the point $p_j$; and the tangent spaces $T_{p_j} S_0$ and $T_{p_j} \me^{(-1)^j \mi c} S_j$ are transverse and satisfy the angle criterion. The points $p_1$ and $p_2$ are the \emph{gluing points} on $S_0$.  Note that the Hopf separation between $S_0$ and $S_1$ is equal but opposite to the Hopf separation between $S_0$ and $S_2$.  Without loss of generality $S_0 := \R^{n+1} \times \{ 0 \} \cap \Sph^{2n+1}$.

\subsection{Preliminary Perturbation of the $n$-Spheres}
\label{subsec:prelimpert}

In order to connect $S_0$ to $S_1$ and $S_2$ at the gluing points with the best possible estimate of the mean curvature of the connected sum, it is first necessary to deform these $n$-sphere in a canonical manner.  The point is to `prepare' them for the connected sum operations by giving them a catenoidal shape near the points of connection.  To this end, let $\eps$ be a small, positive number (which will depend on $c$ in a manner to be determined below) and choose the unique distributional solution $G_0 : S_0 \setminus \{p_1, p_2\} \rightarrow \R$ of the equation $\Delta_S (G_0) = \delta_1 - \delta_2$ satisfying $\int_{S_0} G_0 = 0$, where $\delta_j$ is the Dirac $\delta$-mass at $p_j$ and $\Delta_S$ is the Laplacian of the induced metric of $S_0$.  Choose similar functions for each $S_j$, but using the gluing points where $S_j$ is to be connected to its two neighbours (one of which is $S_0$).  Henceforth consider only $S_0$ and the function $G_0$ since the analysis on the neighbouring $n$-spheres is the same. Note that
\begin{equation}
	\label{eqn:greengrowth}
	G_0(x) = 
	\begin{cases}
		 \dfrac{(-1)^{j+1}}{[ \mathrm{dist}(x, p_j) ]^{n-2}} + O(1) &\qquad n \geq 3 \\[3ex]
		 (-1)^j \log [ \mathrm{dist}(x, p_j) ] + O(1) &\qquad n=2
	\end{cases}
\end{equation}
in a neighbourhood of each $p_j$.

The idea is now to replace each $S_0$ by a Legendrian perturbation of $S_0$ generated by $\eps^{n} G_0$.  The way to do this is fairly simple and depends on the fact that the intersection of a Lagrangian cone in $\C^{n+1}$ with $\Sph^{2n+1}$ is a Legendrian submanifold.  Let $\bar G_0 : \R^{n+1} \times \{0 \} \rightarrow \R$ denote the degree-two homogeneous extension of $G_0$ to the $(n+1)$-plane containing $S_0$, namely the function defined by $\bar G_0( x) := \| x \|^2 G_0 (x / \| x \|)$ for $x \in \R^{n+1} \times \{0 \}$.  Therefore the submanifold
$$\overline{\mathit{Pert}_\eps(S_0)} := \{ (x, \eps^n \mathring \nabla \bar G_0(x) ) : x \in \R^{n+1} \} \, , $$
where $\mathring \nabla$ denotes the Euclidean gradient, is Lagrangian and a cone; and thus intersecting with $\Sph^{2n+1}$ yields a Legendrian submanifold.

\begin{defn}
	The \emph{preliminary perturbation} of $S_0$ with parameter $\eps$ is the Legendrian submanifold $\mathit{Pert}_\eps( S_0 ) := \overline{\mathit{Pert}_\eps(S_0)} \cap \Sph^{2n+1}$.
\end{defn}

The induced metric and second fundamental form of $S_0$ are trivial (i.e.~the standard metric and zero, respectively) and this changes when $S_0$ is replaced by $\mathit{Pert}_\eps( S_0 ) $.  But if $\eps$ is sufficiently small, then the geometry of $\mathit{Pert}_\eps( S_0 )$ changes only slightly, and this is true even relatively near to the gluing points.  The purpose of the following proposition is to quantify this statement.  The key to deriving the following results is to exploit the relationship between a geometric object on $\mathit{Pert}_\eps( S_0 )$ and the corresponding geometric object on $\overline{\mathit{Pert}_\eps( S_0 )}$.  Also $\| \cdot \|^\ast_{C^k(\mathcal O)}$ denotes the so-called \emph{scale-invariant} norm on a neighbourhood $\mathcal O$ in which the $j^{\mathit{th}}$ derivative is weighted by the factor $r^j$ where $r := \mathit{diam}(\mathcal O)$.

\begin{prop}
	\label{prop:exest}
	Let $r_0$ denote a radius so that \eqref{eqn:greengrowth} is valid inside $B_{r}(p_j)$ and set $S' := \mathit{Pert}_\eps( S_0 ) \setminus \big[ B_{r_0} (p_1) \cup B_{r_0}(p_2) \big]$.  Then the following estimates for geometric objects on $\mathit{Pert}_\eps(S_0)$ hold.
	\begin{itemize}
		\item The second fundamental form satisfies the estimate
		$$ \| B \|^\ast_{C^2(B_r(p_j))} \leq \frac{C \eps^n}{r^{n+1}} \quad \forall  \: r \in (0, r_0) \qquad \mbox{and} \qquad \| B \|_{C^2(S')} \leq C \eps^n \, .$$
		
		\item The mean curvature satisfies the estimate
		\begin{gather*}
			| \nabla \cdot H |^\ast_{C^1(B_r(p_j))} \leq \frac{C \eps^{3n}}{r^{3n+2}} \quad \forall \: r \in (0, r_0) \qquad \mbox{and} \qquad | \nabla \cdot H |_{C^1(S')} \leq C \eps^{3n} \\[2ex]
			\| H \|_{C^2(B_r(p_j))}^\ast \leq C \frac{\eps^n}{r^{n-1}} \qquad \mbox{and} \qquad \| H \|_{C^2(S')} \leq C \eps^n \, .
		\end{gather*}
		
		\item The Laplacian satisfies
		\begin{gather*}
			| \Delta (u) - \Delta_S (u) |^\ast_{C^1(B_r(p_j))} \leq \frac{C \eps^{2n}}{r^{2n+2}} | u |^\ast_{C^2(B_r(p_j))} \quad \forall \: r \in (0, r_0) \\[1ex]
			\mbox{and} \\
			| \Delta (u) - \Delta_S (u) |_{C^1(S')} \leq C \eps^{2n} |u|_{C^2(S')} 
		\end{gather*}
		where $\Delta_S$ is the Laplacian of the standard $n$-sphere.
	\end{itemize}
\end{prop}

\begin{proof}
 Let $\bar h, \bar \nabla, \bar \Delta, \bar B$ and $\bar H$ denote the induced metric, covariant derivative, Laplacian, second fundamental form and mean curvature of $\overline{\mathit{Pert}_\eps(S_0)}$ while $h, \nabla, \Delta, B$ and $H$ denote the same things for $\mathit{Pert}_\eps(S_0)$.  For typographical convenience, set $\bar G := \eps^n \bar G_0$ and $G := \eps^n G_0$.
 
Simple computations show that on $\overline{\mathit{Pert}_{\eps}(S_0)}$ the induced metric is $\bar h_{ij} := \delta_{ij} + \sum_s \bar G_{,si} \bar G_{,sj}$, the Christoffel symbols are $\bar \Gamma_{ijk} := \sum_s \bar G_{, ijs} \bar G_{,sk}$ and the components of the second fundamental form are $\bar B_{ijk} := \bar G_{,ijk}$.  Thus the inverse of $\bar h_{ij}$ equals the identity plus $O(\| \bar \nabla^2 \bar G \|^2)$ terms so that
$$\| \bar \nabla^k \bar B \| = O( \| \bar \nabla^{3+k} \bar G \| ) \, ,$$
all other terms being smaller.  Also, 
\begin{align*}
	\bar \nabla \cdot \bar H &= \bar h^{kl} \bar h^{ij} (\bar B_{,ijk})_{;l} \\
	&= \big[ \delta^{kl} + O(\| \bar \nabla^2 \bar G \|^2) \big] \, \big[ \delta^{ij} + O(\| \bar \nabla^2 \bar G \|^2) \big] \, \big[ \bar G_{, ijkl} - 2 \bar G_{,sjk} \bar \Gamma^{s}_{il} - \bar G_{,ijs} \bar \Gamma^{s}_{kl} \big] \\
	&= O(\| \bar \nabla^2 \bar G \|^2 \| \bar \nabla^4 \bar G \| ) + O( \| \bar \nabla^3 \bar G \|^2 \| \bar \nabla^2 \bar G \| ) \, .
\end{align*}
This is because $\delta^{kl} \delta^{ij} \bar G_{,ijkl} = \mathring \Delta ( \mathring \Delta \bar G)$ where $\mathring \Delta$ is the Euclidean Laplacian, and the choice of $G$ implies $\mathring \Delta ( \mathring \Delta \bar G) = \frac{1}{\| x \|^2} \Delta( \Delta + 2(n+1))(G) = 0$.  An estimate of $\| H \|$ follows in a similar way to the previous two estimates, and yields 
$$ \| H \| = O( \| \bar \nabla \bar G \|) \, .$$
This is not quite as good as before because the leading term in the expansion of $\| H \|$ is the derivative of the Laplacian of the function $\bar G$, which is the derivative of $\bar G$. Finally,
\begin{align*}
	\bar \Delta (u) - \mathring \Delta(u) &= \bar h^{ij} \big[ u_{,ij} - u_{,s} \bar \Gamma^{s}_{ij} \big] \\
	&= O( \|\bar  \nabla^2 \bar G \|^2) \cdot \bar \nabla^2 u + O( \| \bar \nabla^3 \bar G \| \cdot  \| \bar \nabla^2 \bar G \| ) \cdot \bar \nabla u
\end{align*}
where $O(\ast) \, \cdot $ denotes a linear operation on tensors with coefficients bounded by $O(\ast)$. 

Now using the homogeneity of $\bar G$ and the fact that $\overline{\mathit{Pert}_{\eps}(S_0)}$ is a cone over ${\mathit{Pert}_{\eps}(S_0)}$, one can replace $\| \bar \nabla^k \bar G \|$ by $\sum_{s=0}^k \| \nabla^s G \|$ to obtain expressions valid on $\mathit{Pert}_{\eps}(S_0)$.  Substituting $G = \eps^n G_0$ and using \eqref{eqn:greengrowth} then yields the desired result.  Higher derivative estimates follow in the same way.
\end{proof}

\subsection[Normal Coordinates for $\Sph^{2n+1}$]{Normal Coordinates for \boldmath $\Sph^{2n+1}$}

 \mylabel{subsec:normcoord}

In order to proceed with the Legendrian connected sum procedure, canonical coordinates for $\Sph^{2n+1}$ near the gluing points are needed in which all the relevant structures have an extremely simple form.  These coordinates will be called \emph{Legendrian normal coordinates} and are defined by first constructing a `standard' coordinate map and then transplanting it to a neighbourhood of a desired point $p$ by means of an $SU(n+1)$-rotation. 

Endow $\Cno$ with coordinates $z := (z^0, z')$ where $z' := (z^1, \ldots, z^n)$ and define the half-space $\mathcal{H} =   \{ z \in \Cno : \Re(z^0) \leq 0 \}$.   Next, endow $\C^n$ with coordinates $(w^1, \ldots, w^n)$ and let $B_1(0) \subseteq \C^n$ be the unit ball.  Define the map $K : B_1(0) \times (-\pi,\pi) \rightarrow \Sph^{2n+1} \setminus \mathcal H$ by $K(w, t) := \me^{\mi t} \big( f(w), w \big)$,  where $f : \C^n \rightarrow \R$ is given by $f(w) :=  (1 - \Vert w \Vert^2)^{1/2}$.  Then $K$ is invertible and the inverse is given by  $K^{-1}(z) =  \big( \me^{- \mi \arg(z^0)} z' , \arg(z^0) \big)$.   This is the `standard' map which is essentially graphical projection from the contact plane at $e_0$ into $\Sph^{2n+1}$ followed by motion along the Hopf fiber.  

Now let $p \in \Sph^{2n+1}$ be any point and suppose $\Pi$ is any Legendrian $n$-plane in $T_p \Sph^{2n+1}$.  Let  $e_0, e_1, \ldots, e_n$ be the standard basis of $\Cno$ and choose $V_{p,\Pi} \in \sun$ taking $e_0$ to $p$ and the real linear span of $e_1, \ldots, e_n$ to $\Pi$.  Note that $V_{p,\Pi}$ is unique up to orthogonal transformations and has the property that the contact hyperplane at $e_0$ (which is the complex linear span of $e_1, \ldots, e_n$) is mapped to the contact hyperplane at $p$ and the Reeb direction at $e_0$ (which is the real line spanned by $\mi e_0$) is mapped to the Reeb direction at $p$ (which is the real real line spanned by $\mi p$). The transformation $V_{p,\Pi}$ allows the map $K$ to be transplanted to the point $p$ in a manner adapted to $\Pi$.  Composing these maps leads to the desired coordinates.

\begin{defn}
	\mylabel{defn:normcoord} 
	The \emph{Legendrian normal coordinate chart} at $p \in \Sph^{2n+1}$ adapted to the Legendrian $n$-plane $\Pi \subseteq T_p \Sph^{2n+1}$ is the map $$\psi_{p,\Pi} : \Sph^{2n+1} \setminus V_{p,\Pi} (\mathcal{H}) \rightarrow B_1(0) \times (-\pi,\pi)$$ defined by $\psi_{p,\Pi} := K^{-1} \circ V_{p,\Pi}^{-1}$.
\end{defn}

The properties of Legendrian normal coordinates that will be used in the remainder of this paper are gathered in the following propositions.  Denote the set $\Sph^{2n+1} \setminus V_{p,\Pi} (\mathcal{H})$ by $\mathcal{H}_{p,\Pi}$.  Endow $\C^{n}$ with the standard symplectic structure $\omega_0 =   \frac{1}{2 \mi} \sum \dif \bar w^k \wedge \dif w^k$  and endow $\C^{n} \times \R$ with the contact structure defined by the contact form $\alpha_0 = \frac{1}{2} \dif t + \sigma_0$ where $\sigma_0 = \frac{1}{4 \mi}\sum \big( \bar{w}^k \dif w^k - w^k \dif \bar{w}^k \big)$.  Note that $\omega_0 = \dif \sigma_0$.  Finally, define the projection $\pi_0 : \C^n \times \R \rightarrow \C^{n}$ by $\pi_0(w,t) = w$, and for each $c \in \R$ define the injection $i_c : \C^n \rightarrow \C^n \times \R$ by $i_c(w) = (w, c)$.  The proof of the next proposition is just straightforward calculation and diagram chasing.

\begin{prop}
	\mylabel{prop:legnormcoordprops}
	Suppose $p \in \Sph^{2n+1}$ is and $\Pi \subseteq T_p \Sph^{2n+1}$ is Legendrian $n$-plane.  The Legendrian normal coordinate chart  $\psi_{p,\Pi} : \mathcal{H}_{p,\Pi} \rightarrow B_1(0) \times (-\pi,\pi)$ at $p$ has the following properties.  
	\begin{enumerate}
	
		\item The map $\psi_{p,\Pi} : (\Sph^{2n+1},\alpha) \rightarrow (\C^n \times \R, \alpha_0)$ is a contactomorphism.
		
		\item  If $\Lambda$ is a Legendrian submanifold tangent to the $n$-plane $\me^{\mi c} \Pi$ at $\me^{\mi c} p$ then $\psi_{p,\Pi}(\Lambda)$ is a Legendrian submanifold of $\C^n \times \R$ tangent to the plane $\R^n \times \{c\}$ at $(0,c)$.   Moreover, $\psi_{p,\Pi}$ takes $\un$-rotated Legendrian $n$-spheres passing through $\me^{\mi c} p$ to Legendrian $n$-planes passing through $(0,c)$.
		
		\item The map $\hat{\psi}_{[p],[\Pi]}  : (\cp^n, \hat{\omega}) \rightarrow (\C^n,\omega_0) $ defined by $\hat{\psi}_{[p,\Pi]}:= \big( \pi_H \circ \psi_{p,\Pi}^{-1} \circ i_0 \big)^{-1}$ is a symplectomorphism and the diagram
		\begin{equation} \mylabel{eqn:diagram}
			\begin{CD}
				\mathcal{H}_{p,\Lambda} \subseteq ( \Sph^{2n+1}, \dif \alpha) @>\psi_{p,\Pi} >>  
				(\C^n \times \R, 
				\dif \alpha_0)\\
				@VV\pi_H V \hspace{-1.5ex} @V\pi_0 VV \hspace{-6.5ex} @AA i_0 A\\
				 [\mathcal{H}_{p,\Pi}] \subseteq (\cp^n, \hat{\omega}) @>\hat{\psi}_{[p,\Pi]}
				   >> (\C^n, \omega_0)
			\end{CD}
		\end{equation}
		commutes.   The indicated differential forms correspond under pull-back  by the appropriate mappings.

	\item Let $g$ be the induced metric of $\Sph^{2n+1}$ and $\hat{g}$ be the Fubini-Study metric of $\cp^n$; and let $g_0 := \big(\psi_{p,\Lambda}^{-1} \big)^\ast g$ and $\hat g_0 := \big(\hat{\psi}_{[p],[\Lambda]}^{-1} \big)^\ast \hat{g}$ be the corresponding pull-back metrics.   Then $$\pi_0 : \big( \C^n \times \R, g_0 \big) \rightarrow \big( \C^n, \hat g_0 \big)$$ is a Riemannian submersion.  Moreover, these metrics are given by
	\begin{equation*}
		\begin{aligned}
			g_0 &= \dif t^2 + \delta + \dif f \otimes \dif f + 2 \big( \dif t \otimes \sigma_0 + \sigma_0 \otimes \dif t \big) \\
			\hat g_0  &=  \delta + \dif f \otimes \dif f - 4 \sigma_0 \otimes \sigma_0
		\end{aligned}
	\end{equation*}
	where $\delta$ is the Euclidean metric of $\C^n$. 		
	\end{enumerate}
\end{prop}

It will be necessary to estimate quite precisely the intrinsic and extrinsic geometry of a Legendrian submanifold $\Lambda$ in $\C^n \times R$ with respect to the induced metric $\hat g_0$.  A consequence of Lemma \ref{thm:meancurvlift} is that the necessary estimates can be gotten by estimating $ [\Lambda] \subseteq \C^n$ with respect to the metric $\hat g_0$, which amounts to a significant simplification since $\hat g_0$ is much simpler than $g_0$.  A further simplification becomes available if the Lagrangian projection of the submanifold to be estimated in this way is contained in a very small neighbourhood of the origin.   This is because $\hat g_0$ is in Riemannian normal form at the origin (namely $(\hat g_0)_{ij} = \delta_{ij}$ and $(\hat g_0)_{ij,k} = 0$ at the origin, for all $i, j, k$) by virtue of the fact that $g = \delta + Q$ where the coefficients of the tensor $Q$ satisfy $|Q_{ij}(z)| + \| z \|\,  | Q_{ij,k}(z) | + \| z \|^2 \, |Q_{ij,kl}(z)| = O(\|z \|^2)$.  Thus one would expect that the geometry of $[\Lambda]$ with respect to the metric $\hat g_0$ is uniformly close to the geometry of $[\Lambda]$ with respect to $\delta$.  The following proposition makes this idea rigorous. The following notation will be used: $H^g$, $B^g$, $\nabla^g$, $\Delta^g$, $\langle \cdot, \cdot \rangle_g$ and $\| \cdot \|_g$ denote the mean curvature, second fundamental form, covariant derivative, Laplacian, inner product and norm induced by a metric $g$.  

\begin{prop}
	\label{prop:meancurvcomparison}
	There exists a small $r_0>0$ so that the following holds for all $0< r < r_0$.  Let $\hat \lambda : \mathcal O \rightarrow \C^n$ be a Lagrangian embedding of a neighbourhood $\mathcal O \subseteq \R^n$ such that $\hat \lambda (\mathcal O) \subseteq B_r (0)$ and put $h = \hat \lambda ^\ast \hat g_0$.  
	\begin{itemize}
	
		\item The mean curvature and second fundamental form of $\mathcal O$ satisfy the following estimates.  	
	\begin{equation*}
		\begin{aligned}
			\| B^{\hat g_0} \|^\ast_{C^2(\mathcal O, h)}  &= \| B^{\delta} \|^\ast_{C^2(\mathcal O, \hat \lambda^\ast \delta)} + O(r) \\
			\| H^{\hat g_0} \|^\ast_{C^2(\mathcal O, h)} &= \| H^{\delta} \|^\ast_{C^2(\mathcal O,  \hat \lambda^\ast \delta)}  + O(r) + O(r^2) \| B^\delta \|^\ast_{C^2(\mathcal O, \hat \lambda^\ast \delta)} \, .
		\end{aligned}
	\end{equation*}
	
		\item The Laplacian of $\mathcal O$ satisfies 
		\begin{equation*}
			\Delta^h u = \Delta^\delta u + O(r) \cdot\nabla^\delta u + O(r^2) \cdot \nabla^\delta \nabla^\delta u  \, .
		\end{equation*}
		where $O(\ast) \, \cdot $ denotes a linear operation on tensors with coefficients bounded by $O(\ast)$.
	\end{itemize}
	
\end{prop}

\begin{proof}  The proof of this result is a rather diabolical exercise in Riemannian geometry that will be abbreviated here for the sake of the reader.  First, pick any point $p \in \lambda(\mathcal O)$ and choose a $\delta$-orthonormal frame $E_i$ for $T_p \lambda (\mathcal O)$.  Extend this to a $\delta$-orthonormal frame for $T_p \C^n$ by adjoining the vectors $J E_i$.  One can assume that $\langle \nabla^\delta_{E_i} E_j , E_k \rangle_\delta = 0$ and that $B_{i j k}^\delta = \langle \nabla^\delta_{E_i} E_j , J E_k \rangle_\delta$. 

The second fundamental form of $\lambda(\mathcal O)$ with respect to the metric $ h$ can be found by computing $B_{ijk}^{ h} = \langle \nabla^{ h}_{E_i} E_j , N^{ h}_k \rangle_{ h}$ where $N_k^{ h}$ is a $h$-orthonormal frame for the  $h$-orthogonal complement of $T_p \lambda(\mathcal O)$.  To complete this calculation, one uses the form of the metric to compute $\langle \nabla^{ h} _X Y, Z \rangle_h  = \langle \nabla^\delta _X Y, Z \rangle_\delta + O(r)\|X\|_\delta \|Y\|_\delta \|Z\|_\delta$ for any vectors $X, Y, Z$ whereas the $N_k^{ h}$ can be written as linear combinations of $J E_k + c^i_k E_i$ where $c^i_k = O(r^2)$. Substituting these quantities into the expression for $B_{ijk}^{ h}$ and its covariant derivatives yields
\begin{align*}
	B_{ijk}^{ h} &= B_{ijk}^\delta + O(r) \\
	\nabla_s^{ h} B_{ijk}^{ h} &= \nabla^\delta_s B_{ijk}^\delta  + O(1) \\
	\nabla^{ h}_s \nabla^{ h}_t B_{ijk}^{ h}&= \nabla^\delta_s \nabla^\delta_t B^\delta_{ijk} + O(r^{-1})
\end{align*}
where $O(r^\ast)$ can now represents a tensor of the required type having coefficients bounded by $r^\ast$.  Next, the form of the metric $h$ implies that the inverse matrix has coefficients $ h^{i j} = \delta^{ij } + O(r^2)$.  One thus obtains the desired estimates by taking the $ h$- and $\delta$-norms of these expressions and making the necessary comparisons.  The expressions for the mean curvature follow by first taking the $h$- and $\delta$-traces and then taking norms.  The calculations for the Laplacian are similar.
\end{proof}

\subsection{The Legendrian Lawlor neck}
\label{subsec:leglawlorneck}

A \emph{Lagrangian Lawlor neck} is an embedded special Lagrangian submanifold of $\C^n$ with the topology of a cylinder $\R \times \Sph^{n-1}$ that is asymptotic to a pair of transverse special Lagrangian planes.  It was discovered by Lawlor in \cite{lawlor1,lawlor2}, and it was shown that any given pair of transverse special Lagrangian planes $\Pi_1$ and $\Pi_2$ can be realized as the asymptotic planes of a Lagrangian Lawlor neck if and only if the \emph{characteristic angles} of $\Pi_1$ and $\Pi_2$ satisfy the \emph{angle criterion} as described in the introduction.   A truncated and re-scaled Lawlor neck can be lifted from $\C^n$ into $\C^n \times \R$ by the standard Legendrian lifting procedure to produce an embedded, cylindrical Legendrian submanifold asymptotic to a pair of Legendrian planes.  This object can then be embedded into the sphere as a Legendrian submanifold using the Legendrian normal coordinates.  The resulting object is the \emph{Legendrian Lawlor neck}.   The purpose of this section is to define this object precisely and to compute good estimates of its mean curvature and other geometric quantities.  

\paragraph{The Lagrangian Lawlor neck.}  The discussion begins with the definition of the Lagrangian Lawlor neck, which proceeds as follows.  First, let $A := (a_{1}, \ldots, a_{n})$ be a vector of positive real numbers and let $P : \R^{n} \times \R \rightarrow \R$ be the function given by
\begin{equation}\label{eqn:P}
	P(A, \lambda) := \frac{\left( 1+a_{1} \lambda^{2} \right)  \cdots \left( 1+a_{n} \lambda^{2} \right) - 1}{\lambda^{2}} \, .
\end{equation}
Next, set
\begin{equation}\label{eqn:angle}
    \theta_{k}(A, \lambda) := \int_{-\infty}^{\lambda} \frac{ \dif s} {(\frac{1}{a_{k}} + s^{2}) \sqrt{P(A, s)}} \, .
\end{equation}
It is easy to see that the integrals \eqref{eqn:angle} are well-defined and converge as $\lambda \rightarrow \infty$.  Moreover, $P(A, \lambda) = O(\lambda^{2n-2})$ and so $\theta_k(A, \lambda) = \theta_{k}(A)  + O(\lambda^{-n})$ for large $\lambda$, where $\theta_{k}(A) := \lim_{\lambda \rightarrow \infty} \theta_{k}( A, \lambda)$.  The embedding giving the Lagrangian Lawlor neck in $\C^n$ is now defined as follows.

\begin{defn}
	\label{defn:leglawlor}
	For every $A = (a_1, \ldots, a_n) \in \R^{n}$, with $a_{k} > 0$ for all $k$, the \emph{Lagrangian Lawlor embedding} with parameter $A$ is the map $\hat F_{A} : \R \times \Sph^{n-1} \rightarrow \C^n$ given by
	\begin{equation*}
		\hat F_{A} (\lambda, \mu) := \left( z^{1}(\lambda, \mu), \ldots , z^{n}(\lambda, \mu) \right)
	\end{equation*}
   	where 
    	\begin{equation*}
		z^{k}(\lambda, \mu) = \mu^{k} \, \sqrt{ \tfrac{1}{a_{k}} + \lambda^{2} } \, \exp\big( {\mi \theta_k(A, \lambda)} \big)
	\end{equation*}
    	and $\mu = (\mu^{1}, \ldots, \mu^{n}) \in \R^{n}$ satisfies $\sum (\mu^{k})^{2} = 1$ so that $\mu \in \Sph^{n-1}$.  The \emph{Lagrangian Lawlor neck} with parameter $A$ is the submanifold $\hat N_{A} := \hat F_{A}(\R \times \Sph^{n-1})$.   The \emph{truncated and re-scaled Lagrangian Lawlor} neck with parameters $\eps, R, A$ is the submanifold $\hat N_{\eps, R, A} := \eps \hat F_A \left( [-R, R] \times \Sph^{n-1} \right)$. 
\end{defn}

\noindent Observe that the two ends of $\hat N_A$ are asymptotic to the real plane $\R^n \times \{0\}$ (when $\lambda \rightarrow -\infty$) and to the plane $D_A (\R^n \times \{0\})$, where $D_A$ is the diagonal $n \times n$ matrix having the entries $\me^{\mi \theta_1(A)}, \ldots, \me^{\mi \theta_n(A)}$ on the diagonal (when $\lambda \rightarrow \infty$).  Moreover, $\sum_{k=1}^n \theta_k(A) = 0 \mod 2 \pi$ so that $D_A \in \sun$. 

\paragraph*{Lifting to a Legendrian submanifold.} In order to lift the Lagrangian submanifold $\hat N_A$ of $\C^n$ to a Legendrian submanifold of $\C^n \times \R$, it is necessary to find a function $S_A(\lambda, \mu)$ so that $F_A( \lambda, \mu) := \big(\hat F_A(\lambda, \mu), S_A(\lambda, \mu) \big)$ defines a Legendrian embedding.  If $\alpha_0 = \frac{1}{2} \dif t + \sigma_0$ is the contact form, then the condition $F_A^\ast \alpha_0 = 0$ implies that the required function must satisfy $\dif S_A(\lambda, \mu) = -2 \hat F_A^\ast \sigma_0$.  One then easily calculates that the required function must be 
\begin{equation}
	\label{eqn:Sfunction}
	S_A(\lambda) := - n \int_0^\lambda \frac{\dif s}{\sqrt{P(A, s)}}
\end{equation}
up to the addition of a constant factor.  It is important to note that $S_A(\lambda) = S_A + O(\log(\lambda))$ for large $\lambda$ in dimension $n=2$, and $S_A(\lambda) = S_A + O(\lambda^{-n+2})$ for large $\lambda$ in dimension $n \geq 3$, where $S_A$ is a constant depending only on $A$ and the dimension.

\begin{defn}
	For every $A = (a_1, \ldots, a_n) \in \R^{n}$ with $a_{k} > 0$ for all $k$, the \emph{Legendrian Lawlor embedding} with parameter $A$ is the map $F_{A} : \R \times \Sph^{n-1} \rightarrow \C^n \times \R$ given by
	\begin{equation*}
		F_{A} (\lambda, \mu) := \left(\hat F_A(\lambda, \mu), S_A(\lambda)  \right) \, .
	\end{equation*}
  	The \emph{Legendrian Lawlor neck} with parameter $A$  is the submanifold $N_{A} := F_{A}(\R \times \Sph^{n-1})$. 
\end{defn}

Although the Lagrangian Lawlor embedding into $\C^n$ can be scaled in a straightforward way, scaling the Legendrian Lawlor embedding into $\C^n \times \R$ requires handling the $\R$-coordinate and the $\C^n$-coordinates differently.   Indeed, to preserve the contact form it is necessary to scale the Legendrian Lawlor embedding as $F_{\eps, A}(\lambda, \mu) := \left( \eps \hat F_{A} (\lambda, \mu) , \eps^2 S_A(\lambda) \right)$.  The \emph{truncated and re-scaled Legendrian Lawlor neck} with parameters $\eps, R, A$ is the submanifold $N_{\eps, R, A} := F_{\eps, A} \left( [-R, R] \times \Sph^{n-1} \right)$.

\paragraph{Properties of the Legendrian Lawlor neck.}  The following propositions gather the relevant properties of $N_{\eps, R, A}$ that will be used in the sequel.  The first of these propositions  gives a more detailed picture of the asymptotics of $N_{\eps, R, A}$.  The analogous results for $\hat N_{\eps, R, A}$ have been proved in \cite{me1}, and their extension to the present case follows easily from Proposition \ref{prop:leg}.  Let $\Pi_0:= \linspanr \{e_1, \ldots, e_n\}$ and $\Pi_A := \linspanr \{ \me^{\mi \theta_1(A)} e_1, \ldots, \me^{\mi \theta_n(A)} e_n\}$.

\begin{prop}
	\mylabel{prop:lawlorembeddingsize}
	If $\eps>0$ is sufficiently small, there exist radii $R \gg 1$ and $r \ll 1$ satisfying $R = O(r/\eps)$ so that the following hold.  
	\begin{enumerate}        
	\item $\hat{N}_{\eps,R,A} \subseteq B_r(0)$ in $\C^n$ and 
	$N_{\eps,R,A} \subseteq B_r(0) \times (0, \eps^2 S_A)$ in $\C^n \times \R$.
	
	\item	$N_{\eps,R,A} \setminus B_{r/2}(0)$ consists of two disconnected components projecting onto $\Pi_0 \cap \annr$ and $\Pi_{A} \cap \annr$, respectively, by nearest-point projection.  Furthermore, there exists a function $G_\eps : \mathit{Ann}(r/2,r) \rightarrow \R$ so that
\begin{equation}
	\mylabel{eqn:graph}
	N_{\eps,R,A} \cap \annr =  \big( \Lambda_{G_\eps} \cap \annr \big) \bigcup \left( D_{A,\eps^2 S_A} \big( \Lambda_{G_\eps} \big) \cap \annr \right)
\end{equation}
in the notation of Proposition \ref{prop:leg}, where $D_{A,\eps^2 S_A}$ is the diagonal $SU(n)$-rotation bringing $\Pi_0$ to $\Pi_A$ followed by upward translation by $\eps^2 S_A$ in the Reeb direction.

\item The function $G_\eps$ exhibits the scaling behaviour $G_\eps(x) = \eps^2 G_1(x/\eps)$, and if $x \in \annr$, then it satisfies the estimates 
\begin{equation}
	\label{eqn:primitiveest}
	\Vert \nabla^k G_\eps(x) \Vert = O \left( \frac{\eps^n}{ \Vert x \Vert^{n-2+k}}  \right) 
\end{equation} 
for all $n\geq 2$ and $k\geq 0$, except $n=2$ and $k=0$ with respect to the metric induced by the ambient Euclidean metric.  In the exceptional case, the estimate is $|G_\eps(x)| = O(\eps^2 |  \log \|x\| \, | )$. 
\end{enumerate}
\end{prop}

The next proposition derives the intrinsic and extrinsic geometry of $N_{\eps, R, A}$ induced by the metric $g_0$.  Recall that this is the same as the geometry of $\hat N_{\eps, R, A} \subseteq \C^n$ that is induced by the metric $\hat g_0$, and that this metric is close to being Euclidean sufficiently close to the origin.  The first lemma describes the geometry of $\hat N_{\eps, R, A}$ with respect to the Euclidean metric of $\C^n$.  The formul\ae\ below all result from straightforward computations. 

\begin{lemma}
	\label{lemma:laglawlorgeom}
	The induced metric of the truncated Lagrangian Lawlor neck $\hat N_{\eps, R, A}$ with respect to the Euclidean metric of $\C^n$ is given by 
	\begin{equation}
		\label{eqn:lawlormetric}
		h_{\eps, A} = \eps^2 \left[\sum_{k=1}^{n} \left( \frac{(\mu^{k})^{2}} {\frac{1}{a_{k}} + \lambda^{2} } \right) \left( \lambda^{2} + \frac{1}{P(A, \lambda)} \right)  \dif \lambda^{2} + \sum_{k=1}^{n} \left( \frac{1}{a_{k}} + \lambda^{2} \right) \bigl( \dif \mu^{k} \bigr)^{2} \right] \, .
	\end{equation}
	The second fundamental form of $\hat N_{\eps, R, A}$ with respect to the Euclidean metric of $\C^n$ is given by 
	\begin{equation}
		\label{eqn:lawlorsecondff}
		B_{\eps, A} = \eps^2 \left[ \frac{\mathrm{Sym}(\dif \lambda \otimes  g_{S})}{\sqrt{P(A, \lambda)}} + b \, \dif \lambda^3 \right]
	\end{equation}
	where $g_{S}$ is the standard metric of $\Sph^{n-1}$ while $\mathrm{Sym}( \dif \lambda \otimes  g_{S})$ is the symmetrization of $\dif \lambda \otimes  g_{S}$ and 
	$$b =  \sum_{k=1}^n \frac{(\mu^k)^2}{\sqrt{P(A, \lambda)} \left( \frac{1}{a_k}+  \lambda^2 \right)^2} \Bigg[ \lambda^2 + \frac{1}{P(A, \lambda)} - \left( \frac{1}{a_k} + \lambda^2 \right) \left(1 + \frac{\lambda P'(A, \lambda)}{2 P(A, \lambda)} \right) \Bigg]  \, .$$
	Finally, the mean curvature of $\hat N_{\eps, R, A}$ with respect to the Euclidean metric of $\C^n$ vanishes.
\end{lemma}

\noindent A consequence of Lemma \ref{lemma:laglawlorgeom} and the methodology of Proposition \ref{prop:meancurvcomparison} is the following estimate for the second fundamental form with respect to the ambient metric induced from the sphere.

\begin{prop}
	\label{prop:neckest}
	The following estimates hold at a point $(\lambda, \mu)$ on $N_{\eps, R, A, c}$. Here $C$ is a constant independent of $\eps$.
	\begin{itemize}
	 	\item The second fundamental form of $N_{\eps, R, A, c}$ measured with respect to the metric $g_0$ satisfies
		\begin{equation*}
			\begin{aligned}
				\| B_{\eps, A} \| + \eps (1 + \lambda^2)^{1/2} \| \nabla B_{\eps, A} \| + \eps^2 (1+\lambda^2) \| \nabla^2 B_{\eps, A} \| & \\
				&\hspace{-20ex} \leq C \left[ \frac{1}{\eps (1+\lambda^2)^{(n+1)/2}} +  \eps(1+\lambda^2)^{1/2} \right] \, .
			\end{aligned}
		\end{equation*}
		
		\item The mean curvature of $N_{\eps, R, A, c}$ measured with respect to the metric $g_0$ satisfies
		\begin{equation*}
			\begin{aligned}
				|H_{\eps, A} | + \eps (1 + \lambda^2)^{1/2} \| \nabla H_{\eps, A} \| + \eps^2 (1+\lambda^2) \| \nabla^2 H_{\eps, A} \| & \\
				&\hspace{-20ex} \leq C \left[  \frac{\eps}{(1+\lambda^2)^{(n-1)/2}} +  \eps(1+\lambda^2)^{1/2} \right] \, .
			\end{aligned}
		\end{equation*}
	\end{itemize}
\end{prop}
	
\begin{proof}
The computation of these estimates proceeds in two steps.  First, one finds the estimate of the second fundamental form of $\hat N_{\eps, R, A}$ with respect to Euclidean metric of $\C^n$.  Then Proposition \ref{prop:meancurvcomparison} and with the fact that the mean curvature with respect to the Euclidean metric of $\C^n$ is zero is used to derive the needed estimates with respect to the ambient metric coming from the sphere.

For the first calculation, it suffices to take $A = (1, \ldots , 1)$ because it is always the case in this paper that the components of $A$ lie in a compact subset of the open positive quadrant of $\R^n$.  Thus one can derive from equations \eqref{eqn:lawlormetric} and \eqref{eqn:lawlorsecondff} that the coefficients of the metric and second fundamental form satisfy
\begin{align*}
	h_{\eps, A} &= \eps^2 \big[ O(1) \, \dif \lambda^2 + (1 + \lambda^2) \, g_S \big] \\
	B_{\eps, A} &= \eps^2 \big[ O \bigl( (1 + \lambda^2)^{(1-n)/2} \bigr) \,  \mathrm{Sym}(\dif \lambda \otimes g_S) + O\bigl( (1 + \lambda^2)^{(-1-n)/2} \bigr)\,  \dif \lambda^3 \big] \, .
\end{align*}
Taking the norm of $B_{\eps, A}$ with $h_{\eps, A}$ as well as the norm of its first two covariant derivatives with respect to $h_{\eps, A}$ yields the estimate
$$ 
\| B_{\eps, A} \| + \eps (1 + \lambda^2)^{1/2} \| \nabla B_{\eps, A} \| + \eps^2 (1+\lambda^2) \| \nabla^2 B_{\eps, A} \| \leq C \eps^{-1} (1+\lambda^2)^{-(n+1)/2}
\, . 
$$
Proposition \ref{prop:meancurvcomparison} can now be used to convert this into an estimate of the second fundamental form and mean curvature with respect to the metric $\hat g_0$, which yields exactly the required estimate. Needed for this task is the fact that the ambient coordinate $r$ is related to the neck coordinate $\lambda$ by $r = \eps (1+\lambda^2)^{1/2}$.
\end{proof}	

%




%

\subsection{Performing the Connected Sum}
\label{subsec:legconnect}

Let $S_0, S_1$ and $S_2$  be as in Section \ref{subsec:prelimpert} and replace these by $\mathit{Pert}_{\eps}(S_0), \mathit{Pert}_{\eps_1}(S_1)$ and $\mathit{Pert}_{\eps_2}(S_2)$.  It will now be shown how $\mathit{Pert}_{\eps}(S_0)$ can be connected to $\mathit{Pert}_{\eps_1}(S_1)$ and $\mathit{Pert}_{\eps_2}(S_2)$ by choosing $\eps, \eps_1, \eps_2$ well and inserting the appropriate Legendrian Lawlor necks.

Consider only the connected sum construction near the point $p_1 \in S_0$ and $\me^{\mi c} p_1 \in S_1$, the construction at $p_2$ and $\me^{-\mi c} p_2 \in S_2$ being similar.  Choose Legendrian normal coordinates $\psi  : \mathcal U \rightarrow \C^{n} \times \R$ in a neighbourhood $\mathcal{U}$ of $p_1$ and adapted to $S_0$ by choosing $\psi := \psi_{p_1, T_{p_1} S_0}$.  Then $\psi(\mathit{Pert}_{\eps} (S_0) \cap \mathcal{U})$ is a Legendrian submanifold graphical over the Legendrian plane $\Pi_0 := \R^n \times \{0\}$ near the origin, and $\psi(\mathit{Pert}_{\eps_1} (S_1) \cap \mathcal U)$ is a Legendrian submanifold graphical over a Legendrian plane of the form $\Pi_A := D_{A,c} \big( \Pi_0 \big)$ near the origin.  Here  $D_{A,c}$ is rotation in the  $\C^{n}$ coordinates by the diagonal matrix having $\me^{\mi \theta_k(A)}$ on the diagonal, where $\theta_1(A), \ldots, \theta_n(A)$ are the asymptotic angles for some Legendrian Lawlor neck $N_{A}$, followed by translation by $c$ in the $\R$ coordinate.  More concretely, there is some fixed small number $r>0$ along with functions $f_0 : B_r(0) \rightarrow \R$ and $f_1 : B_r(0) \rightarrow \R$ so that
\begin{equation*}
	\psi(\mathit{Pert}_{\eps} (S_0) \cap \mathcal{U}) \cap B_r(0)  = \Lambda_{f_0} \qquad \mbox{and} \qquad \psi(\mathit{Pert}_{\eps_1} (S_1)\cap \mathcal U) \cap B_r(0) = D_{A,c} \big( \Lambda_{f_1} \big)
\end{equation*}
where $\Lambda_f$ denotes the Legendrian submanifold generated by the function $f$ as in the notation of Proposition \ref{prop:leg}.  Moreover, by \eqref{eqn:greengrowth} and a change of variables, one can assume that the functions $f_0$ and $f_1$ can be expanded as
\begin{equation}
	\label{eqn:graphfn}
	f_j (x) =
	\begin{cases}
		\dfrac{C_j \eps_j^n}{\| x \|^{n-2} } + O (\eps_j^n \| x \|^2)  &\qquad n \geq 3 \\[3ex]
		\eps_j^2 \big( C_j^0 + C_j^1 \log(\| x \|) \big) + O(  \eps_j^2 \| x \|^2) &\qquad n=2
	\end{cases} 
\end{equation}
in a small neighbourhood of the origin. 

One can now attach $\Lambda_{f_0}$ to $\Lambda_{f_1}$ using a truncated and re-scaled Legendrian Lawlor neck as follows. Choose $\eps, \eps_1, \eps_2$ to satisfy $\eps^n c_0 = \eps_1^n c_1 = \eps_2^n c_2$ as well as the equation $c := \eps^2 S_A$ (or $c := \eps^2 S_A - \eps_1^2 C_1^0)$ in the $n=2$ case). Next, set 
$$r_\eps := \eps^{s} \qquad \mbox{for} \quad s \in (0,1)$$ 
to be the radius at which the gluing takes place.  The parameter $s$ will be chosen at the end of the proof; in the mean time, it suffices to know $\eps \ll r_\eps$ when $\eps$ is sufficiently small. Let $\eta : [0, 1] \rightarrow \R$ be a smooth, monotone cut-off function with support in $[0, 1/2]$.  Define the following two functions
\begin{align*}
	G_0(x) &:= \eta (\|x\| / r_\eps ) G_\eps (x) + f_0(x) \big( 1- \eta( \|x\| / r_\eps) \big) \\
	G_1(x) &:= \eta ( \|x\| / r_\eps) \big( G_\eps (x) + \eps^2 S_A \big) + \big( f_1(x) + c \big)  \big( 1- \eta( \|x\| / r_\eps) \big) 
\end{align*}
where $G_{\eps}$ is the graphing function of the Legendrian Lawlor neck from Proposition \ref{prop:lawlorembeddingsize} (with the appropriate modification for the $n=2$ case).  These two functions transition smoothly between the values of $f_0$ and $f_1+ \eps^2 S_A$ outside $B_{r_\eps}(0)$ and the values of $G_\eps$ inside $B_{r_\eps/2}(0)$, respectively.  Now define the submanifolds
$$ T_0: = \Lambda_{G_0} \cap \mathrm{Ann}(r_\eps/2, r_\eps) \qquad \mbox{and} \qquad T_1 :=  D_A \left( \Lambda_{G_1} \right) \cap  \mathrm{Ann}(r_\eps/2, r_\eps)$$
which transition smoothly between the Legendrian Lawlor neck $N_{\eps, R, A, \eps^2 S_A} \cap B_{r_\eps/2}(0)$ and the submanifolds $\Lambda_{f_0}$ and $\Lambda_{f_1}$ outside of $B_{r_\eps}(0)$, respectively.  Finally, define the submanifold
\begin{equation}
	\Lambda_\eps := \big( \Lambda_{f_0} \cap \mathrm{Ann}(r_\eps, r) \big) \cup T_0 \cup \left( N_{\eps, R, A, \eps^2 S_A} \cap B_{r_\eps/2}(0) \right) \cup T_1 \cup \big( D_{A, \eps^2 S_A}(\Lambda_{f_1}) \cap \mathrm{Ann}(r_\eps, r) \big) 
\end{equation}
which combines these five pieces into a smooth, Legendrian submanifold of $\C^n \times \R$.  Applying the inverse of the Legendrian normal coordinates then places this submanifold back in $\Sph^{2n+1}$.

\begin{defn}
	\label{defn:legconnect}
	Choose $\eps, \eps_1, \eps_2$ as above.  The \emph{Legendrian connected sum} of $S_0$ and $S_1$ is
	$$\lcs_\eps (S_0, S_1) := \big[\mathit{Pert}_{\eps} (S_0) \cup \mathit{Pert}_{\eps_1} (S_1) \setminus \psi^{-1}(B_r(0) ) \big]  \cup \psi^{-1} (\Lambda_\eps) \, .$$
	Call $\mathcal T := \psi^{-1}(T_0)\cup \psi^{-1}(T_1)$ the \emph{transition region} of $\lcs_\eps (S_0, S_1)$.
\end{defn}
	
The submanifold $\lcs_\eps (S_0, S_1)$ is identical to $\mathit{Pert}_{\eps} (S_0)$ or $\mathit{Pert}_{\eps_1} (S_1)$ outside of $ \psi^{-1} \big( B_r(0) \big)$ and identical to $\psi^{-1} \big( N_{\eps, R, A} \big)$ inside the neck region $\psi^{-1} \big( B_{r_\eps/2}(0) \big)$.  Therefore the estimates of the geometry of  $\lcs_\eps (S_0, S_1)$ are provided by Proposition \ref{prop:neckest} in the neck region and Proposition \ref{prop:exest} in the exterior region.   The following proposition gives the needed estimates in the transition region.

\begin{prop}
	\label{prop:transmeancurvest}
	The following estimates hold in the transition region $\mathcal T$ of $\mathit{LCS}_\eps (S_0, S_1)$.
	\begin{itemize}
	
		\item The second fundamental form satisfies
		$$ \| B \|^\ast_{C^2(\mathcal T)} \leq C \left( \frac{ \eps^n}{r_\eps^{n+1}} +  r_\eps \right)$$
		for some constant $C$ independent of $\eps$.

		\item The mean curvature satisfies
		\begin{equation*}
			\begin{aligned}
				| \nabla \cdot H |^\ast_{C^1(\mathcal T)} &\leq C \left( \frac{\eps^{n+1}}{r_\eps^{n+3}}  + 1 \right) \\[3ex]
				\| H \|_{C^2(\mathcal T)}^\ast &\leq C \left(\dfrac{\eps^n}{r_\eps^{n-1}} + r_\eps \right) 
				\end{aligned}
		\end{equation*}
		for some constant $C$ independent of $\eps$.
				
		\item The Laplacian satisfies
		$$ \Delta (u) = \mathring \Delta (u) + O \left( \frac{\eps^{2n}}{ r_\eps^{2n} } + r_\eps^2 \right) \cdot \mathring \nabla^2 u + O \left( \frac{\eps^{2n}}{ r_\eps^{2n+1}} + r_\eps \right) \cdot \mathring \nabla u $$
		where $O(r^\ast) \, \cdot $ denotes a linear operation on tensors with coefficients bounded by $O(r^\ast)$. 
		
	\end{itemize}
\end{prop}

\begin{proof}
The computation of these estimates will be carried out in Legendrian normal coordinates centered on the point of gluing.  One begins by performing the calculations for the Lagrangian projection $[\mathit{LCS}_\eps (S_0, S_1)] \subseteq \C^n$ using the ambient Euclidean metric.  The conclusion then follows from Theorem \ref{thm:meancurvlift} and Proposition \ref{prop:meancurvcomparison} as before.   

Suppose now that $f := \eta f_1 + (1-\eta) f_2$ where $f_1$ is the graphing function for the image of $[ \mathit{Pert}_\eps (S_0) ]$ in the Legendrian normal coordinate chart being used here, and $f_2$ is the graphing function for the end of the Lagrangian Lawlor neck that connects to it.  Let $\eta$ be the cut-off function appearing in the construction of $\mathit{LCS}_\eps(S_0, S_1)$.  Set 
$$L_{f_\ast} := \{ (x, \mathring \nabla f_\ast (x) ) : x \in B_{r_\eps}(0) \setminus B_{r_\eps/2}(0) \}$$
where $\mathring \nabla$ is the Euclidean gradient and $\ast$ refers to 1, 2 or nothing.  Since the components of the second fundamental form of $L_f$ with respect to the Euclidean metric are $B_{ijk} = f_{,ijk}$ then
$$\| B \| \leq \| \mathring \nabla^3 f_1 \| + \| \mathring \nabla^3 f_2 \| + \| T \|$$
where $T$ is a tensor of order three whose components are linear combinations of the quantities $\eta_{,i} ( f_{1,jk} - f_{2,jk})$, $\eta_{,ij} (f_{1, k} - f_{2,k})$ and $\eta_{, ijk} (f_{1} - f_{2})$ with $O(1)$ coefficients.  To proceed, recall the expansions \eqref{eqn:primitiveest} and \eqref{eqn:graphfn} which yield
$$\| \mathring \nabla^k f_j (x) \| = O \left( \frac{\eps^n}{\| x\|^{n-2+k}} \right) \qquad \mbox{and} \qquad \| \mathring \nabla^k (f_1 - f_2) (x) \| = O \left( \frac{\eps^{n+1}}{\| x\|^{n-1+k}} \right)  + O(\eps^n \| x \|^2) \, .$$
Hence with $r_\eps/2 \leq \| x \| \leq r_\eps$, one finds
$$ \| B \|_{C^0(\mathcal T)} \leq C \left( \frac{ \eps^n}{r_\eps^{n+1}} + \frac{\eps^{n+1}}{r_\eps^{n+2}} + \frac{\eps^n}{r_\eps}  \right) \, .$$
The estimates of the higher derivatives of $B$ are similar, though more involved.  Using Proposition \ref{prop:meancurvcomparison} and the fact that $\eps < r_\eps$, the desired estimate follows.

Next, the divergence of the mean curvature of $L_f$ with respect to the background Euclidean metric is $\nabla \cdot H = h^{kl} h^{ij} ( f_{,ijk} )_{;l}$ where $h^{ij}$ are the coefficients of the inverse of the induced metric $h_{ij} := \delta_{ij} + \sum_s f_{,is} f_{, js}$  and the covariant derivatives are taken with respect to this metric.  The Christoffel symbols are $\Gamma_{ijk} := \sum_s f_{,ijs} f_{,ks}$.  Expanding in terms of $f_1$ and $f_2$ yields after some work
\begin{equation*}
	| \nabla \cdot H | \leq | \nabla_1 \cdot H_1 | + | \nabla_2 \cdot H_2 | + \| T' \|  + O ( \| \mathring \nabla^2 f \|^2 \| \mathring \nabla^4 f \| ) + O( \| \mathring \nabla^3 f \|^2 \| \mathring \nabla^2 f \| )
\end{equation*}
where $\nabla_j \cdot H_j$ is the divergence of the mean curvature of $L_{f_j}$ with respect to the induced metric of $L_{f_j}$ with respect to the background Euclidean metric and this time $T'$ is a tensor of order four whose components are linear combinations of the quantities $\eta_{,i} ( f_{1,jkl} - f_{2,jkl})$, $\eta_{,ij} (f_{1, kl} - f_{2,kl})$, $\eta_{, ijk} (f_{1,l} - f_{2,l})$ and $\eta_{,ijkl} (f_1 - f_2)$ with $O(1)$ coefficients.   Taking $r_\eps / 2 \leq \|x \| \leq r_\eps$ now gives
$$| \nabla \cdot H |_{C^0( \mathcal T ) } \leq C \left(  \frac{ \eps^{3n}}{r_\eps^{3n+2}} + \frac{\eps^{n+1}}{r_\eps^{n+3}} + \frac{\eps^n}{r_\eps^2} \right)
$$
using the estimate $|\nabla_1 \cdot H_1|*_{C^1(\mathcal T)} \leq C \eps^{3n} r_\eps^{-3n-2}$ from Proposition \ref{prop:exest} and the fact that $H_2 = 0$.   The second derivatives of $H$ can be estimated in the same way.  Combined with the second fundamental form estimate above along with  Proposition \ref{prop:meancurvcomparison} and the fact that $\eps < r_\eps$, the first of the desired estimates follows.  The estimate of $\| H \|_{C^2(\mathcal T)}^\ast$ follows in a similar way to the previous two estimates, though yields a result that is not quite as good because the leading term is the derivative of the Laplacian of the function $G$ which is of size $O( \| \mathring \nabla G \|) = O(\eps^n r_\eps^{-n+1})$.

Finally, the Laplacian is $\Delta (u) := h^{ij} \big( u_{,ij} - \Gamma_{ij}^s u_{,s} \big)$ so that expanding as before yields
$$\Delta (u) = \mathring \Delta (u) + O( \| \mathring \nabla^2 f \|^2 ) \mathring \nabla^2 u + O( \| \mathring \nabla^3 f \| \cdot \| \mathring \nabla^2 f \| ) \cdot \mathring \nabla u \, .$$
The desired estimate follows by estimating $\mathring \nabla^k f$ as before and invoking Proposition \ref{prop:meancurvcomparison}.
\end{proof}

\section{The Approximate Solutions}
\label{sec:approxsol}

\subsection{The Initial Configurations}
\label{sec:examples}

Examples of initial configurations of Legendrian $n$-spheres to which the Main Theorem applies will now be given.  Denote the standard real basis of $\C^{n+1}$ by $P_1 := (1, \ldots, 0) \ldots  P_{n+1} := (0, \ldots, 1)$ in what follows, treating these as points in $S_0$ or as vectors in $\C^{n+1}$  depending on the context. 

\paragraph*{Example 1.}  The first example of an initial configuration can be described as follows.  Choose complex numbers $\xi_s := \exp (2 \pi \mi / m_s)$ for $s = 2, \ldots, n$ and $\bar \xi_{n+1} := \exp (2 \pi \mi \sum_{s=2}^{n} 1/m_s)$, where $m_s$ are odd integers with $|m_s| \geq 4$.   The n $U$ is defined to be
\begin{equation*}
	U := \mbox{ \footnotesize $ \left(
	\begin{array}{cccc}
		1 &&&\\
		& \xi_2&&\\
		&&\ddots & \\[-1ex]
		&&&\xi_{n+1}
	\end{array}
	\right) $}
\end{equation*}
Next, let $N := \lcm(m_2, \ldots m_{n})$, choose $k$ and $a$ odd satisfying $\gcd (a, kN) = 1$ and define $\zeta := \exp( 2 \pi \mi a / kN)$.  Then $(\zeta U)^{k N} = I$ and no power of $\zeta U$ smaller than $kN$ is the identity.

The $U^s (S_0)$ all intersect at the points $\pm P_1$ so that the $(\zeta U)^s (S_0)$ acquire Hopf separation on the order of  $\mathit{Arg}(\zeta)$.  In fact, if $s_1 := (kN+1)/2$ and $s_2 := (kN-1)/2$ then $(\zeta U)^{s_1}(P_1) = \zeta^{1/2} P_1$ and $(\zeta U)^{s_2}(P_1) = \zeta^{-1/2} P_1$ with Hopf separation $\pm \pi a / k N$.  Moreover $(\zeta U)^s (S_0)$ is closest to $(\zeta U)^{s_1}(S_0)$ and $(\zeta U)^{s_2}(S_0)$ at these points.  Note that $T_{P_1} U^{s_1} (S_0) = \linspanr \{ \xi_2^{1/2} P_2, \ldots, \xi_{n+1}^{1/2} P_{n+1} \}$ and $T_{-P_1} U^{s_2} (S_0) = \linspanr \{ \xi_2^{-1/2} P_2, \ldots, \xi_{n+1}^{-1/2} P_{n+1} \}$ and it is clear that these tangent planes are always transverse.  One can also check that the angle criterion holds for these tangent spaces for a sufficiently large choice of $m_k$, depending on the dimension.   The gluing should be performed by connecting the points given above.  One can check that the resulting submanifold is closed and embedded (because the Hopf separation between any pair of $n$-spheres is never smaller $\pi a / k N$).  The Legendrian angle function advances by $\pi a / kN$ from one $n$-sphere to the next and acquires a period of $2 \pi a$ around the entire configuration.

\paragraph{Example 2.}  The next example is different from the first, in that gluing occurs at a pair of non-antipodal points.  Once again, choose complex numbers $\xi_s := \exp (2 \pi \mi / m_s)$ for $s = 2, \ldots, n$ but now $\bar \xi_{n+1} := - \exp ( 2 \pi \mi \sum_{s=2}^{n} 1/m_s)$, where $m_s$ are integers with $|m_s| \geq 4$.   Then $U$ is defined to be
\begin{equation*}
	U := \mbox{ \footnotesize $ \left(
	\begin{array}{cc|ccc}
		0 & \xi_2 &&&\\
		1 & 0 &&&\\
		\hline &&\xi_3&& \\[-1ex]
		&&&\ddots & \\[-1ex]
		&&&&\xi_{n+1}
	\end{array}
	\right) $}
\end{equation*}
Next, let $N := \lcm(m_2, \ldots m_{n})$, choose $k$ and $a$ satisfying $\gcd (a, kN) = 1$ and define $\zeta := \exp( 2 \pi \mi a / kN)$.  Then $(\zeta U)^{2 k N} = I$ and no power of $\zeta U$ smaller than $2kN$ is the identity.

In the present case, $U(S_0) \cap S_0 = \{ \pm P_2 \}$ and $U^{-1}(S_0) \cap S_0 = \{ \pm P_1 \}$.  Indeed, one finds that $U^{2s} (P_1) = \xi_2^s P_1 = U^{2s-1} (P_2)$ and $U^{2s} P_2 = \xi_2^s P_2 = U^{2s + 1} (P_1)$.  Multiplying by $\zeta$ creates the desired Hopf separation.  In fact $\zeta^{2s} \xi_2^s P_1 \in (\zeta U)^{2s} (S_0)$ is closest to $\zeta^{2s-1} \xi_2^s P_1 \in (\zeta U)^{2s-1}(S_0)$ and $\zeta^{2s} \xi_2^s P_2 \in (\zeta U)^{2s} (S_0)$ is closest to $\zeta^{2s+1} \xi_2^s P_2 \in (\zeta U)^{2s+1}(S_0)$.  The Hopf separation is $\pm 2\pi a / k N$, respectively.  Note that $T_{ P_1} S_0 = \linspanr \{ P_2, \ldots, P_{n+1} \}$ and $T_{P_2} S_0 = \linspanr \{ P_1, P_3, \ldots, P_{n+1} \}$ while $T_{P_2} U (S_0) = \linspanr \{ \xi_2 P_1, \xi_3 P_3, \ldots, \xi_{n+1}^{1/2} P_{n+1} \}$ and $T_{P_1} U^{-1} (S_0) = \linspanr \{ \bar \xi_2 P_2, \ldots, \bar \xi_{n+1}^{1/2} P_{n+1} \}$. It is clear that these tangent planes are always transverse.  One can also can check that the angle criterion holds for these tangent spaces for a sufficiently large choice of $m_k$, depending on the dimension. The gluing should be performed by connecting the points given above.  One can check that the resulting submanifold is closed but embedded only for appropriate choices of $\xi_k$ and $\zeta$.  The Legendrian angle function acquires a period of $4 \pi a$ around the entire configuration.

\subsection{Assembling the Approximate Solutions}
\label{subsec:assembly}

The approach outlined in the introduction for creating an approximately contact-stationary Legendrian submanifold from $\un$-rotated copies of $S_0$ will now be implemented here in the two examples just given.  Broadly speaking, these examples will be constructed as follows.  

\begin{enumerate}

	\item Choose $U \in \sun$ generating a cyclic subgroup of $\sun$ of integer order $N$.  Suppose that $S_0 \cap U(S_0)$ has contact-transverse intersection at a pair of antipodal points $\pm p \in S_0$ in such a way that the tangent spaces at the intersection points satisfy the angle criterion. 
		
	\item Let $\zeta := \me^{2 \pi \mi a / k N}$ for some integer $a$ and large integer $k$ satisfying $\gcd(a, kN) = 1$.   
		
	\item Glue each $(\zeta U)^s (S_0)$ to its nearest neighbours at one of the points $(\zeta U)^s(\pm p)$ using the Legendrian connected sum procedure with scale parameter $\eps_k$ corresponding to the Hopf separation $2 \pi a / k N$. 
	
\end{enumerate}

\noindent The \emph{approximate solution} of the contact-stationary deformation equation is the submanifold constructed by means of the three steps above and will be denoted by $\approxsol$.

  The submanifold $\approxsol$ is smooth and Legendrian, and if $U$ and $\zeta$ are chosen properly, then it is embedded and each $(\zeta U)^s (S_0)$ is connected to exactly two neighbouring $\zeta U$-rotated copies of $S_0$ in such a that $\approxsol$ is a closed and compact Legendrian submanifold with the topology of $\Sph^1 \times \Sph^{n-1}$.  Since the Legendrian angle function $\Theta_{U, \zeta}$ on the unperturbed part of $(\zeta U)^s(S_0)$ is equal to $2 \pi s a / kN$, then $\Theta_{U, \zeta}$ increases by $2 \pi a$ for every circuit around $\approxsol$. 

The following terminology will be used in the sequel.  Each $(\zeta U)^s(S_0)$ is attached to two neighbours at two distinct points chosen from the set $\{ (\zeta U)^t(\pm p) : t = 0, \ldots, kN - 1\} \cap (\zeta U)^s(S_0)$.   These points are the gluing points on $(\zeta U)^s(S_0)$ and these shall be denoted by $p_{1s}$ and $p_{2s}$.  (In the case $s=0$, denote them simply by $p_1$ and $p_2$.)  Now choose $r \in (0, r_0]$ where $r_0$ is sufficiently small but independent of $\eps$, and define the disjoint union of balls centred on the gluing points
\begin{equation}
	\label{eqn:balls}
	\mathcal B_r := \bigcup_{s=0}^{kN- 1} \left[ B_{r} \big( ( \zeta U)^s(p_{1s}) \big) \cup B_{r} \big( ( \zeta U)^s(p_{2s}) \big) \right] \cap \approxsol \, .
\end{equation}
Finally, subdivide $\approxsol$ into regions of three distinct types: the union of all the neck regions, the union of all the transition regions and the union of all the exterior regions of $\approxsol$ and denote these regions by $\mathcal N$, $\mathcal T$ and $\Lambda'$.  That is, set
\begin{equation*}
	\begin{aligned}
		\mathcal N &:= \mathcal B_{r_\eps / 2} \\
		\mathcal T &:=  \mathcal B_{r_\eps} \setminus \mathcal B_{r_\eps/2} \\
		\Lambda' &:= \approxsol \setminus \Big( \mathcal N \cup \mathcal T \Big) \, .
	\end{aligned}
\end{equation*}

\subsection{Symmetries Satisfied by the Approximate Solutions}
\label{sec:additional}

The approximate solution $\approxsol$ constructed above is invariant with respect to the transformation $\zeta U$ for two reasons.  First, the collection of $n$-spheres $(\zeta U)^s(S_0)$ is invariant under this group for obvious reasons.  This shows that all of $\approxsol$ except perhaps the neck and transition regions are invariant under this group.  Second, since every neck region is uniquely determined by the tangent planes at the points of gluing and the Legendrian normal coordinates are equivariant with respect to $\zeta U$ (in fact, $\zeta U$ can be used to transplant the coordinates from the point $p$ to every other gluing point), the neck regions are permuted amongst themselves by the action of $\zeta U$.  Third, one can easily check that the Legendrian connected sum procedure keeps the transition regions of $\approxsol$ invariant with respect to this group as well.

There are also additional symmetries satisfied by $\approxsol$ that will be important later on.  These are specific to each of the two choices of $U$ made above.

\paragraph{Example 1.}  Consider the transformations $K_j : \C^{n+1} \rightarrow \C^{n+1}$ defined by
\begin{equation}
	\label{eqn:extrasym}
	K_j (z) := \big( z^1, z^2, \ldots ,- z^j, \ldots, z^{n+1} \big) \quad \mbox{for $j = 1, \ldots, n+1$} \, .
\end{equation}
These transformations are not $\un$-rotations, but rather are orthogonal transformations of $\C^{n+1}$ (viewed as $\R^{2n+2}$) that still preserve the contact structure.  One can easily verify that $K_j^{-1} \circ \zeta U \circ K_j = \zeta U$ for all $j$, so that each $K_j$ descends to an isometry of the fundamental cell of $\approxsol$ consisting of the spherical region, a transition region and a neck region.  It remains to see that each of these regions is invariant under $K_j$.

The spherical region $\mathit{Pert}_\eps(S_0) \setminus \big(B_{r_\eps} (-P_1) \cup B_{r_\eps}(P_1) \big)$ is clearly invariant under these symmetries.  Next, combining the Legendrian normal coordinate map from  Section \ref{subsec:normcoord} with the definition of the Legendrian Lawlor neck from Section \ref{subsec:leglawlorneck} shows that the neck $\psi(N_{\eps, R, A, 0})$ is embedded into $\Sph^{2n+1}$ by the mapping 
$$\psi \circ F_{\eps, R, A, 0}( \lambda, \mu) = \exp ( \mi S_A(\lambda) ) \left( R_A(\lambda, \mu), \cdots \mu^j \sqrt{\tfrac{1}{a_j} + \lambda^2} \exp (\mi \theta_j(A, \lambda) ) \cdots \right) \, ,$$
where $R_A (\lambda, \mu) = \left[ 1 - \sum_{j=1}^n (\mu^j)^2 \left( \tfrac{1}{a_j} + \lambda^2 \right) \right]^{1/2}$.  Now it is easy to see that
\begin{equation}
	\label{eqn:neckpresj}
	K_j \circ \phi \circ F_{\eps, R, A, 0}(\lambda, \mu) = \phi \circ F_{\eps, R, A, 0}(\lambda, \mu^1, \ldots, -\mu^j, \ldots, \mu^n) \qquad j= 2, \ldots, n+1
\end{equation}
showing that the transformations $K_j$ for all $j = 2, \ldots, n+1$ preserve the neck.  These transformations also preserve the transition regions of $\approxsol$ because the transition regions are rotationally symmetric about the vertical axis in the Legendrian normal coordinates. 

Finally, the transformation $K_1$ preserves $\mathit{Pert}_\eps(S_0)$ because the Green's function of the Laplacian of $S_0 \setminus \{ P_1, - P_1\}$ is odd.  Moreover, $K_1$ sends $P_1$ to $-P_1$ and acts by reflection across the horizontal hyperplane in the Legendrian normal coordinate chart.  Hence the neck at $P_1$ is sent to a reflection of itself across the horizontal hyperplane.  But this is exactly the neck that is inserted at $-P_1$.  Consequently $K_1$ preserves the neck regions of $\approxsol$.  It also preserves the transition regions of $\approxsol$ by similar arguments.

\paragraph{Example 2.} The transformations $K_j $ for $3 \leq j \leq n+1$ defined in \eqref{eqn:extrasym} from Example 1 are additional symmetries of $\approxsol$ in this case as well.  One additional symmetry involves the $z^1$ and $z^2$ coordinates.  Define the transformations $K_{12} : \C^{n+1} \rightarrow \C^{n+1}$ by 
$$K_{12}(z^1, z^2, z^3, \ldots, z^{n+1}) := ( -\bar z^2, -\bar z^1, \bar z^3, \ldots, \bar z^{n+1}) \, .$$
Then it is easy to check that $K_{12} ^{-1} (\zeta U) K_{12} = (\zeta U)^{-1}$ so that $K_{12}$ descends to an isometry of the fundamental cell of $\approxsol$.  The spherical part of the fundamental cell is clearly preserved by $K_{12}$.  Moreover, calculating with the embedding of the neck as in \eqref{eqn:neckpresj} shows that $K_{12}$ maps the neck $N_A$ to the neck $(\zeta U)^{-1} (N_A)$, but embedded up-side down with a reflection into the sphere.  Consequently, the neck regions are permuted amongst each other under $K_{12}$.  Similar arguments show that the transition regions are permuted in the same way under $K_{12}$ as well.

\section{Setting Up the Analysis}

\subsection{The Banach Space Inverse Function Theorem}
\label{sec:strategy}

In the previous two sections, an approximately contact-stationary submanifold $\approxsol$ has been constructed using the Legendrian connected-sum procedure.  The tasks ahead are to parametrize small contact deformations of $\approxsol$ over a suitable Banach space of $C^{4, \beta}$ functions of $\approxsol$ and solve the equation $\Phi_{\approxsol} (f) = 0$, where $\Phi_{\approxsol}$ is the contact-stationary deformation operator of $\approxsol$ as defined abstractly in Definition \ref{defn:cminleg}.   Henceforth, denote $\Phi_{U, \zeta} := \Phi_{\approxsol}$ and denote also the linearization of this operator  at zero by $\mathcal L_{U, \zeta}$ for typographical convenience.  

The theorem that will be invoked to find $f$, albeit in a slightly subtle way, is the \emph{Banach space inverse function theorem}.   This fundamental result will now be stated in fairly general terms \cite{amr}. 

\begin{nonumthm}[IFT]
	Let $\Phi : X \rightarrow Z$ be a smooth map of Banach spaces, set $\Phi(0) := E$ and denote the linearized operator $\Dif \Phi(0)$ by $\mathcal L$.  Suppose that $\mathcal L $ is bounded, surjective, and possesses a right inverse $\mathcal R : Z \rightarrow X$ satisfying the estimate 
	\begin{equation}
		\label{eqn:iftestone}
		\Vert \mathcal R(z) \Vert \leq C \Vert z \Vert
	\end{equation}
	holds for all $z \in Z$.  Choose $R$ so that if $y \in B_R(0) \subseteq X$, then 
	\begin{equation}
		\label{eqn:iftesttwo}
		\Vert \mathcal L ( x ) - \Dif \Phi(y)( x ) \Vert \leq \frac{1}{2C}  \Vert x \Vert
	\end{equation}
for all $x \in X$.  Then if $z \in Z$ is such that 
	\begin{equation}
		\label{eqn:iftestthree}
		\Vert z - E \Vert \leq \frac{R}{2C} \, ,
	\end{equation}
there exists a unique $x \in B_R(0)$ so that $\Phi(x) = z$.  Moreover, $\Vert x \Vert \leq 2C \Vert z - E \Vert$.
\end{nonumthm}

\subsection{Weighted Schauder Norms}
\label{subsec:funcspaces}

A preliminary step in applying the Banach space inverse function theorem is to have norms for the Banach spaces $X$ and $Z$ which make explicit the dependence on the parameter $k$ in all the estimates.  The reason this is necessary is because when $k \rightarrow \infty$ and $\eps_k \rightarrow 0$, then the approximate solution becomes singular in each neck region and the contact-stationary Legendrian equation degenerates, making it impossible to invoke the Banach space inverse function with more `conventional' norms.  

\begin{rmk}
The neck size parameter $\eps_k$ is more geometric that the parameter $k$. Therefore the dependence of all estimates on $\eps_k$ will henceforth be tracked, rather than their dependence on $k$.  In the sequel, denote $\eps := \eps_k$ for simplicity.
\end{rmk}

The $\eps$-dependence of the estimates will be tracked by using \emph{weighted Schauder norms} for the spaces $C^{l,\beta}(\approxsol)$.   To define the \emph{weight function}, first define a `regularized' distance function $\rho_0$ on $\Lambda'$ by requiring $\rho_0(x) := \mathrm{dist}(x, p_{js})$ in a neighbourhood of each of the gluing points $p_{js}$ and allowing $\rho_0$ to transition smoothly to one outside these neighbourhoods.   Here the terminology from Section \ref{subsec:assembly} for the various regions of $\approxsol$ is being used.   Now make the following definition.

\begin{defn} 
	\mylabel{defn:weightfn}
	Define the \emph{weight function} $\rho_\eps : \approxsol \rightarrow \R$ by
	\begin{equation}
		\mylabel{eqn:weightfn}
		\rho_\eps(x) := 
		\begin{cases}
			\rho_0(x)  &\quad x \in \Lambda'\\
			\mbox{Interpolation} &\quad x \in \mathcal T \\
			\eps \sqrt{1 + \lambda^2}  &\quad x = \psi_s^{-1}(\lambda, \mu) \in \mathcal N
		\end{cases}
	\end{equation}
where $\psi_s$ is the Legendrian normal coordinate chart for the neck region connecting $( \zeta U)^s(S_0)$ to $( \zeta U)^{s+1}(S_0)$.  Furthermore, one can assume that the function $\rho_\eps$ is invariant under the all symmetries satisfied by $\approxsol$.  
\end{defn}

The weighted Schauder norm is now defined as follows.  Set $\Lambda'' := \{ x \in \Lambda' : \rho_0(x) = 1 \}$.  Let $q$ be any tensor on $\approxsol$ and let $\mathcal O \subseteq \approxsol$ be any open subset.  Recall the notation 
\begin{equation*}
	\Vert q \Vert_{0, \mathcal O} = \sup_{x \in \mathcal O} \Vert q(x) \Vert \\
	\qquad \mbox{and} \qquad [q]_{\beta, \mathcal O} = \sup_{x,y \in \mathcal O} \frac{\Vert q(x) - \mathit{PT}(q(y)) 	\Vert}{\mathrm{dist}(x,y)^\beta} \, ,
\end{equation*}
where the norms and the distance function that appear are taken with respect to the induced metric of $\approxsol$, while $\mathit{PT}$ is the parallel transport operator from $x$ to $y$ with respect to this metric.  Now choose some large $r_0$ so that $x \in \Lambda''$ whenever $x$ has a distance larger than $r_0$ from any of the gluing points.  For any $0< r \leq 2 r_0$, define the annular region $A_r := \mathcal B_{r} \setminus \mathcal B_{r/2}$ as well as the norm on any subset $\mathcal O \subseteq A_r$
\begin{equation*}
	\vert f \vert_{l, \beta, \gamma, \mathcal O \cap \mathcal A_r} := \rho_\eps^{-\gamma}(r) \vert f \vert_{0, \mathcal O \cap \mathcal A_r} + \ \cdots + \rho_\eps^{-\gamma + l}(r) \Vert \nabla^l f \Vert_{0, \mathcal O \cap \mathcal A_r} + \rho_\eps^{-\gamma + l + \beta}(r) [ \nabla^l f ]_{\beta, \mathcal O \cap \mathcal A_r} \, .
\end{equation*}
Again, the derivatives which appear here are taken with respect to the induced metric of $\approxsol$.   Now make the following definition.

\begin{defn} 
	Let $\mathcal O \subset \approxsol$. The $C^{l, \beta}_\gamma$ norm on $\mathcal O \subseteq \approxsol$ is given by
\begin{equation}
	\mylabel{eqn:weightnorm}
	| f |_{C^{l, \beta}_\gamma (\mathcal O)} := \sum_{i=0}^l \| \nabla^i f \|_{0, \mathcal O \cap \Lambda''} + [ \nabla^l f ]_{\beta, \mathcal O \cap \Lambda''} + \sup_{r \in (0,2r_0]} | f |_{l, \beta, \gamma, \mathcal O \cap \mathcal A_r} \, . 
\end{equation}  
The notation for this norm will be abbreviated $| \cdot |_{C^{l, \beta}_\gamma}$ when there is no cause for confusion. 
\end{defn}

\subsection{Legendrian Deformations and the Contact-Stationary Legendrian Equation}
\label{subsec:def}

A second preliminary step is to have a concrete way of associating a Legendrian deformation $\phi_f$ of $\approxsol$ to every function $f$ in $X$.  The key is to construct the association `by hand' in such a way to have explicit control of $\phi_f$ in terms of $f$.  Broadly speaking, the deformation associated to $f$ is obtained as follows: one first extends $f$ to a neighbourhood of $\approxsol$ in a canonical way and multiplies by a suitable cut-off function; then one integrates the contact vector field generated by this extension up to time one.  

Here are the details of this construction along with the necessary estimates.  Observe first that the exponential map $\exp: N\approxsol \rightarrow \Sph^{2n+1}$ of the normal bundle of $\approxsol$ is a diffeomorphism in a tubular neighbourhood $\mathcal{V}$ of width $O(1)$ in $\mathcal T \cup \Lambda''$ and of width transitioning to $O(\eps)$ in the narrowest part of $\mathcal N$. Let $\eta$ be a smooth, monotone cut-off functions with support in $\mathcal V$ and define the extension operator $E_\eps : C^{4,\beta}(\approxsol) \rightarrow C^{4,\beta}(\Sph^{2n+1})$ by
$$E_\eps(f) (x) := 
\begin{cases}
	f(\exp^{-1}(x)) \cdot \eta \left( \dfrac{\mathrm{dist}(x,\approxsol)}{ \rho_\eps( \exp^{-1}(x))} \right) &\qquad x \in \mathcal{V} \\
	0 &\qquad x \not\in \mathcal{V}
\end{cases}$$
where $\rho_\eps$ is the weight function defined in Definition \ref{defn:weightfn}.  Consequently, $E_\eps(f)$ coincides with $f$ on $\approxsol$.  Now recall that the function $E_\eps(f)$ defines a contact vector field $X_{f}$ by means of the equations
\begin{equation*}
	\alpha(X_f) = E_\eps(f) \qquad \mbox{and} \qquad X_f \elbow \dif \alpha = \dif E_\eps(f) \big|_\Xi \, , 
\end{equation*}
and that integrating $X_f$ yields a one-parameter family of contactomorphisms $\phi_f^t : \Sph^{2n+1} \rightarrow \Sph^{2n+1}$ that is normal to $\approxsol$.  The time-one flow of this family is simply denoted $\phi_f$.

\begin{defn}
	\label{defn:contacto}
	The desired association of functions on $\approxsol$ to Legendrian deformations of $\approxsol$ is $f \mapsto E_\eps(f) \mapsto X_{f} \mapsto \phi_f$.  
\end{defn}

\subsection{The Strategy of the Proof}
\label{subsec:strategy}

The deformation procedure developed above can be fed into the abstract definition of $\Phi_{U, \zeta}$ and leads to the concrete partial differential equation that must be solved to find the contact-stationary Legendrian submanifold near $\approxsol$.  That is, a function $f \in \cfbg(\approxsol)$ must be found so that 
$$\Phi_{U, \zeta}(f) := \nabla \cdot \left( H_{\phi_f (\approxsol)} \elbow \dif \alpha \Big|_{\phi_f (\approxsol)} \right) = 0 \, $$  
where $\phi_f$ is the contactomorphism defined above. 

The strategy that will henceforth be used to apply the Banach space inverse function theorem to $\Phi_{U, \zeta}$ will take into account two important observations.  The first observation is that $\Phi_{U, \zeta}$ is equivariant with respect to the symmetries of $\approxsol$ induced by global isometries of $\Sph^{2n+1}$ preserving the Legendrian condition.  That is, if $\sigma \in O(2n+2)$ is such an isometry then $\Phi_{U, \zeta} (f \circ \sigma) = \Phi_{U, \zeta} (f) \circ \sigma$.  The second observation is that $\Phi_{U, \zeta}$ is not a full-rank mapping.  To see this, recall that the first variation formula \eqref{eqn:firstvar} applied to the generator of a contact isometry implies that every Legendrian submanifold $\Lambda$ satisfies 
$$\int_{\phi_f(\Lambda)} \Phi_{\Lambda}(f) \cdot q_H \big|_{\phi_f(\Lambda)} = 0$$
for every function $f$ and every Hermitian, harmonic, homogeneous polynomial $q_H$ generating the $U(n+1)$-rotations of $\Sph^{2n+1}$.  Consequently, the image of $\Phi_{U, \zeta}$ is constrained by $(n+1)^2$ conditions.  The correct interpretation of these conditions is to say that the graph $\{ (f, \Phi_{U, \zeta}(f)) : f \in X \}$ is contained in the Banach submanifold $\{ ( f , u ) : \int_{\phi_f(\approxsol)} u \cdot q_H \big|_{\phi_f(\approxsol)} = 0 \:\: \forall \:\: q_H \}$ of $Z \times Z$.  Therefore it suffices to show that the equation $\pi \circ \Phi_{U, \zeta}(f) = 0$ has a solution, where $\pi$ is the $L^2$-projection to the orthogonal complement of the subspace spanned by the functions $q_H \big|_{\approxsol}$.  Note that the linearization of $\pi \circ \Phi_{U, \zeta}$ at zero maps \emph{into} this orthogonal complement, and thus $\Dif \big( \pi \circ \Phi_{U, \zeta} \big) (0) = \mathcal L_{U, \zeta}$ holds. 

These observations suggest that the Banach space inverse function theorem should be applied in the following way.  In what follows, the subscript $\mathit{sym}$ denotes invariance with respect to the group of symmetries of $\approxsol$ from Section \ref{sec:additional}, so that if $u$ belongs to a $\mathit{sym}$-subscripted space then $u \circ \sigma = u$ for all such symmetries $\sigma$.  Also, the superscript $\perp$ denotes the $L^2$-orthogonal complement of a subspace of functions.  Let $\mathcal Q_0 := \mathrm{span}_\R \{ q_H \big|_{\approxsol} \}$ and define the following Banach spaces.

\begin{defn}
	\label{def:bspace}
	$X := \big[ C^{4, \beta}_{\gamma}(\approxsol) \big]_{\mathit{sym}}$ and $Z :=  \big[ C^{0, \beta}_{\gamma-4}(\approxsol) \cap \mathcal Q_0^\perp \big]_{\mathit{sym}} $.
\end{defn}

\noindent By the two observations above, the operator $\pi \circ \Phi_{U, \zeta} : X \rightarrow Z$ is well-defined and it is sufficient to find a solution of the equation $\pi \circ \Phi_{U, \zeta}(f) = 0$ in these spaces.  

Proving that a solution of this equation exists requires verifying that the conditions of the Banach space inverse function theorem hold for $\pi \circ \Phi_{U, \zeta}$.   In particular, one must construct an appropriately bounded right inverse for $\mathcal L_{U, \zeta}$.  Of course, as already mentioned, the primary obstruction to constructing a right inverse for $\mathcal L_{U, \zeta}$ with a small enough upper bound is the possibility that $Z$ is not sufficiently transverse to the approximate co-kernel of $\mathcal L_{U, \zeta}$ generated by the Jacobi fields.   However, an important part of the proof to follow is to show that the orthogonality and symmetry conditions built into $Z$ conspire to rule out these obstructions.

\section{The Linear Analysis}

\subsection{The Jacobi Fields of the Approximate Solution}
\label{subsec:apker}

It is necessary to identify a subspace of functions that approximates the invariant Jacobi fields of $\mathcal{L}_{U, \zeta}$ in $Z$.  As mentioned earlier, good approximating functions can be constructed as follows.  One takes the exact Jacobi fields of the linearized operator of one of the spherical constituents of $\approxsol$ and multiplies them by a cut-off function vanishing in a neighbourhood of the gluing points connecting this constituent to its neighbours.  Such functions can be understood as generating transformations which act by $U(n+1)$-rotation of the constituent of $\approxsol$ in question while leaving the others fixed.  Then one chooses all $\zeta U$-invariant linear combinations of these functions.  

More concretely, let $\chi_{\ext, r}$ be a smooth, monotone cut-off function that equals one in $\mathit{Pert}_{\eps}(S_0) \setminus [B_{r} (p_1) \cup B_{r} (p_2) ]$ for some $r$ and transitions to zero elsewhere.  Suppose $\chi_{\ext}$ is invariant under the symmetries of $\approxsol$. The functions just described are of the form
$$\tilde q :=  \left[ \sum_{s=0}^{kN-1} \chi_{\ext} \cdot q \big|_{\mathit{Pert}_\eps (S_0)} \circ (\zeta U)^{-s} \right] $$
where $q$ is any linear combination of the subset of polynomials in the list \eqref{eqn:ambientpolys} that do not vanish on $S_0$.  These functions, in coordinates $(x^1, x^2, \ldots, x^{n+1}) \in \R^{n+1}$ with $\sum_{s=1}^{n+1} (x^s)^2 = 1$, are
\begin{equation}
	\label{eqn:spherepolys}
	\begin{aligned}
		q_0 (x) &:= 1 \\
		q_1 (x) &:= (n+1) (x^1)^2 - 1 \\
		q_j(x) &:= (x^j)^2 - (x^{n+1})^2 \qquad \mbox{for $j = 2, \ldots, n$} \\
		q_{\textsc{im}, jk}(x) &:= 2 x^j x^k \qquad \mbox{for $1 \leq j < k \leq n+1$} \, .
	\end{aligned}
\end{equation}
Note that there are $D := \frac{1}{2} (n+1)(n+2)$ such functions. 

As mentioned earlier, one can not expect to find an appropriately bounded right inverse for $\mathcal L_{U, \zeta}$ on any subspace of functions that is not `sufficiently transverse' to the space of these invariant approximate Jacobi fields.  However, the principle that underlies the construction in this paper is that the orthogonality and symmetry conditions built into the definitions of $X$ and $Z$ guarantee this degree of transversality.  The way in which this principle will be used in the construction of the right inverse of $\mathcal L_{U, \zeta}$ is established in the following lemma, which shows that one has `approximate' orthogonality on each constituent of $\approxsol$.  

\begin{lemma}
	\label{lemma:orthog}
	Suppose $f : \approxsol \rightarrow \R$ is $L^2$-orthogonal to the functions $q_1, \ldots q_D$ and invariant with respect to all the symmetries of $\approxsol$.  Then 
	$$ \left| \int_{\approxsol} f  \cdot \tilde q_j \right| \leq C r^{\delta + n} |f|_{C^0_\delta}$$ 
	for all $j = 1, \ldots, D$.
\end{lemma}

\begin{proof}
This estimate is derived in two steps.  First, by invariance with respect to $\zeta U$, one finds
\begin{equation}
	\label{eqn:orthcalc}
	\left| \, \int_{\approxsol}  f \cdot \tilde q_j \right| \leq \left| \, \int_{\approxsol} f \cdot \frac{1}{kN} \sum_{s=0}^{kN-1}q_j \circ (\zeta U)^s \, \right| + C  r^{\delta + n} |f|_{C^0_\delta} \, .
\end{equation}
Now compute $\sum_{s=0}^{kN-1} q_j \circ (\zeta U)^{s}$ for each $q_j$.  This calculation depends upon the specific choice of $U$ made in Example 1 or in Example 2.  In the first example one finds:
\begin{equation}
	\label{eqn:nonzero1}
	\begin{aligned}
		\frac{1}{kN} \sum_{s=0}^{kN-1} q_0 \circ (\zeta U)^s (x) &= 1 \\
		\frac{1}{kN} \sum_{s=0}^{kN-1} q_{1} \circ (\zeta U)^s (x) &= (n+1) (x^1)^2 - 1 \\
		\frac{1}{kN} \sum_{s=0}^{kN-1} q_{j} \circ (\zeta U)^s (x) &= \frac{(x^j)^2 - (x^{n+1})^2}{N} \qquad \mbox{for $j = 2, \ldots, n+1$}
	\end{aligned}
\end{equation}
whereas all others vanish.  In the second example, it is better to use a slightly different basis for the space of functions spanned by $q_1, \ldots, q_D$.  That is, take the matrices $I, H_3, \ldots, H_{n+1}, H_{jk}$ and $H'_{jk}$ from \eqref{eqn:ambientpolys} and replace $H_1$ and $H_2$ by
\begin{equation*}
	H_1' = \mbox{ \footnotesize $ \left(
	\begin{array}{cc|ccc}
		n-1&&&&\\
		&n-1&&&\\
		\hline &&-2&&\\[-1ex]
		&&&\ddots&\\[-1ex]
		&&&&-2
	\end{array} \right) $ } \qquad \mbox{and} \qquad
	H_2' = \mbox{ \footnotesize $ \left(
	\begin{array}{cc|ccc}
		1&&&&\\
		&-1&&&\\
		\hline &&0&&\\[-1ex]
		&&&\ddots&\\[-1ex]
		&&&&0
	\end{array} \right) $} \, .
\end{equation*}
The Hermitian, harmonic, homogeneous polynomials corresponding to these matrices, restricted to $S_0$,  yield the functions  $q_1'(x) = (n+1)\big( (x^1)^2 + (x^2)^2 \big) - 2$ and $q_2'(x) = (x^1)^2 - (x^2)^2$.  One now finds
\begin{equation}
	\label{eqn:nonzero2}
	\begin{aligned}
		\frac{1}{kN} \sum_{s=0}^{kN-1} q_0 \circ (\zeta U)^s (x) &= 1 \\
		\frac{1}{kN} \sum_{s=0}^{kN-1} q_{1} \, \!\! ' \circ (\zeta U)^s (x) &= (n+1) \big( (x^1)^2 + (x^2)^2 \big) - 2  \\
		\frac{1}{kN} \sum_{s=0}^{kN-1} q_{j} \circ (\zeta U)^s (x) &= \frac{(x^j)^2 - (x^{n+1})^2}{N} \qquad \mbox{for $j = 3, \ldots, n+1$} \\
		\frac{1}{kN} \sum_{s=0}^{kN-1} q_{\textsc{im}, 12} \circ (\zeta U)^s (x) &= \frac{ x^1 x^2 }{N} \, ,
	\end{aligned}
\end{equation}
whereas all others vanish.  Consequently, the integral term in \eqref{eqn:orthcalc} vanishes for $q_j$ in the list \eqref{eqn:nonzero1} for Example 1 and for $q_j$ in the list \eqref{eqn:nonzero2} for Example 2, and the estimate holds for these $q_j$.

In each example, the lists \eqref{eqn:nonzero1} and \eqref{eqn:nonzero2} do not form a basis for the space of functions spanned by $q_1, \ldots, q_D$.  The remaining basis functions are 
$$q_{\textsc{im}, jk} = x^j x^k \qquad \mbox{for} \: 1 \leq j < k \leq n+1$$
in Example 1, and  
\begin{align*}
	q_{2} &= (x^1)^2 - (x^2)^2 \\
	q_{\textsc{im}, 1j} &= x^2 x^j  \qquad \mbox{for $3 \leq j  \leq n+1$}  \\
	q_{\textsc{im}, 2j} &= x^1 x^j  \qquad \mbox{for $3 \leq j  \leq n+1$} 
\end{align*}
in Example 2.   But for each of these functions, there is an additional symmetry $\sigma \in O(2n+2)$  under which $\approxsol$ is invariant and $q \circ \sigma = - q$, as can be readily verified from the definition of the additional symmetries in Section \ref{sec:additional}.  One can thus compute directly that 
$$\int_{\approxsol} f \cdot \tilde q_j = \int_{\approxsol} (f \cdot \tilde q_j) \circ \sigma = - \int_{\approxsol} f \cdot \tilde q_j = 0 \, ,$$
thereby verifying the estimate once again.  
\end{proof}

\subsection{The Linear Estimate}
\label{sec:linest}

This task at hand is to find an appropriately bounded right inverse for $\mathcal L_{U, \zeta}: X \rightarrow Z$ and a number of steps are needed to reach this goal.  The necessary steps become more complicated as the dimension $n$ decreases because of the nature of the indicial roots of the principal part of $\mathcal L_{U, \zeta}$, which is the bi-Laplacian operator, in low dimensions.  The material below presents the simplest proofs available in three cases: for dimensions $n\geq 5$, for dimensions $n=3,4$ and for dimension $2$.  This presentation has been chosen since it highlights the essential reason why the bounded right inverse exists, which is most obvious in the $n\geq5$ case.  The alternative --- finding a unified proof for all dimensions --- is somewhat possible, but would end up being just as long and would obscure the $n\geq 5$ case with unnecessary machinery.

The first of the step in the analysis leading to the bounded right inverse in all dimensions is a lemma showing that the lower order term $Q_{U, \zeta}:= 2\,  \nabla \cdot  B_ {\approxsol} (H_{\approxsol}, \nabla u )  - ( H_ {\approxsol} \cdot \nabla )^2 u $ in the linearized operator ${\mathcal L}_{U, \zeta} = \Delta_{U, \zeta} ( \Delta_{U, \zeta} + 2 (n+1)) + Q_{U, \zeta}$ as given in Corollary \ref{cor:linop} is uniformly small.

\begin{lemma}
	\label{lemma:loworder}
	The operator $Q_{U, \zeta} : X \rightarrow Z$ satisfies the estimate
	\begin{equation*}
		| Q_{U, \zeta}(u) |_{\cobg} \leq C\left(    \frac{\eps^n}{ r_\eps^{n-2}} + \eps^2 + r_\eps^4  \right)  | u |_{\cfbg}
	\end{equation*}
	for all $u \in X$.
\end{lemma}

\begin{proof}
	The operator $Q_{U, \zeta}$ satisfies the pointwise estimate
	\begin{equation*}
		\bigl| [\rho_\eps(x)]^{4-\gamma}  Q_{U, \zeta}(u)(x) \bigr| \leq C | \rho_\eps(x)|^2 \, \bigl| H_{\approxsol} (x) \bigr| \left(   \bigl\| B_{\approxsol}(x) \bigr\| + \bigl\| \rho_\eps(x)  \nabla B_{\approxsol}(x) \bigr\| \right)   |u |_{\cfbg} 
	\end{equation*}
for any $x \in \approxsol$.  If $x$ belongs to the exterior region $\Lambda' \subseteq \approxsol$, Proposition \ref{prop:exest} gives the estimate 
$$| Q_{U, \zeta}(u) |_{\cobg(\Lambda')} \leq C\eps^{2n} r_\eps^{-2n+2} |u|_{\cfbg} \, .$$
If $x$ belongs to the neck region $\mathcal N \subseteq \approxsol$ where $\rho_\eps(\lambda) = \eps \sqrt{1 + \lambda^2} :=  \rho(\lambda)$ Proposition \ref{prop:neckest} gives 
\begin{align*}
	| Q_{U, \zeta}(u) |_{\cobg(\mathcal N)} &\leq \sup_{ | \lambda| \leq r_\eps/\eps} C \, \eps^2 [ \rho(\lambda)]^2 \left( \frac{\eps}{[\rho(\lambda)]^{n-1}} + \eps \rho(\lambda) \right) \left(\frac{1}{\eps [\rho(\lambda)]^{n+1}} + \eps \rho(\lambda) \right) |u|_{\cfbg} \\
	&\leq C \big( \eps^2 + r_\eps^4 \big)  \, . 
\end{align*}
Finally, if $x$ belongs to the transition region $\mathcal T \subseteq \approxsol$, Proposition \ref{prop:transmeancurvest} gives 
\begin{align*}
	| Q_{U, \zeta}(u) |_{\cobg(\mathcal T)} &\leq C \, r_\eps^2 \left( \frac{\eps^n}{r_\eps^{n-1}} + r_\eps \right) \left(\frac{ \eps^n}{r_\eps^{n+1}} +  r_\eps \right) \\
	&\leq  C \left(\frac{\eps^n}{r_\eps^{n-2}} + r_\eps^4 \right) 
\end{align*}
using the fact that $\eps < r_\eps$.  Consolidating these separate results yields the desired estimate for the supremum norm.  The H\"older norm estimate is similar.
\end{proof}

Henceforth, set $L_{U, \zeta} := \Delta_{U, \zeta} ( \Delta_{U, \zeta} + 2(n+1) )$.  The existence of an appropriately bounded right inverse in dimensions $n \geq 5$ is already a consequence of Lemma \ref{lemma:loworder}.  It uses a fairly straightforward contradiction argument together with the fact that $\mathcal L_{U, \zeta}$ and  $L_{U, \zeta}$ are self-adjoint.

\begin{prop}
	\label{prop:linfive}
	Suppose $n \geq 5$ and choose $\gamma \in (4-n, 0)$.  If $\eps$ is sufficiently small then there is a constant $C$ independent of $\eps$ so that the operator $\mathcal  L_{U, \zeta} : X \cap \mathcal Q_0^\perp \rightarrow Z$ is bijective and satisfies $|u|_{\cfbg} \leq C | \mathcal L_{U, \zeta} (u) |_{\cobg}$ for all $u \in X$ where $C$ is a constant independent of $\eps$.
\end{prop}

\begin{proof}
The strategy of this proof is as follows.  It will be shown below that $L_{U,\zeta}$ is injective on $X \cap \mathcal Q_0^\perp$ and satisfies $|u|_{\cfbg} \leq C_0 | L_{U, \zeta} (u) |_{\cobg}$ for some constant $C_0$ independent of $\eps$.  Since $L_{U,\zeta}$ is self-adjoint then $L_{U,\zeta} : X \cap \mathcal Q_0^\perp \rightarrow Z$ is surjective.  The estimate of Lemma \ref{lemma:loworder} then shows that the operator $\mathcal L_{U, \zeta} $ is injective on $X \cap \mathcal Q_0^\perp$ as well when $\eps$ is sufficiently small, and satisfies the desired bound with $C := 2 C_0$.  Therefore $\mathcal L_{U, \zeta} $ is surjective, again by self-adjointness.

	The injectivity bound for $L_{U,\zeta}$ can be found in this way.  By standard elliptic theory and scaling arguments, there is already some constant $C$ independent of $\eps$ so that 
	\begin{equation}
		\label{eqn:ellest}
		|u|_{\cfbg} \leq C \big( | { L}_{U, \zeta} (u) |_{\cobg} + |u|_{C^{0,\beta}_{\gamma}} \big)
	\end{equation}
for all $u \in \cfbg(\approxsol)$.  Hence it is enough to obtain a contradiction from the assumption that the estimate $|u|_{C^{0,\beta}_{\gamma}} \leq C  | { L}_{U, \zeta} (u) |_{\cobg}$ is false.  In other words, suppose that there are sequences of: 
\begin{itemize}
	\item scale parameters $\eps_i \rightarrow 0$ and corresponding complex numbers $\zeta_i \rightarrow 1$
	
	\item approximate solutions $\Lambda_i := \tilde \Lambda_{U, \zeta_i}$
	
	\item  linear operators $L_i := L_{U, \zeta_i}$
	
	\item  weight functions $\rho_i  := \rho_{\eps_i}$ 
	
	\item  Banach spaces $X_i := \big[ C^{4, \beta}_{\gamma}(\Lambda_i) \cap {\mathcal Q}_0^\perp \big]_{\mathit{sym}}$ and $Z_i :=  \big[ C^{0, \beta}_{\gamma-4}(\Lambda_i) \cap \mathcal Q_0^\perp \big]_{\mathit{sym}}$
	
	\item functions $u_i \in X_i$ satisfying $|u_i|_{C^{0,\beta}_\gamma} = 1$ and $\displaystyle \lim_{i \rightarrow \infty} | L_i (u_i) |_{\cobg} = 0$.  
\end{itemize}
Now let $p_i \in \Lambda_i$ be a point where $u_i(p_i) [\rho_i(p_i)]^{-\gamma} = 1$. There are two cases to consider.
	
\medskip \noindent \scshape Case 1: \upshape The non-concentrating case. \medskip

Let $S_R := S_0 \setminus [ B_R(p_1) \cup B_R(p_2)]$ and $S_0^\ast := S_0 \setminus \{ p_1, p_2 \}$ where $p_1, p_2$ are the gluing points of $S_0$.  Suppose without loss of generality that $p_i$ converges to some point $p \in S_R$ for some non-zero radius $R$.   By the estimate \eqref{eqn:ellest} and the Arzela-Ascoli theorem, there is a subsequence of $u_i$ converging to a non-zero function $u_R \in C^{4, \beta}(S_R)$.  It is now possible to pass to a further subsequence and obtain convergence to a non-zero function $u^\ast \in C^{4, \beta}_{\gamma}(S_0^\ast)$.   This function satisfies $\Delta_{S} (\Delta_{S}  + 2(n+1) ) (u^\ast) = 0$ on $S_0^\ast$, where $\Delta_{S}$ is the Laplacian of the standard sphere.  Since $\gamma-2 \in (2-n, 0)$, then $(\Delta_{S}  + 2(n+1) ) (u^\ast)$ is actually smooth and thus must be constant.  Since $\gamma \in (4-n, 0) \subseteq (2-n, 0)$, then $u^\ast$ is also smooth and consequently must be a linear combination of  the functions in the list \eqref{eqn:spherepolys}.  But the orthogonality and symmetry conditions that $u^\ast$ must satisfy as a limit of functions $u_i \in X_i$ rules this out, as can be seen by allowing $\eps_i \rightarrow 0$ in Lemma \ref{lemma:orthog}.

\medskip \noindent \scshape Case 2: \upshape The case where concentration occurs. \medskip

Suppose without loss of generality that $p_i$ converges to the gluing point $p_1 \in S_0$.  Consequently $p_i$ eventually enters a Legendrian normal neighbourhood containing the neck region connecting $S_0$ to its neighbour at $p_1$.  One can thus use the Legendrian normal coordinates and the Legendrian Lawlor embedding to write $p_i$ as a point $(\lambda_i, \mu_i) \in \R \times \Sph^{n-1}$ and to consider $u_i \in \cfbg(\R \times \Sph^{n-1})$.  The norm on this space is a weighted norm, where the weight function is $\rho(\lambda) := \eps \sqrt{1 + \lambda^2}$.  Up to subsequences, there are now two possibilities: either $(\lambda_i, \mu_i) \rightarrow  (\lambda, \mu)$ or else $\lambda_i \rightarrow \infty$.  Each of these possibilities will be ruled out in turn.  In the first of these possibilities, the estimate \eqref{eqn:ellest} and the Arzela-Ascoli theorem on the neck itself imply that $u_i$ has a subsequence converging to a non-zero function $u \in C^{4,\beta}_{\gamma'}(\R \times \Sph^{n-1})$ for any $\gamma' > \gamma$ and satisfying the equation $\Delta_N (\Delta_N (u)) = 0$ where $\Delta_N$ is the Laplacian of the Lawlor neck.  Since $\gamma \in (4-n, 0)$, this solution is decaying at infinity.  But given the form of the Lawlor neck metric \eqref{eqn:lawlormetric} on  $\R \times \Sph^{n-1}$, there can be no decaying solutions of this equation.  In the second possibility, one can consider $u_i$ on a ball of constant radius about $\lambda_i$ and re-scale by a factor of $\eps_i^{-1}$ to obtain a sequence of functions which possesses a subsequence converging to a non-zero, decaying solution of $\mathring \Delta (\mathring \Delta (u)) = 0$ on $\R^n$ with respect to the Euclidean metric.  Again, there are no decaying solutions of this equation.
\end{proof}

Because of the overlapping nature of the indicial roots of the bi-Laplacian, the technique above is not available in the lower dimensions $n=4, 3, 2$.  In these cases, right inverses will have to be constructed `by hand' and the required estimates derived as part of the construction.  Two lemmata are needed.  The first shows that it is again only necessary to work with the simpler operator ${ L}_{U, \zeta}$.  Denote the operator norm of a map between Banach spaces by  $\| \cdot \|_{\mathit{op}}$.

\begin{lemma}
	\label{lemma:rinv}
	Suppose $ L_{U, \zeta} : X \rightarrow Z$ possesses a right inverse $\mathcal R : Z \rightarrow X$ satisfying the bound $|\mathcal R (f) |_{\cfbg} \leq C |f|_{\cobg}$ where $C := C(\eps)$ is a constant depending in some way on $\eps$.   If  
	$$\displaystyle \limsup_{\eps \rightarrow 0} \, C(\eps)  \| Q_{U, \zeta} \|_{\mathit{op}} < \frac{1}{2}$$
	and $\eps$ is sufficiently small, then $\mathcal L_{U, \zeta}$ also possesses a right inverse satisfying the same sort of bound, except with $C$ replaced by $2 C$.
\end{lemma}

\begin{proof}
	Propose $u := \mathcal R (f+g)$ as the solution of the equation $\mathcal L_{U, \zeta} (u) = f$.  Since $\mathcal L_{U, \zeta} =  L_{U, \zeta} + Q_{U, \zeta}$ then $g$ satisfies $g = - Q_{U, \zeta} \circ \mathcal R (f+g)$.  The mapping $T_f : Z \rightarrow Z$ given by $T_f(g) := - Q_{U, \zeta} \circ \mathcal R (f+g)$ satisfies both 
	$$|T_f(g_1) - T_f(g_2)|_{\cobg} \leq C  \| Q_{U, \zeta} \|_{\mathit{op}} |g_1 - g_2|_{\cobg} \quad \mbox{and} \quad |T_f(g) |_{\cobg} \leq C  \| Q_{U, \zeta} \|_{\mathit{op}} | f + g |_{\cobg}$$
	so that $T_f$ is a contraction mapping on $B_R(0)$ when $|f|_{\cobg} < R$ and $\eps$ is sufficiently small.  Therefore there is a unique fixed point $g_f$. The desired right inverse is the map $f \mapsto \mathcal R(f + g_f)$, for which the desired bound is clear.
\end{proof}

\noindent The second lemma describes the kernel of the Laplacian on the Legendrian Lawlor neck.  Note that the non-constant, bounded harmonic function in dimension $n=2$ has logarithmic growth.  This will necessitate a different proof of the construction of the right inverse in the $n=2$ case.

\begin{lemma}
	\label{lemma:laplacekernel}
	The following facts about the kernel of the Laplacian of the Lawlor neck are true.
	\begin{itemize}
		\item The harmonic functions of sub-linear growth are spanned by the constant function $1$ and the function $S_A(\lambda)$ given in equation \eqref{eqn:Sfunction}.
		
		\item Suppose $f = a_i^+ \mu^i | \lambda |  + O(|\lambda|^{-1})$ for $\lambda \gg 1$ and $f = a_i^- \mu^i | \lambda |  + O(|\lambda|^{-1})$ for $\lambda \ll -1$ in the coordinates of the neck, where $a_i^\pm \in \R$.  Then there is a linearly growing harmonic function which equals $f$ up to $O(|\lambda|^{-1})$ on the ends of the neck.
	\end{itemize}
\end{lemma}

\begin{proof}
	By separation of variables, one expects to find $2$ linearly independent harmonic functions depending only on $\lambda$, as well as $2n$ linearly independent and linearly growing harmonic functions.  The constant function $1$ is clearly an example of the first kind.   The procedure of Corollary \ref{cor:jacorigin}, or at least its Lagrangian analogue in $\C^n$ can be used to find the remaining functions.  (That is, by the computation of Corollary \ref{cor:jacorigin} a Hamiltonian function that generates a symplectic isometry and then restricted to a minimal Lagrangian submanifold gives a harmonic function on the submanifold.)  The function $S_A$ arises in this way by considering the infinitesimal generator of dilations.  The linearly growing functions arise in this way by considering the infinitesimal generators of translations parallel to one or the other of the asymptotic $n$-planes of the neck.  Since these $n$-planes are transverse, the $2n$ resulting functions are linearly independent and form a basis for linear growing functions on each end of the neck separately.  
\end{proof}

The estimates for the lower dimensions $n = 4, 3, 2$ can now be derived.  Before beginning the proofs of the estimates, define the following smooth, monotone cut-off functions:
\begin{align*}
	\chi_{\ext, r} &:= 
	\begin{cases}
		1 &\qquad x \in \mathcal B_r^c \\
		\mbox{Interpolation} &\qquad x \in \mathcal B_r \setminus \mathcal B_{r/2} \\
		0 &\qquad x \in \mathcal B_{r/2}
	\end{cases} \\[1.5ex]
	\chi_{\neck, r} &:= 1 - \chi_{\ext, r} 
\end{align*}
where $\mathcal B_r$ is the disjoint union of balls of radius $r$ containing the neck regions of $\approxsol$ as defined in equation \eqref{eqn:balls}.  One can assume that these functions are invariant with respect to all the symmetries satisfied by $\approxsol$.

\begin{prop}
	\label{prop:linfour}
	Suppose $n = 4$ and choose $\gamma \in (0,1)$.  Then the operator $L_{U, \zeta} : X \rightarrow Z$ possesses a right inverse $\mathcal R : Z \rightarrow X$ satisfying the bound $| \mathcal R(f) |_{\cfbg} \leq C  |f|_{\cobg}$ where $C$ is a constant independent of $\eps$.
\end{prop}

\begin{proof}
Using a contradiction argument as in the proof of Proposition \ref{prop:linfive}, one can show that $\Delta_{U, \zeta} : [ C^{2, \beta}_{\gamma - 2}(\approxsol) \cap \mathcal Q_0 ]_{\mathit{sym}} \rightarrow Z$ is bijective and satisfies the estimate $|w|_{\ctbg} \leq C | \Delta_{U, \zeta} (w) |_{\cobg}$ for all $w \in [ C^{2, \beta}_{\gamma - 2}(\approxsol) \cap \mathcal Q_0 ]_{\mathit{sym}}$ where $C$ is a constant independent of $\eps$.  The reason the same argument works is because $\gamma - 2 \in (-2, -1)$ is once again contained in the range of weights for which the Laplacian on $S_0 \setminus \{p_1, p_2\}$ and on $N$ is bijective.  Moreover, the self-adjointness of $\Delta_{U, \zeta}$ implies that the orthogonality conditions satisfied by $f$ are inherited by $w := \Delta^{-1}_{U, \zeta} (f)$.  Therefore it remains only to solve the equation $(\Delta_{U, \zeta} + 2(n+1) )(u) = w$ for $u \in \cfbg(\approxsol)$ and $w \in [ C^{2, \beta}_{\gamma - 2}(\approxsol) \cap \mathcal Q_0 ]_{\mathit{sym}}$ for $\gamma \in (0,1)$. 

The technique that will be used to find a solution of $\big( \Delta_{U, \zeta} + 2(n+1) \big) (u_{\ext}) = 0$ will be to patch together local solutions on the exterior regions with local solutions on the neck regions of $\approxsol$ to construct a good approximate solution.  This process will then be iterated to yield an exact solution. To begin the patching process, choose four radii satisfying $0 < r_3 \ll r_4 < r_2 < r_1$ in such a way that the supports of $\nabla \chi_{\ext, r_j}$ do not overlap.   It will also be necessary to choose $r_3 = O(\eps)$ and $r_1, r_2, r_4$ small \emph{a priori} but independent of $\eps$.  

\paragraph{Step 1: The exterior regions.} Define the functions $\hat q_i := \sum_{s=0}^{kN-1} q_i \big|_{S_0} \circ (\zeta U)^s$ as well as $\hat q_i^{\mathit{approx}} := \hat q_i \chi_{\ext, r_1}$ and
$$w_{\ext} := w \chi_{\ext, r_1} -   \sum_{i,j = 1}^D M^{ij} \left(  \int_{S_0} w \chi_{\ext, r_1} \hat q_i \right) \hat q_j^{\mathit{approx}} $$
where $M^{ij}$ is the inverse of the matrix whose coefficients are $\int_{S_0} \hat  q_i \hat q_j^{\mathit{approx}}$.  This is the $L^2$-projection of $w_{\ext}$ to the orthogonal complement of the space of functions spanned by $ \hat q_1, \ldots, \hat q_D$.  Furthermore, $w_{\ext} \big|_{S_0}$ is a smooth function of compact support on $S_0^\ast$ which is, by the calculations of Lemma \ref{lemma:orthog}, orthogonal to $q_1 \big|_{S_0}, \ldots, q_D \big|_{S_0}$.   The estimate of Lemma \ref{lemma:orthog} implies $| w_{\ext} |_{\ctbg} \leq C |w|_{\ctbg}$ with $C$ independent of $\eps$ because $\gamma \in (0,1)$.  By $\zeta U$-invariance, it is necessary only to find the solution $u_{\ext} \big|_{S_0}$ and then extend this solution by $\zeta U$-invariance to all of $\approxsol$.  By the theory of the Laplace operator on $S_0$ in smooth spaces, there exists a solution $u_{\ext} \in C^{4,\beta}(S_0)$ of the equation $(\Delta_S + 2(n+1))(u) = w_{\ext}$.  By considering the Taylor expansion of $u_{\ext}$ at the points $p_1$ and $p_2$, one can write
\begin{equation}
	\label{eqn:modifone}
	u_{\ext} := v_{\ext} + a \big(  \eta^1_{\neck,r_2} +  \eta^2_{\neck,r_2} \big)
\end{equation}
where $a := u_{\ext}(p_1) = u_{\ext}(p_2)$ (which are the same by symmetry) and $v_{\ext} \in \cfbg(S_0 \setminus \{p_1, p_2\})$ (since $\gamma \in (0,1)$ and $v_{\ext}$ grows linearly with distance from each $p_i$).  Here $  \eta^i_{\neck,r_2} $ is a smooth, monotone cut-off function that equals one near $p_i$ and transitions to zero elsewhere, and equals $\chi_{\neck, r_2}$ on the region of overlap.  Moreover, the estimate
$$ | v_{\ext}|_{\cfbg(S_0 \setminus \{p_1, p_2\})} + | a |  \leq C | w_{\ext} |_{C^{2,\beta} (S_0 \setminus \{p_1, p_2\})} \leq C r_1^{\gamma - 2}  |w_{\ext}|_{\ctbg} $$
is valid for some constant $C$ independent of $\eps$.  

At present, the function $u_{\ext}$ can be viewed as a function defined on $\approxsol \setminus \mathcal B_r$ for some small $r$.  To extend $u_{\ext}$ to a function defined on all of $\approxsol$, first subtract $a J$ where $J$ is the function in the list \eqref{eqn:nonzero1} in Example 1 and or in the list \eqref{eqn:nonzero2} in Example 2 which is even with respect to the symmetry exchanging $p_1$ and $p_2$ and has $J(p_i) = 1$.  Note that on can write $J = ( \eta^1_{\neck,r_2} +  \eta^2_{\neck,r_2}) + \tilde J$ where $\tilde J \in \cfbg(S_0)$.   Then set
$$\bar u_{\ext} := \chi_{\ext, r_2} \big( v_{\ext} - a  \tilde J \big) \, .$$
Finally, one deduces the estimate
$$|\bar u_\ext|_{\cfbg} \leq C  |w|_{\ctbg}$$
where the constant $C$ is independent of $\eps$.

\paragraph{Step 2.  The neck regions.}  Define the function 
$$w_{\neck} := \chi_{\neck, r_3}\,  \big ( \, w -( \Delta_{U, \zeta} +2(n+1))(\bar u_{\ext})\,  \big)$$
and consider $w_{\neck}$ as a function of compact support defined on the $\eps$-scaled Lawlor neck $N_\eps$.  One now looks for a solution of the equation $\Delta_N (u) = w_{\neck}$ where $\Delta_N$ is the Laplacian of $N_\eps$.  Set $\gamma' := -2-\gamma$.  Since $w_{\neck}$ has compact support, then $w_{\neck} \in C^{2,\beta}_{\gamma' -2}(N_\eps)$ where the norm in this case is the weighted H\"older norm on $N_\eps$ with the weight function $\rho_\eps(\lambda) := \eps \sqrt{1+\lambda^2}$ given in the coordinates of the Lawlor neck.  By the theory of the Laplace operator on asymptotically flat manifolds, the operator
$$\Delta_N : C^{4,\beta}_{\gamma'}(N_\eps) \oplus \mathcal D \rightarrow C^{2,\beta}_{\gamma' - 2}(N_\eps)$$
is surjective with a two-dimensional kernel.  Here $\mathcal D$ is the \emph{deficiency space} given  by
$$\mathcal D := \mathrm{span}_{\R} \left\{ \frac{\chi_{\ext, r_3}}{| \lambda |^2} , \chi_{\ext, r_3}\right\} \, .$$
The reason for this is because the range $(\gamma', \gamma)$ contains rates of growth of the indicial roots of the Laplacian corresponding to constant functions and functions which decay like the Green's function at infinity, and the symmetry conditions require $f$ to be even with respect to the transformation $\lambda \mapsto - \lambda$.  Therefore one can find a solution in this space, and by adding a constant this solution has the form
\begin{equation}
	\label{eqn:modiftwo}
	u_{\neck} := v_{\neck} +  b \,\frac{ \chi_{\ext, r_3} }{| \lambda |^2} 
\end{equation}
where $v_{\neck} \in C^{4,\beta}_{\gamma'}(N_\eps)$ and satisfies the estimates
\begin{equation}
	\label{eqn:linone}
	|v_{\neck}|_{C^{4,\beta}_{\gamma'}(N_\eps)} +  \eps^{2 + \gamma} |b|  \leq C r_3^{2 + 2 \gamma} |w_{\neck}|_{\ctbg}
\end{equation}
for some constant $C$ independent of $\eps$.  This is because $w_{neck}$ is supported in the region of radius $O(r_3)$.   Finally, extend $u_{\neck}$ to all of $\approxsol$ simply by setting 
$$\bar u_{\neck} := \chi_{\neck, r_4} u_{\neck}$$
and invoking invariance with respect to the group generated by composition with $\zeta U$.  Finally, one deduces the estimate
\begin{equation*}
	\begin{aligned}
		| \bar u_{\neck} |_{\cfbg} & \leq |v_{\neck}|_{\cfbg(N_\eps)} +   \sup_{r_3/\eps \leq |\lambda | \leq r_4/\eps} [\rho_{\eps} (\lambda)]^{-\gamma} \, | b | \,|\lambda|^{-2}  \\[1ex]
		& \leq C  \left[ \left( \frac{r_3}{\eps} \right)^{2+ 2 \gamma}   +  \left( \frac{r_3}{\eps} \right)^{\gamma}\right] |w_{\neck}|_{\ctbg} \\
		&\leq C |w|_{\ctbg}
	\end{aligned}
\end{equation*}
using equation \eqref{eqn:linone} and the estimate for $|\bar u_{\ext}|_{\cfbg}$ given in Step 1 along with the values of the radii $r_i$ and the estimate for $| w_{\ext} |_{\ctbg}$ in terms of $|w|_{\ctbg}$.

\paragraph{Step 3.  Estimates and convergence.}

Set $\bar u := \bar u_{\ext} + \bar u_{\neck}$.  The estimates from Step 1 and Step 2 show that $| \bar u |_{\cfbg} \leq C  |w|_{\ctbg}$.   The function $\bar u$ should be seen as the approximate solution of the equation $(\Delta_{U,\zeta} + 2(n+1))(u) = w$ satisfying the appropriate estimate.  To justify this, it must be shown that $| (\Delta_{U,\zeta} + 2(n+1))(\bar u) - w |_{\ctbg} \leq \frac{1}{2} | w |_{\ctbg}$ when $\eps$ is sufficiently small.  If this holds, then the procedure of Step 1 and Step 2 can be iterated to yield an exact solution in the limit and satisfying the appropriate estimate.

A straightforward calculation shows that
\begin{equation}
	\label{eqn:lintwo}
	\begin{aligned}
		(\Delta_{U,\zeta} + 2(n+1))(\bar u) - w &= 2(n+1) \bar u_{\neck} \\
		&\qquad 	+ (\Delta_{U,\zeta} - \Delta_N) (\bar u _{\neck})  + \chi_{\ext, r_3} ( \Delta_{U, \zeta} - \Delta_S) (\bar u_{\ext}) \\
		&\qquad + [\Delta_N, \chi_{\neck, r_4}](u_{\neck})   \\
		&\qquad +  \chi_{\ext, r_2} ( w - w_{\ext} )  \, .
	\end{aligned}
\end{equation}
where $[\Delta, \chi](u)$ is notation for $\Delta(\chi u) - \chi \Delta (u)$.  The various terms in \eqref{eqn:lintwo} are all small multiples of $| w |_{\ctbg}$ when measured with respect to the $\ctbg$ norm, but for different reasons. 
\begin{itemize}
	\item The term on the first line is supported in the region $r < r_4$.  Thus its $\ctbg$ norm acquires the factor $r_4^{2}$ so that $2(n+1) | \bar u_{\neck}|_{\ctbg} \leq C r_4^2  |w|_{\ctbg}$.  This is can be made small with a small enough initial choice of $r_4$.    
	
	\item In the terms on the second line, the operators $\Delta_{U, \zeta} - \Delta_S$ and $\Delta_{U, \zeta} - \Delta_N$ are uniformly small in the regions where their arguments are supported, as can be deduced from Proposition \ref{prop:exest}, Proposition \ref{prop:meancurvcomparison} and Proposition \ref{prop:neckest}.  
	
	\item In the term on the third line, the function $u_{\neck}$ has been engineered to have strong decay in the support of $\nabla \chi_{\neck, r_4}$.  Indeed, the estimate \eqref{eqn:linone} implies that  $|v_{\neck} |_{\ctbg ( \mathcal B_{r_4} \setminus \mathcal B_{r_4/2})} \leq C(r_3/r_4)^{ 2 + 2\gamma}$,  Since $r_3 / r_4 = O(\eps)$, this quantity can be made as small as desired by choosing $\eps$ and the $r_i$ small enough.  A similar analysis holds for the $b \chi_{\ext, r_3} |\lambda|^{-2}$ term in $\bar u_{\neck}$.  
	
	\item The remaining term is handled by Lemma \ref{lemma:orthog} and becomes small when $r_1$ is small enough.
	
\end{itemize}

The arguments of the preceding paragraphs show that $| (\Delta_{U,\zeta} + 2(n+1))(\bar u) - w |_{\ctbg} \leq \frac{1}{2} | w |_{\ctbg}$ provided that $\eps$ is chosen sufficiently small.  As indicated above, this is enough to construct the right inverse that satisfies the desired bound.  
\end{proof}

\begin{prop}
	\label{prop:linthree}
	Suppose $n = 3$ and choose $\gamma \in (1,2)$.  Then the operator $L_{U, \zeta} : X \rightarrow Z$ possesses a right inverse $\mathcal R : Z \rightarrow X$ satisfying the bound $| \mathcal R(f) |_{\cfbg} \leq C  |f|_{\cobg}$ where $C$ is a constant independent of $\eps$.
\end{prop}

\begin{proof}
The strategy of this proof is similar to that of the previous proof.  In fact, the proof of the existence of the right inverse of $\Delta_{U, \zeta} : \ctbg(\approxsol) \rightarrow \cobg(\approxsol)$, bounded above by a constant independent of $\eps$, still follows from an argument by contradiction as before because $\gamma - 2 \in (-1,0)$ is in the correct range.  What is different is that now the construction of the right inverse for  $\Delta_{U, \zeta} + 2(n+1) : \cfbg(\approxsol) \rightarrow \ctbg(\approxsol)$ has to be modified to take into account the range $\gamma \in (1, 2)$. Indeed, equation \eqref{eqn:modifone} must be modified by expanding $u_{\ext}$ up to order two and cancelling the linear term in the expansion with the linearly growing functions in the kernel of $\Delta_S + 2(n+1)$ on $S_0$ given in the lists \eqref{eqn:nonzero1} or \eqref{eqn:nonzero2}.  (One can check that under the symmetry conditions on $u_{\ext}$ there are enough such functions to accomplish this.)  Then the solution procedure on the neck must be modified by enlarging the deficiency space to include the functions with linear growth on the neck as well.  By Lemma \ref{lemma:laplacekernel}  one can add a linear combination of these functions to eliminate the linearly growing term in the solution on the neck.   This leads to equation \eqref{eqn:modiftwo} but with $O( |\lambda|^{-1} )$ terms.  The estimate of the approximate solution and the convergence of the iteration leading to the exact solution are unchanged.
\end{proof}

The final estimate in the $n=2$ case is more complicated still.   This is because there is no range of weights for which the Laplacian on punctured manifolds is bijective when $n=2$ so that one can no longer use a contradiction argument to derive the existence of an appropriately bounded right inverse for $\Delta_{U, \zeta} : [C^{2,\beta}_{\gamma - 2} (\approxsol) ]_{\mathit{sym}} \rightarrow Z$.  The patching technique of the previous two propositions must thus be used for $\Delta_{U, \zeta}$ as well as for $\Delta_{U, \zeta} + 2 (n+1)$.  And both of these tasks are complicated by the fact that the odd harmonic function depending only on the neck coordinate $\lambda$ has logarithmic growth.  This fact makes it necessary to exploit the symmetries of $\approxsol$ even further.

\begin{prop}
	\label{prop:lintwo}
	Suppose $n = 2$ and choose $\gamma \in (1,2)$.  Then the operator $L_{U, \zeta} : X \rightarrow Z$ possesses a right inverse $\mathcal R : Z \rightarrow X$ satisfying the bound $| \mathcal R(f) |_{\cfbg} \leq C \eps^{\gamma - 2} |f|_{\cobg}$ where $C$ is a constant independent of $\eps$. 
\end{prop}

\begin{proof}
  	The proof has two parts: the construction of a right inverse for $\Delta_{U, \zeta} : [C^{2,\beta}_{\delta} (\approxsol) ]_{\mathit{sym}} \rightarrow [C^{0,\beta}_{\delta-2} (\approxsol) ]_{\mathit{sym}}$ with $\delta \in (-1,0)$ and one for $\Delta_{U, \zeta} + 6 : [C^{4,\beta}_{\gamma} (\approxsol) ]_{\mathit{sym}} \rightarrow [C^{2,\beta}_{\gamma - 2} (\approxsol) ]_{\mathit{sym}}$ with $\gamma \in (1,2)$.  Begin with the first of these constructions, in which the equation $\Delta_{U, \zeta} (w) = f$ for  $f \in [C^{0,\beta}_{\delta-2} (\approxsol) ]_{\mathit{sym}}$ will be solved according to the following four steps.   The condition $\int_{\approxsol} f = 0$ which is built into the space $[C^{0,\beta}_{\delta-2} (\approxsol) ]_{\mathit{sym}}$ will be used in a critical way below.
	
\newcommand{\trans}{\mathit{trans}}	
	
\paragraph{Step 0. Re-balancing the mass of \boldmath $f$.} It will first be shown that one can find a solution of the equation $\Delta_{U, \zeta} (w) = \bar f_{\trans}$ where $\bar f_{\trans}$ is some function that coincides with $f$ in an annular region around each of the gluing points.  Then one will be able to complete the construction of the right inverse by finding local solutions on the neck and on the exterior regions of the equation $\Delta_{U, \zeta} (w) = f - \bar f_{\trans}$, which has the advantage that $f - \bar f_{\trans}$ vanishes in the annular region.  The usefulness of this will become apparent in Step 3, where it provides the additional precision required to match the local solutions properly.

To begin, consider the gluing points $p_1, p_2 \in S_0$. For any pair of radii $0< r_0 < r_4$ let $\chi_{\trans, r_0, r_4}$ be a smooth, monotone cut-off function that equals one in $\mathit{Pert}_\eps (S_0) \cap \bigl[ \bigcup_k B_{r_4}(p_k) \setminus B_{ r_0}(p_k) \bigr]$ and transitions to zero in $\mathit{Pert}_\eps (S_0) \cap \bigl[ \bigcup_k B_{2 r_4}(p_k) \setminus B_{ r_0/2}(p_k) \bigr]$.  Fix two such radii and define 
$$f_{\trans} := f \chi_{\trans, r_0, r_4} \, .$$
An approximate solution of the equation $\Delta_{U, \zeta} (u) = f_{\trans}$ can be found as follows.  View $f_{\trans}$ as a symmetric function defined on $S_0^\ast$ carrying a small perturbation of the standard metric, and consider the equation $\Delta_S (u) = f_{\trans}$ on $S_0^\ast$.  By the theory of the Laplacian on punctured manifolds,  $\delta \in (-1, 0)$ is in the range where $\Delta_S : C^{2,\beta}_\delta (S_0^\ast) \rightarrow  C^{0, \beta}_{\delta-2} (S_0^\ast)$ is surjective.  Let $w_{\trans} \in C^{2, \beta}_{\delta}(S_0^\ast)$ be a solution of this equation and extend it to all of $\approxsol$ by setting 
$$\bar w_{\trans} := \chi_{\trans, r_0/4, 4 r_4} w_{\trans} $$ 
and extending by symmetry.  Now compute
\begin{equation}
	\label{eqn:rebal}
	\Delta_{U, \zeta} (\bar w_{\trans}) = ( \Delta_{U, \zeta} - \Delta_S)(\bar w_{\trans}) + [ \Delta_S, \chi_{\trans, r_0/4, 4 r_4} ] (w_{\trans}) + f_{\trans} \, .
\end{equation}
The following estimates are valid.  First, 
$$ | \bar w_{\trans} |_{C^{2, \beta}_\delta} \leq C |f|_{C^{0,\beta}_{\delta - 2}} \qquad \mbox{and} \qquad | [ \Delta_S, \chi_{\trans, r_0/4, 4 r_4} ] (w_{\trans}) |_{C^{0,\beta}_{\delta - 2}} \leq C  |f|_{C^{0,\beta}_{\delta - 2}} 
$$
for a constant $C$ independent of $\eps$; and as usual, the first term in \eqref{eqn:rebal} is small, satisfying 
\begin{equation*}
	( \Delta_{U, \zeta} - \Delta_S)(\bar w_{\trans}) |_{C^{0,\beta}_{\delta - 2}} \leq C \left( \frac{\eps^{2n}}{r_0^{2n}} + r_4^2 \right) |w_{\trans} |_{C^{2, \beta}_{\delta}} \leq \frac{1}{2} |f|_{C^{0, \beta}_{\delta - 2}}
\end{equation*}
provided $\eps$ and $r_4$ are small enough and $r_0 \geq \kappa \eps$ for $\kappa$ sufficiently large.  This follows from Proposition \ref{prop:exest}.  As a result, the procedure above can be iterated to yield a solution of the equation 
$$ \Delta_{U, \zeta} (\bar w_{\trans} ) = \bar f_{\trans} $$
where $\bar f_{\trans} := f_{\trans} + E_{\trans}$ and $E_{\trans}$ has support outside of $\mathit{Pert}_\eps (S_0) \cap \bigl[ \bigcup_k B_{r_4}(p_k) \setminus B_{ r_0}(p_k) \bigr]$.  One has the estimates $|\bar w_{\trans}|_{C^{2, \beta}_{\delta}} + |E_{\trans}|_{C^{0,\beta}_{\delta-2}} \leq C |f|_{C^{0,\beta}_{\delta-2}}$ for a constant $C$ independent of $\eps$.  Note that it is necessarily the case that $\int_{\approxsol} \bar f_{\trans} = 0$.

\paragraph{Step 1. The neck regions.}  The outcome of Step 0 is that it is now only necessary to solve the equation $\Delta_{U, \zeta} (u) = \hat f$ where $\hat f := f - \bar f_{\trans}$ vanishes in $\mathit{Pert}_\eps (S_0) \cap \bigl[ \bigcup_k B_{r_4}(p_k) \setminus B_{ r_0}(p_k) \bigr]$ and satisfies $|\hat f |_{C^{0,\beta}_{\delta - 2}} \leq C | f |_{C^{0,\beta}_{\delta - 2}} $.   Note that it is necessarily the case that $\int_{\approxsol} \hat f = 0$.  

Choose $r_1< r_2 < r_3  \in (r_0, r_4)$ and set $f_{\neck} := \chi_{\neck, r_1} \hat f$.  An approximate solution of the equation $\Delta_{U, \zeta} (w) = f_{\neck}$ will now be found.  To this end, view $f_{\neck}$ as a function of compact support on the $\eps$-scaled neck $N_\eps$ and look for a solution of the equation $\Delta_N(w) = f_{\neck}$ in $C^{2,\beta}_{\delta'}(N_\eps)$ where $\delta' \in (-1,\delta)$.  The purpose is to use the compactness of the support of $f_{\neck}$ to squeeze some extra decay at the ends of the neck out of the solution.  The decomposition results used in Proposition \ref{prop:linfour} are valid and by using the symmetry of $f_{\neck}$, there is a solution of the form 
$$w_{\neck} := v_{\neck} + a \chi_{\ext, r_1} L $$
where $v_{\neck} \in C^{2,\beta}_{\delta'}(N_\eps)$ satisfies the estimate $|v_{\neck}|_{C^{2,\beta}_{\delta'}(N_\eps)} \leq C r_1^{\delta - \delta'} |f|_{C^{0,\beta}_{\delta - 2}}$ and $L(\lambda ) := \log(|\lambda|) + \tilde L(\lambda)$ belongs to the kernel of $\Delta_N$, with $\tilde L (\lambda) = O(|\lambda|^{-1})$.  The coefficient $a$ can be found explicitly: 
$$ a = \frac{1}{4 \pi} \int f_{\neck } \dif \! \mathit{Vol}_N$$
by direct computation, where $\dif \! \mathit{Vol}_N$ is the volume form of the metric on $N_\eps$.  Finally, extend this solution to all of $\approxsol$ by defining $\bar w_{\neck} := \chi_{\neck, r_2} w_{\neck}$ and extending by symmetry.  One has the estimate $|\bar w_{\neck}|_{C^{2,\beta}_{\delta}} \leq C \bigl(  r_1^{\delta - \delta'} + \eps^\delta \bigr) | f_{\neck} |_{C^{2,\beta}_{\delta - 2}}$ where $\eps^\delta$ comes from the estimate $|a| \leq C \eps^\delta |f|_{C^{0,\beta}_{\delta - 2}}$.

\paragraph{Step 2. The exterior regions.}  Set $f_{\ext} := \chi_{\ext, r_3} \hat f$.  An approximate solution of the equation $\Delta_{U,\zeta} (w) = f_{\ext}$ will now be found.   To this end, view $f_{\ext}$ as a function of compact support on $S_0 ^\ast$ and look for a solution of the equation $\Delta_S (w) = f_{\neck}$ in $C^{2,\beta}_{-\delta}(S_0^\ast)$.  The  decomposition results used in Proposition \ref{prop:linfour} are valid and by using the symmetry of $f_{\neck}$, there is a solution of the form 
$$w_{\ext} := v_{\ext} + a' \chi_{\neck, r_3} L'$$
where $v_{\neck} \in C^{2,\beta}_{-\delta}(S_0^\ast)$ satisfies the estimate $|v_{\neck}|_{C^{2,\beta}_{- \delta}(S_0)} \leq C r_3^{2 \delta} |f|_{C^{0,\beta}_{\delta - 2}}$ and $L'$ has logarithmic growth near $p_1$ and $p_2$.  Without loss of generality, one can assume that $L'(\lambda ) := \log(|\lambda|) + \tilde L'(\lambda)$ with $\tilde L' (\lambda) = O(|\lambda|^{-1})$ in the common $\lambda$-coordinate on the region of overlap between the neck and $\mathit{Pert}_{\eps} (S_0)$.  The coefficient $a'$ can again be found explicitly: 
$$ a' = - \frac{1}{4 \pi} \int f_{\ext } \dif \! \mathit{Vol}_S$$
by direct computation, where $ \mathit{Vol}_S$ is the volume form of the standard metric on $S_0$.  Finally, extend this solution to all of $\approxsol$ by defining $\bar w_{\ext} := \chi_{\ext, r_2} w_{\ext}$ and extending by symmetry.  One has the estimate $|\bar w_{\ext}|_{C^{2,\beta}_{\delta}} \leq C \bigl( r_3^{2 \delta} + \eps^\delta \bigr)| f_{\ext} |_{C^{2,\beta}_{\delta - 2}}$.

\paragraph{Step 3. Estimates and convergence.}  Set $\bar w := \bar w_{\neck} + \bar w_{\ext}$ and consider $\bar w$ as the approximate solution of the equation $\Delta_{U, \zeta} (\bar w) = \hat f$.  Before justifying this, it is straightforward to show that 
\begin{align*}
	a &= \int_{\approxsol} f_{\neck} \dif \! \mathit{Vol} + \int_{\approxsol} f_{\neck} ( \dif \! \mathit{Vol}_N - \dif \! \mathit{Vol} ) = A + O( r_0^{2+\delta}) \\[1ex]
	a' &=  - \int_{\approxsol} f_{\ext} \dif \! \mathit{Vol} - \int_{\approxsol} f_{\ext} ( \dif \! \mathit{Vol}_S - \dif \! \mathit{Vol} ) = A + O( \eps^2) 
\end{align*}  
where $A := \int_{\approxsol} f_{\neck}$ and $\int_{\approxsol} f_{\ext} = - A$.  Now perform the computation
\begin{align}
	\label{eqn:lowdimone}
	\Delta_{U, \zeta} (\bar w) - \hat f &= \big( \Delta_{U, \zeta} - \Delta_N \big)(\bar w_{\neck}) +  \big( \Delta_{U, \zeta} - \Delta_S \big)(\bar w_{\ext}) \notag \\
	&\qquad + [\chi_{\neck, r_2} , \Delta_N] (w_\neck) + [\chi_{\ext, r_2} , \Delta_S] (w_\ext) \notag \\
	&\qquad +  \chi_{\neck, r_2} f_{\neck} + \chi_{\ext, r_2} f_{\ext} - \hat f \notag \\[1ex]
	&  
	\begin{aligned}
		&= \big( \Delta_{U, \zeta} - \Delta_N \big)(\bar w_{\neck}) +  \big( \Delta_{U, \zeta} - \Delta_S \big)(\bar w_{\ext}) \\
		&\qquad + [\chi_{\neck, r_2} , \Delta_N] (v_\neck) - [\chi_{\neck, r_2} , \Delta_S] (v_\ext) \\
		&\qquad  +  [\chi_{\neck, r_2}, \Delta_N](a \tilde L) - [\chi_{\neck, r_2}, \Delta_S](  a' \tilde L') \\
		&\qquad + [\chi_{\neck, r_2}, \Delta_N - \Delta_S] (a \log(|\lambda|))  - [\chi_{\neck, r_2}, \Delta_S]( (a'-a)\log(|\lambda|))
	\end{aligned}
\end{align}
using the formul\ae\ for $a$ and $a'$ as well as the fact that $\hat f$ vanishes between radii $r_1$ and $r_3$.  The usual analysis can now be invoked to show that each term in \eqref{eqn:lowdimone} is small a small enough multiple of $|f|_{C^{0,\beta}_{\delta - 2}}$ in the $C^{0,\beta}_{\delta - 2}$ norm by suitable \emph{a priori} choice of $r_i$ and small enough $\eps$.  

The conclusion to be drawn from the work above is that iteration produces a solution of the equation $\Delta_{U, \zeta} (\bar w) = \hat f$ satisfying the estimate $|\bar w |_{C^{2,\beta}_{\delta}} \leq C\eps^\delta |f|_{C^{0,\beta}_{\delta - 2}}$ where $C$ is a constant independent of $\eps$.  Coupled with the result of Step 0, one now has a right inverse $\mathcal R_1 :   [C^{0,\beta}_{\delta-2} (\approxsol) ]_{\mathit{sym}} \rightarrow [C^{2,\beta}_{\delta} (\approxsol) ]_{\mathit{sym}}$ satisfying the bound  $| \mathcal R_1 (f) |_{C^{2,\beta}_\delta} \leq C \eps^\delta |f|_{C^{0,\beta}_{\delta - 2}}$.

It remains to construct the right inverse for $\Delta_{U, \zeta} + 6 : [C^{4,\beta}_{\gamma} (\approxsol) ]_{\mathit{sym}} \rightarrow [C^{2,\beta}_{\gamma - 2} (\approxsol) ]_{\mathit{sym}}$ with $\gamma \in (1,2)$.   This task is simpler given the range for $\gamma$ and the fact that the $L^1$ norm of a $C^{2,\beta}_{\gamma - 2}$ function in the neck region is a small multiple of the $C^{2,\beta}_{\gamma - 2}$ norm of the function provided the neck region is sufficiently small.  In fact, the procedure of Proposition \ref{prop:linthree} can be used almost verbatim, except replacing $w_{\ast}$ with $w_{\ast}^\perp$ where $\ast$ refers to either $\neck$ or $\ext$ and $\perp$ refers to the $L^2$-projection perpendicular to the constant functions.  Thus logarithmic terms will not appear in the solutions $u_{\ast}$.  The discrepancy $| w_{\ast} - w_{\ast}^\perp |_{C^{2, \beta}_{\gamma - 2}}$ is sufficiently small not to spoil the convergence of the iteration.  The result of this analysis is a right inverse $\mathcal R_2 :   [C^{2,\beta}_{\gamma-2} (\approxsol) ]_{\mathit{sym}} \rightarrow [C^{4,\beta}_{\gamma} (\approxsol) ]_{\mathit{sym}}$ satisfying the bound  $| \mathcal R_2 (w) |_{C^{4,\beta}_{\gamma }} \leq C  |w|_{C^{2,\beta}_{\gamma - 2}}$.  
\end{proof}

\section{The Non-Linear Analysis}

\subsection{The Non-Linear Estimates}
\label{subsec:nonlinest}

The remainder of the proof of the Main Theorem is devoted to establishing the final two fundamental estimates  \eqref{eqn:iftesttwo} and \eqref{eqn:iftestthree} needed to invoke the inverse function theorem.  The first of these estimates is the measurement of the size of $\Phi_{U, \zeta}(0)$ in the $\cobg$ norm.  The second of these is the measurement of the variation in $\Dif \Phi_{U, \zeta}(f)$ in the $\cobg$ norm as $f$ varies.  Choose $\gamma \in (4-n, 0)$ for $n \geq 5$ as well as $\gamma \in (4-n, 5-n)$ for $n= 4, 3$ and $\gamma \in (1, 2)$ for $n=2$.

\begin{prop}
	\label{prop:nonlinsizeest}
	The approximate solution $\approxsol$ satisfies the following estimate.  There exists some constant $C$ independent of $\eps$ so that 
	\begin{equation}
		\label{eqn:nonlinsizeest}
		\Vert \Phi_{U, \zeta} (0) \Vert_{C^{0,\beta}_{\gamma - 4}(\approxsol)} \leq  C r_\eps^{4-\gamma}  \, .
	\end{equation}
\end{prop}

\begin{proof}
 	By Proposition \ref{prop:exest} the divergence of the mean curvature in the exterior region $\Lambda'$ satisfies
	$$[\rho_\eps (x)]^{4-\gamma} | \nabla \cdot H (x)| \leq
	\begin{cases}
		C \eps^{3n} &\qquad x \in \Lambda'' \\[1ex]
		\dfrac{C \eps^{3n}}{r^{3n-2+\gamma}} &\qquad \mathrm{dist}(x, p_i) = r \: \: \mbox{and} \: \: r \geq r_\eps \, .
	\end{cases}
	$$
	By Proposition \ref{prop:neckest} the divergence of the mean curvature in the neck region $\mathcal N$ satisfies
	$$ \eps^{4-\gamma} (1+ \lambda^2)^{(4-\gamma)/2}  | \nabla \cdot H (\lambda, \mu) | \leq C \eps^{4-\gamma} \left[ \frac{1}{(1+ \lambda^2)^{(n - 4 + \gamma)/2}} + (1+ \lambda^2)^2 \right] \qquad \mbox{for} \: \: |\lambda| \leq r_\eps/\eps \, .$$
	By Proposition \ref{prop:transmeancurvest} the divergence of the mean curvature in the transition region $\mathcal T$ satisfies
	$$ r_\eps^{4-\gamma} | \nabla \cdot H | \leq C \left[ \frac{\eps^{n+1}}{r_\eps^{n-3+\gamma}}  + r_\eps^{4-\gamma} \right]$$
	The desired supremum estimate is obtained by taking the supremum above and using the ranges for $\gamma$ and the fact that $\eps < r_\eps$.  The estimate of the H\"older coefficient follows similarly.
\end{proof}

\begin{prop}
	\label{prop:nonlindifest}
	The linearization of the contact-stationary Legendrian operator near $\approxsol$ satisfies the following estimate.  If $\eps$ is sufficiently small, then there is a constant $C$ independent of $\eps$ so that 
	\begin{equation}
		\label{eqn:nonlin}
		\big\Vert \Dif \Phi_{U, \zeta}(f) (u) - \mathcal L_{U, \zeta} (u) \big\Vert_{C^{0, \beta}_{\gamma - 4}} \leq C \eps^{-2 + \gamma} \| f \|_{C^{4, \beta}_\gamma} \| u \|_{C^{4, \beta}_\gamma}
	\end{equation}
	for all $u, f \in C^{4, \beta}_\gamma(\approxsol)$.  
\end{prop}

\begin{proof}
	The desired estimate will be proved using scaling arguments as in \cite{me1,leey}.   By compactness, the estimate is certainly true in the region $\Lambda'$ of $\approxsol$.  So consider a small annular region $\Lambda_\sigma := A_\sigma \cap (\mathcal N \cup \mathcal T)$ in $\approxsol$. It is sufficient to perform the calculations for the Lagrangian projection $[\Lambda_\sigma]$ in $\cp^n$.  Consider the operator 
\begin{equation}
	\label{eqn:nonlindifest}
	\Phi_\sigma (f) := \nabla \cdot \Big( H_{\hat \phi_{E_\sigma(f)} [\Lambda_\sigma]} \elbow \omega_0 \Big)
\end{equation}
	where $E_\sigma(f)$ is the extension of a function $f: \Lambda_\sigma \rightarrow \R$ to a tubular neighbourhood of $[\Lambda_\sigma]$ and $\hat \phi_{E_\sigma(f)}$ is the associated time-one Hamiltonian flow.
		
	The next step is to determine how all the objects in \eqref{eqn:nonlindifest} scale with $\sigma$.  First, $[\Lambda_\sigma] = \sigma [\Lambda_1]$ where $\Lambda_1 = A_1 \cap  \sigma^{-1} (\mathcal N \cup \mathcal T)$.  Moreover, $[\Lambda_1]$ is some Lagrangian submanifold whose geometry is bounded by a universal constant since the norm of the second fundamental form of $[\Lambda_\sigma]$ in the annulus $A_\sigma$ is $O(\sigma^{-1})$.  Suppose that $[\Lambda_1]$ carries the metric $g_1$ so that $[\Lambda_\sigma]$ carries the metric $g_\sigma = \sigma^2 g_1$.  Since $E_\sigma (f) = f \circ \exp^{-1}_{\Lambda_\sigma}$ near $\Lambda_\sigma$,  then one can check that $\hat \phi_{E_\sigma(f)} ([\Lambda_\sigma]) = \sigma \hat \phi_{E_1(f / \sigma^2)}([\Lambda_1])$.  Consequently,
\begin{equation}
	\label{eqn:nonlindifesttwo}
	\Phi_\sigma (f) = \nabla \cdot \Big( H \big(\sigma \hat \phi_{E_1(f / \sigma^2)}([\Lambda_1]) \big) \elbow \omega_0 \Big) = \frac{1}{\sigma^2} \Phi_1\left( \frac{f}{\sigma^2} \right)
\end{equation}
using the scaling property of the mean curvature and the covariant derivative under conformal transformation.  All quantities on the right hand side of \eqref{eqn:nonlindifesttwo} refer to the metric $g_1$.

One can now derive the desired estimate.  Equation \eqref{eqn:nonlindifesttwo} implies that the linearization of $\Phi_\sigma$ in $A_\sigma$ must satisfy
\begin{align}
	\label{eqn:nonlindifestthree}
	\big\vert \Dif \Phi_\sigma (f) (u) - \Dif \Phi_\sigma (0) (u) \big\vert &=  \frac{1}{\sigma^2} \left\vert \Dif \Phi_1\left( \frac{f}{\sigma^2} \right) \left( \frac{u}{\sigma^2} \right) -  \Dif \Phi_1(0) \left( \frac{u}{\sigma^2} \right) \right\vert \notag \\
	&\leq \frac{C}{\sigma^6} \| f \|_{C^4(\Lambda_1)} \, \| u \|_{C^4(\Lambda_1)}
\end{align}
where $C$ is some universal constant pertaining to $\Phi_1$ on $\Lambda_1$.  Multiplying by $\sigma^{4 - \gamma}$ and reversing the scaling in equation \eqref{eqn:nonlindifestthree} gives
\begin{equation}
	\label{eqn:nonlindifestfour}
	\big\Vert \Dif \Phi_\sigma (f) (u) - \Dif \Phi_\sigma (0) (u) \big\Vert_{C^4_{\gamma - 4}(\Lambda_\sigma)} \leq C \sigma^{-2 +\gamma} \|f\|_{C^4_{\gamma}(\Lambda_\sigma)} \, \|u\|_{C^4_\gamma(\Lambda_\sigma)} \, ,
\end{equation}
where all quantities in \eqref{eqn:nonlindifestfour} refer to the metric $g_\sigma$.  One can now piece the above estimates together for different $A_\sigma$, using the fact that the smallest $\sigma $ can be is $O(\eps)$ in the centre of the neck region, and obtain the desired supremum estimate.  The estimate of the H\"older coefficient follows similarly.
\end{proof}

\subsection{The Proof of the Main Theorem}

The three fundamental estimates needed to invoke the inverse function theorem for $\Phi_{U, \zeta}$ when $k$ is sufficiently large and $\eps := \eps_k$ is sufficiently small have now been established and the proof of the Main Theorem is at hand.  The theorem is re-formulated here in the more technical language of the preceding sections.  Its proof is a re-organization and summary of all the results above.

\begin{thm}
	\label{thm:mainthm}
	Let $U \in SU(n+1)$ be as in Example 1 or Example 2 and let $\zeta = \me^{2 \pi \mi / k}$ be such that it is possible to construct an approximate solution $\approxsol$ as in Section \ref{sec:approxsol}.  Let $X, Z$ be the Banach subspace of functions on $\approxsol$ defined in Definition \ref{def:bspace}.  Choose $\gamma \in (4-n, 0)$ for $n \geq 5$ as well as $\gamma \in (4-n, 5-n)$ for $n= 4, 3$ and $\gamma \in (1,2)$ for $n =2$.   Set $r_\eps := \eps^s $.  Then there is $s \in (0,1)$ so that a solution of the contact-stationary Legendrian equation $\Phi_{U, \zeta}(f) = 0$ can be found with $f \in X$ satisfying the bound $|f|_{\cfbg} \leq C r_\eps^{4-\gamma}$.  
\end{thm}

\begin{proof}
To prove this theorem by means of the Banach space inverse function theorem and the analysis contained within this paper, it is necessary to establish the following facts.  First, one must show that $\limsup_\eps C(\eps) \| Q_{U, \zeta} \|_{\mathit{op}} < 1/2$ where $\| Q_{U, \zeta} \|_{\mathit{op}}$ is the operator norm of $Q_{U, \zeta}$ calculated in Lemma \ref{lemma:loworder} and $C(\eps)$ is the upper bound for the right inverse of $\mathcal L_{U, \zeta}$ found in Section \ref{sec:linest}.  Since $| \Phi_{U, \zeta}(0) |_{\cobg} = O (  r_\eps^{4-\gamma} )$ found in Proposition \ref{prop:nonlinsizeest} and $R(\eps) = O(\eps^{2-\gamma})$ found in Proposition \ref{prop:nonlindifest}, one must then show that $r_\eps^{4-\gamma} \eps^{\gamma - 2} C(\eps)$ can be made as small as desired, where  $C(\eps)$ is the upper bound for the right inverse of $\mathcal L_{U, \zeta}$ found in Section \ref{sec:linest}.
\begin{itemize}
	\item When $n \geq 3$ then $C(\eps) = O(1)$ and 
	$$\lim_{\eps \rightarrow 0} \| Q_{U, \zeta} \| \leq C \cdot  \lim_{\eps \rightarrow 0} \bigl( \eps^{n - s(n-2)} + \eps^{2} + \eps^{4s} \bigr) =  0 $$
	provided $s < n/(n-2)$.  Also, $r_\eps^{4-\gamma} \eps^{\gamma - 2} = \eps^{(4-\gamma)s + \gamma - 2} \rightarrow 0$ when $s > (2-\gamma) / (4-\gamma)$.  By choice of $\gamma$ and $n$, the range of such $s \in (0,1)$ is non-empty and thus the theorem is true.
	
	\item When $n=2$ then $C(\eps) = O(\eps^{\gamma - 2})$ and
	$$\lim_{\eps \rightarrow 0} C(\eps) \| Q_{U, \zeta} \| \leq C \cdot \lim_{\eps \rightarrow 0} \bigl(  \eps^{2} + \eps^{4s} \bigr) \eps^{\gamma - 2}  =  0 $$
	provided $s > (2-\gamma)/4$.  Also, $r_\eps^{4-\gamma} \eps^{\gamma - 2} = \eps^{(4-\gamma)s + \gamma - 2} \rightarrow 0$ when $s > (4-2\gamma) / (4-\gamma)$.  By choice of $\gamma$ the range of such $s \in (0,1)$ is non-empty and thus the theorem is true.

\end{itemize}
This completes the proof of the theorem.
\end{proof}

\subsection{Embeddedness of the Solutions}

Denote the solution of the contact-stationary Legendrian problem constructed from $U \in \sun$ with $\zeta = \me^{2 \pi \mi / k} $ in the previous sections by $f_k$.  Let $\phi_{f_k}$ be the contact deformation constructed from $f_k$ using the method developed in Section \ref{subsec:def}.   A simple chain of reasoning shows that the deformed submanifold $\phi_{f_k} ( \approxsol)$ is embedded whenever $\approxsol$ itself is, so long as $k$ is sufficiently large.  

Suppose that $\approxsol$ is embedded in $\Sph^{2n+1}$.  Then $\approxsol$ is contained in some non-self-intersecting tubular neighbourhood of itself.  The width of this tubular neighbourhood is clearly larger in the parts of $\approxsol$ that are subsets of some $U^s( S_0)$ and smaller in the neck regions of $\approxsol$.  In fact, one can argue based on scaling that the width of the tubular neighbourhood at a point $p \in \approxsol$ is $O(\rho_\eps(p))$ where $\rho_\eps$ is the weight function from Definition \ref{defn:weightfn} that is used to define the $\clbg$ norm.  The question of embeddedness can now be re-phrased in terms of this tubular neighbourhood: $\phi_{f_k}(\approxsol)$ fails to be embedded if either $\phi_{f_k}(\approxsol)$ has local self-intersection somewhere within the tubular neighbourhood, or else $\phi_{f_k}(\approxsol)$ intersects itself by leaving the tubular neighbourhood somewhere and re-entering it somewhere else.

In order to decide if $\phi_{f_k}(\approxsol)$ intersects itself in one of these two ways, one must understand how `far' the contact deformation $\phi_{f_k}$ can move the points of $\approxsol$.  Recall that $\phi_{f_k}$ is the time-one flow of the contact vector field corresponding to the function $E(f_k)$ that is an extension of $f_k$ orthogonal to $\approxsol$.  Thus the distance of $\phi_{f_k}(p)$ from $p$ is governed by the size of this vector field, which in turn is governed by the size of $f_k$ (in the Hopf direction) and the first derivative of $f_k$ (in contact directions) so long as these quantities are sufficiently small.  Moreover, $\phi_{f_k}( \approxsol)$ remains graphical over $\approxsol$ so long as the derivative of $\phi_{f_k}$ remains sufficiently small, which in turn requires that $f_k$ is small up to its second derivative.  In fact, it is necessary to have $ |f_k(x)| + |\nabla f_k (x)| = O(\rho_\eps(x))$ and $|\nabla^2 f_k(x)| \leq O(1)$ for all $x \in \approxsol$.   It is now a simple matter to verify that the estimates of the size of $f_k$ and its derivatives from Theorem \ref{thm:mainthm} imply that the requirements for embeddedness are met according to the remark following Definition \ref{defn:contacto}.  In summary, the above chain of reasoning leads to the following result.

\begin{prop}
	Let $f_k : \approxsol \rightarrow \R$ be the solution of the equation $\Phi_{U,\zeta} (f_k) = 0$, where $\zeta = \me^{2 \pi \mi / k}$, that was constructed in Theorem \ref{thm:mainthm}.  If $\approxsol$ is an embedded submanifold of $\Sph^{2n+1}$ then the deformed submanifold $\phi_{f_k} (\approxsol)$ is embedded for all sufficiently large $k$.
\end{prop}

\renewcommand{\baselinestretch}{1.0}
\small
\bibliography{cminlag}
\bibliographystyle{amsplain}

\end{document}